\documentclass[12pt,twoside,leqno,openany]{amsart}

\usepackage{amssymb,amsbsy,amsmath,amsfonts,amssymb,amscd,times,
graphics,color,xypic,footmisc,fancyhdr,multicol,fancybox,
graphicx,mathrsfs,rotating,ifthen,wasysym}
\usepackage[all]{xy}

\usepackage[T1]{fontenc}
\newcommand{\C}{\mathbb{C}}

\newcommand{\R}{\mathbb{R}}

\newcommand{\red}{\textcolor{red}}

\newcommand{\HEAD}[2]{%
\pagestyle{fancy}
\fancyhead[RO]{\tiny\sf\thepage}
\fancyhead[CO]{{\tiny\sf #1}}
\fancyhead[LE]{\tiny\sf\thepage}
\fancyhead[CE]{{\tiny\sf #2}}
\fancyfoot{}}

\theoremstyle{definition}

\renewcommand{\det}{\text{\footnotesize\sf det}}

\newcommand{\explain}[1]{\text{\scriptsize\sf [#1]}}

\newcommand{\function}{\text{\footnotesize\sf function}}

\newcommand{\isqrt}{{\scriptstyle{\sqrt{-1}}}}

\renewcommand{\lim}{\text{\footnotesize\sf lim}}

\renewcommand{\mod}{\text{\footnotesize\sf mod}}

\newcommand{\rank}{\text{\footnotesize\sf rank}}

\newcommand{\something}{\text{\scriptsize\sf something}}

\newcommand{\vf}{\vfill\end{document}}

\newcommand{\zero}[1]{\underline{#1}_{\red{\circ}}}

\let\mathcal\mathscr

\begin{document}

$\:$

\bigskip\bigskip

\begin{center}

{\Large\bf Equivalences of $5$-dimensional CR-manifolds V:}

\medskip

{\Large\bf Six initial frames and coframes;}

\medskip

{\Large\bf Explicitness obstacles}

\end{center}

\medskip

\begin{center}
Jo\"el {\sc Merker}
\end{center}

\bigskip

\begin{center}
\begin{minipage}[t]{10.25cm}
\baselineskip =0.32cm 
{\scriptsize
{\bf Abstract.}
Local CR-generic submanifolds of $\C^N$ are in one-to-one
correspondence with their respective graphing functions, but it is
well known that (despite their importance) the
Cartan-Hachtroudi-Chern-Moser invariants and coframes for Levi
nondegenerate hypersurfaces $M \subset \C^{ n+1}$ have been fully
computed in CR dimension $n \geqslant 2$ only in special 
cases which show off a tremendous collapse of
computational complexity in comparison to the general case.  One of
the goals of this Part V is to set up systematic initial data that are
essentially explicit in terms of the concerned graphing functions, for
the six already studied general classes
$\text{\sf I}$,
$\text{\sf II}$,
$\text{\sf III}_{\text{\sf 1}}$,
$\text{\sf III}_{\text{\sf 2}}$,
$\text{\sf IV}_{\text{\sf 1}}$,
$\text{\sf IV}_{\text{\sf 2}}$.
Incredibly, for Class $\text{\sf III}_{\text{\sf 1}}$ CR-generic
submanifolds $M^5 \subset \C^4$ that are the geometry-preserving
deformations of one of the natural models of Beloshapka, even the
initial frame and coframe are not absorbable by an individual personal
computer, for some of the concerned coefficient-functions incorporate
nearly $100\,000\,000$ of monomials in $165$ jet variables, not to
mention that the exploration of biholomorphic equivalences yet
requires to differentiate such functions at least four times.  As will
appear later on, deep (archaic) mathematical links are extant between
the effective Cartan theory and the famous hyperbolicity conjecture of
Kobayashi. 
}
\end{minipage}
\end{center}

\bigskip

\begin{center}
\begin{minipage}[t]{11.75cm}
\baselineskip =0.35cm {\scriptsize

\centerline{\bf Table of contents}

\medskip

{\bf \ref{initial-I}.~$M^3 \subset \C^2$ of general class $\text{\sf I}$: 
initial frame and coframe in local coordinates
\dotfill~\pageref{initial-I}.}

{\bf \ref{initial-II}.~$M^4 \subset \C^3$ of general class $\text{\sf II}$: 
initial frame and coframe in local coordinates
\dotfill~\pageref{initial-II}.}

{\bf \ref{initial-III-1}.~$M^5 \subset \C^4$ of 
general class $\text{\sf III}_{\text{\sf 1}}$: 
initial frame and coframe in local coordinates
\dotfill~\pageref{initial-III-1}.}

{\bf \ref{initial-III-2}.~$M^5 \subset \C^4$ of 
general class $\text{\sf III}_{\text{\sf 2}}$: 
initial frame and coframe in local coordinates
\dotfill~\pageref{initial-III-2}.}

{\bf \ref{initial-IV-1}.~$M^5 \subset \C^3$ of 
general class $\text{\sf IV}_{\text{\sf 1}}$: 
initial frame and coframe in local coordinates
\dotfill~\pageref{initial-IV-1}.}

{\bf \ref{initial-IV-2}.~$M^5 \subset \C^3$ of 
general class $\text{\sf IV}_{\text{\sf 2}}$: 
initial frame and coframe in local coordinates
\dotfill~\pageref{initial-IV-2}.}

}\end{minipage}
\end{center}

\bigskip


\bigskip

\section{\sf $M^3 \subset \C^2$ of general class $\text{\sf I}$: 
\\
initial frame and coframe in local coordinates}
\label{initial-I}
\HEAD{\ref{initial-I}.~$M^3 \subset \C^2$: 
Initial frame and coframe in local coordinates}{
Jo\"el {\sc Merker}, D\'epartement de Math\'ematiques d'Orsay}

\medskip

Consider the {\em most simple case} of:
\[
\Big(
M^3
\,\subset\,
\C^2
\Big)
\,\,\in\,\,
\text{\sf General Class $\text{\sf I}$}.
\]
Representing as before
(\cite{ Merker-Pocchiola-Sabzevari-5-CR-II, Merker-5-CR-III})
$M$ in coordinates:
\[
(z,w) 
= 
\big(
x+\isqrt\,y,\, 
u+\isqrt\,v\big),
\]
as a graph:
\[
v
=
\varphi(x,y,u),
\]
an associated explicit frame for $TM$ is:
\[
\aligned
X
&
=
\frac{\partial}{\partial x}
+
\varphi_x(x,y,u)\,\frac{\partial}{\partial v},
\\
Y
&
=
\frac{\partial}{\partial y}
+
\varphi_y(x,y,u)\,\frac{\partial}{\partial v},
\\
U
&
=
\frac{\partial}{\partial u}
+
\varphi_u(x,y,u)\,\frac{\partial}{\partial v},
\endaligned
\]
viewed here {\em extrinsically}, namely together with the (transversal) 
coordinate $v$ which is not internal to $M$. 

A way of
understanding these fields intrinsically is to introduce the auxiliary
projection of $M$ onto its tangent plane at $0$:
\[
\pi\big\vert_M
\colon\ \ \
M
\,\longrightarrow\,
T_0M
\]
defined as the restriction to $M$ of the coordinate projection:
\[
\aligned
\pi
\colon\ \ \ \ \ \ \ \ \ \ \ \ \ \ \
\R^4
&
\,\longrightarrow\,
\R^3
\\
(x,y,u,v)
&
\,\longmapsto\,
(x,y,u).
\endaligned
\]

Considering then $\pi \vert_M \colon M
\longrightarrow \R^3$ as a {\sl chart} on $M$, 
one recovers the intrinsic representation of the tangential fields
simply as:
\[
\aligned
\pi_*
\bigg(
\frac{\partial}{\partial x}
+
\varphi_x\,
\frac{\partial}{\partial v}
\bigg)
&
=
\frac{\partial}{\partial x},
\\
\pi_*
\bigg(
\frac{\partial}{\partial y}
+
\varphi_y\,
\frac{\partial}{\partial v}
\bigg)
&
=
\frac{\partial}{\partial y},
\\
\pi_*
\bigg(
\frac{\partial}{\partial u}
+
\varphi_u\,
\frac{\partial}{\partial v}
\bigg)
&
=
\frac{\partial}{\partial u},
\endaligned
\]
and this visibly corresponds just to dropping the $\frac{ \partial}{
\partial v}$-components.

\smallskip
 
Next, as was shown in~\cite{Merker-Pocchiola-Sabzevari-5-CR-I}, 
a natural local vector field generator
of $T^{1, 0} M$ is:
\[ 
\mathcal{L} 
:= 
\frac{\partial}{\partial z} + {\tt 
A}\,\frac{\partial}{\partial w}, 
\] 
with:  
\[ 
{\tt A} 
= 
\frac{-\,2\,\varphi_z}{\isqrt+\,\varphi_u}, 
\] 
a coefficient-function
which is thus {\em de facto} a function of only $(x, y, u)$,
independently of $v$.

Restricting $\mathcal L$ to $M^3$, one must simply and
only drop the (extrinsic) vector field $\frac{ \partial}{ \partial
v}$: 
\[ 
\aligned 
\mathcal{L}\big\vert_M 
=
\pi_*\big(\mathcal{L}\big)
& 
= 
\frac{\partial}{\partial z} 
+ 
{\tt A}\, 
\bigg(
\frac{1}{2}\,\frac{\partial}{\partial u} 
- 
\zero{\frac{\isqrt}{2}\,
\frac{\partial}{\partial v}} 
\bigg) 
\\ 
& 
= 
\frac{\partial}{\partial z} 
- 
\frac{\,\varphi_z}{\isqrt+\varphi_u}\, 
\frac{\partial}{\partial u}. 
\endaligned 
\] 
One uses the notation:
\[
\aligned
A
:=
&\,
\frac{\tt A}{2}
\\
=
&\,
-\,\frac{\varphi_z}{\isqrt+\varphi_u}.
\endaligned
\]

Thus intrinsically on $M^3$, the CR-structure induced by the ambient 
$\C^2$ on $M^3$ is encoded by the 
{\em intrinsic} $(1,0)$ complex-valued vector field:
\[ 
\boxed{\,
\mathcal{L}
= 
\frac{\partial}{\partial z} 
-
\frac{\varphi_z}{\isqrt+\varphi_u}\,
\frac{\partial}{\partial u},\,}
\]
together with its conjugate:
\[ 
\boxed{\,
\overline{\mathcal{L}} 
= 
\frac{\partial}{\partial\overline{z}} 
-
\frac{\varphi_{\overline{z}}}{-\isqrt+\varphi_u}\,
\frac{\partial}{\partial u}.\,}
\]

By the assumption that $M$ belongs to the
General Class $\text{\sf 1}$:
\[
\big\{
\mathcal{L},\,\overline{\mathcal{L}},\,
\big[\mathcal{L},\overline{\mathcal{L}}\big]
\big\}
\]
constitutes a frame for $\C \otimes_\R TM$.

\medskip

Set:
\[
\mathcal{T}
:=
\isqrt\,
\big[\mathcal{L},\overline{\mathcal{L}}\big],
\]
and compute:
\[
\aligned
\mathcal{T}
&
:=
\isqrt\,
\bigg[
\frac{\partial}{\partial z} 
-
\frac{\varphi_z}{\isqrt+\varphi_u}\,
\frac{\partial}{\partial u},\,\,
\frac{\partial}{\partial\overline{z}} 
-
\frac{\varphi_{\overline{z}}}{-\isqrt+\varphi_u}\,
\frac{\partial}{\partial u}
\bigg]
\\
&
=
\frac{1}{(\isqrt+\varphi_u)^2\,(-\isqrt+\varphi_u)^2}\,
\bigg\{
2\,\varphi_{z\overline{z}}
+
2\,\varphi_{z\overline{z}}\,\varphi_u\,\varphi_u
-
2\,\isqrt\,\varphi_{\overline{z}}\,\varphi_{zu}
-
\\
&
\ \ \ \ \ \ \ \ \ \ \ \ \ \ \ \ \ \ \ \ \ \ \ \ \ \ \ \ \ \ \ \ \ \
\ \ \ \ \ \ \ \ \ \ \ \ \ \ \ \ \ \ \
-\,
2\,\varphi_{\overline{z}}\,\varphi_{zu}\,\varphi_u
+
2\,\isqrt\,\varphi_z\,\varphi_{\overline{z}u}
+
\\
&
\ \ \ \ \ \ \ \ \ \ \ \ \ \ \ \ \ \ \ \ \ \ \ \ \ \ \ \ \ \ \ \ \ \
\ \ \ \ \ \ \ \ \ \ \ \ \ \ \ \ \ \ \
+
2\,\varphi_z\,\varphi_{\overline{z}}\,\varphi_{uu}
-
2\,\varphi_z\,\varphi_{\overline{z}u}\,\varphi_u
\bigg\}\,
\frac{\partial}{\partial u}.
\endaligned
\]

Abbreviate the appearing coefficient-function as:
\[
\aligned
\ell
&
:=
\frac{1}{(\isqrt+\varphi_u)^2\,(-\isqrt+\varphi_u)^2}\,
\bigg\{
2\,\varphi_{z\overline{z}}
+
2\,\varphi_{z\overline{z}}\,\varphi_u\,\varphi_u
-
2\,\isqrt\,\varphi_{\overline{z}}\,\varphi_{zu}
-
\\
&
\ \ \ \ \ \ \ \ \ \ \ \ \ \ \ \ \ \ \ \ \ \ \ \ \ \ \ \ \ \ \ \ \ \
\ \ \ \ \ \ \ \ \ \ \ \ \ \ \ \ \ \ \
-\,
2\,\varphi_{\overline{z}}\,\varphi_{zu}\,\varphi_u
+
2\,\isqrt\,\varphi_z\,\varphi_{\overline{z}u}
+
\\
&
\ \ \ \ \ \ \ \ \ \ \ \ \ \ \ \ \ \ \ \ \ \ \ \ \ \ \ \ \ \ \ \ \ \
\ \ \ \ \ \ \ \ \ \ \ \ \ \ \ \ \ \ \
+
2\,\varphi_z\,\varphi_{\overline{z}}\,\varphi_{uu}
-
2\,\varphi_z\,\varphi_{\overline{z}u}\,\varphi_u
\bigg\}.
\endaligned
\]

To know the full Lie bracket structure of the frame:
\[
\big\{
\mathcal{L},\,\overline{\mathcal{L}},\,
\big[\mathcal{L},\overline{\mathcal{L}}\big]
\big\},
\]
it yet remains to compute:
\[
\aligned
&
\big[\mathcal{L},\mathcal{T}\big],
\\
&
\big[\overline{\mathcal{L}},\mathcal{T}\big],
\endaligned
\]
the second being the conjugate of the first, for, 
as always, $\mathcal{ T}$ is real:
\[
\overline{\mathcal{T}}
=
\overline{
\isqrt\,\big[\mathcal{L},\,\overline{\mathcal{L}}\big]
}
=
-\,\isqrt\,\big[\overline{\mathcal{L}},\,\mathcal{L}\big]
=
\isqrt\,\big[\mathcal{L},\,\overline{\mathcal{L}}\big]
=
\mathcal{T}.
\]

Since:
\[
\mathcal{T}
=
\ell\cdot
\frac{\partial}{\partial u},
\]
one has:
\[
\aligned
\big[\mathcal{L},\mathcal{T}\big]
&
=
\bigg[
\frac{\partial}{\partial z} 
-
\frac{\varphi_z}{\isqrt+\varphi_u}\,
\frac{\partial}{\partial u},\,\,
\ell\,
\frac{\partial}{\partial u}
\bigg]
\\
&
=
\bigg\{
\mathcal{L}(\ell)
+
\mathcal{T}
\bigg(
\frac{\varphi_z}{\isqrt+\varphi_u}\,
\bigg)
\bigg\}\,
\frac{\partial}{\partial u},
\endaligned
\]
hence replacing:
\[
\frac{\partial}{\partial u}
=
\frac{1}{\ell}\,
\mathcal{T},
\]
one gets:
\[
\big[\mathcal{L},\mathcal{T}\big]
=
\frac{
\mathcal{L}(\ell)
+
\mathcal{T}
\big(
\frac{\varphi_z}{\isqrt+\varphi_u}\,
\big)}{\ell}
\cdot
\mathcal{T}
\]

Call this coefficient-function:
\[
\boxed{\,
P
:=
\frac{
\mathcal{L}(\ell)
+
\mathcal{T}
\big(
\frac{\varphi_z}{\isqrt+\varphi_u}\,
\big)}{\ell}.\,}
\]

\medskip\noindent{\bf Lemma.}
{\em On a hypersurface:}
\[
\Big(
M^3
\,\subset\,
\C^2
\Big)
\,\,\in\,\,
\text{\sf General Class $\text{\sf I}$},
\]
{\em graphed as:}
\[
v
=
\varphi(x,y,u),
\]
{\em with natural generator for $T^{1, 0}M$:}
\[
\mathcal{L}
= 
\frac{\partial}{\partial z} 
-
\frac{\varphi_z}{\isqrt+\varphi_u}\,
\frac{\partial}{\partial u}
\]
{\sl the associated frame for $\C \otimes_\R TM$:}
\[
\big\{
\mathcal{L},\,\overline{\mathcal{L}},\,
\isqrt\,\big[\mathcal{L},\overline{\mathcal{L}}\big]
\big\}
=
\big\{
\mathcal{L},\,\overline{\mathcal{L}},\,\mathcal{T}
\big\}
\]
{\em enjoys the Lie structure:}
\[
\aligned
\big[\mathcal{L},\overline{\mathcal{L}}\big]
&
=
-\,\isqrt\,\mathcal{T},
\\
\big[\mathcal{L},\mathcal{T}\big]
&
=
P\cdot\mathcal{T},
\\
\big[\overline{\mathcal{L}},\mathcal{T}\big]
&
=
\overline{P}\cdot\mathcal{T},
\endaligned
\]
where:
\[
\aligned
P
=
\frac{P^{\sf \,numerator}}{
P^{\sf \,denominator}},
\endaligned
\]
where:
\[
\aligned
P^{\sf \,denominator}
&
=
\big(1+\varphi_u\varphi_u\big)
\big(\isqrt+\varphi_u\big)
\Big(
\varphi_{z\overline{z}}
+
\varphi_{z\overline{z}}\varphi_u\varphi_u
-
\isqrt\,\varphi_{\overline{z}}\varphi_{zu}
-
\\
&
\ \ \ \ \
-
\varphi_{\overline{z}}\varphi_{zu}\varphi_u
+
\isqrt\,\varphi_z\varphi_{\overline{z}u}
-\varphi_z\varphi_{\overline{z}u}\varphi_u
+
\varphi_z\varphi_{\overline{z}}\varphi_{uu}
\Big),
\endaligned
\]
and where:
\[
\aligned
P^{\sf \,numerator}
&
=
\isqrt\,
\varphi_{zz\overline{z}}
-
\varphi_{zz}\varphi_{\overline{z}u}
+
\varphi_{z\overline{z}}\varphi_{zu}
+
\isqrt\,\varphi_z\varphi_{z\overline{z}}\varphi_{uu}
+
2\,\isqrt\,\varphi_z\varphi_{\overline{z}}
\varphi_{zuu}
+
\\
&
\ \ \ \ \
+
\varphi_{\overline{z}}\varphi_{zzu}
+
\isqrt\,\varphi_{zz}\varphi_{\overline{z}}\varphi_{uu}
-
\varphi_z\varphi_z\varphi_{\overline{z}}\varphi_{uuu}
-
2\,\varphi_z\varphi_{z\overline{z}u}
-
\\
&
\ \ \ \ \
-
\isqrt\,\varphi_{\overline{z}}\varphi_{zu}\varphi_{zu}
-
\isqrt\,\varphi_z\varphi_z\varphi_{\overline{z}uu}
-
\isqrt\,\varphi_z\varphi_{\overline{z}u}\varphi_{zu}
+
\endaligned
\]
\[
\aligned
&
+
\varphi_u\Big(
-\,2\,\isqrt\,\varphi_{\overline{z}}\varphi_{zzu}
-
4\,\isqrt\,\varphi_{z\overline{z}}\varphi_{zu}
+
4\,\varphi_{uu}\varphi_{uu}\varphi_z\varphi_z\varphi_{\overline{z}}
+
3\,\varphi_z\varphi_{z\overline{z}}\varphi_{uu}
+
2\,\varphi_z\varphi_{\overline{z}}\varphi_{zuu}
-
\\
&
\ \ \ \ \ \ \ \ \ \ \
-\,
8\,\isqrt\,\varphi_{zu}\varphi_z\varphi_{\overline{z}}\varphi_{uu}
+
4\,\isqrt\,\varphi_{uu}\varphi_z\varphi_z\varphi_{\overline{z}u}
+
\varphi_{zz\overline{z}}
+
\varphi_z\varphi_z\varphi_{\overline{z}uu}
-
5\,\varphi_{\overline{z}}\varphi_{zu}\varphi_{zu}
+
\\
&
\ \ \ \ \ \ \ \ \ \ \
+
5\,\varphi_z\varphi_{\overline{z}u}\varphi_{zu}
+
\varphi_{zz}\varphi_{\overline{z}}\varphi_{uu}
\Big)
+
\endaligned
\]
\[
\aligned
&
+
\varphi_u^2\Big(
-\,\isqrt\,\varphi_z\varphi_{\overline{z}u}\varphi_{zu}
-
\isqrt\,\varphi_z\varphi_z\varphi_{\overline{z}uu}
-
8\,\varphi_{zu}\varphi_z\varphi_{\overline{z}}\varphi_{uu}
-
4\,\varphi_z\varphi_{z\overline{z}u}
-
\varphi_z\varphi_z\varphi_{\overline{z}}\varphi_{uuu}
-
\\
&
\ \ \ \ \ \ \ \ \ \ \
-\,
4\,\varphi_{uu}\varphi_z\varphi_z\varphi_{\overline{z}u}
-
2\,\varphi_{zz}\varphi_{\overline{z}u}
+
\isqrt\,\varphi_{zz}\varphi_{\overline{z}}\varphi_{uu}
+
7\,\isqrt\,\varphi_{\overline{z}}\varphi_{zu}\varphi_{zu}
+
\\
&
\ \ \ \ \ \ \ \ \ \ \
+
\isqrt\,\varphi_z\varphi_{z\overline{z}}\varphi_{uu}
+
2\,\isqrt\,\varphi_{zz\overline{z}}
-
2\,\varphi_{z\overline{z}}\varphi_{zu}
+
2\,\isqrt\,\varphi_z\varphi_{\overline{z}}\varphi_{zuu}
\Big)
+
\endaligned
\]
\[
\aligned
&
+
\varphi_u^3\Big(
-\,4\,\isqrt\,\varphi_{z\overline{z}}\varphi_{zu}
+
\varphi_z\varphi_z\varphi_{\overline{z}uu}
-
2\,\isqrt\,\varphi_{\overline{z}}\varphi_{zzu}
+
3\,\varphi_{\overline{z}}\varphi_{zu}\varphi_{zu}
+
2\,\varphi_{zz\overline{z}}
+
\\
&
\ \ \ \ \ \ \ \ \ \ \
+
5\,\varphi_z\varphi_{\overline{z}u}\varphi_{zu}
+
\varphi_{zz}\varphi_{\overline{z}}\varphi_{uu}
+
2\,\varphi_z\varphi_{\overline{z}}\varphi_{zuu}
+
3\,\varphi_z\varphi_{z\overline{z}}\varphi_{uu}
\Big)
+
\endaligned
\]
\[
\aligned
&
+
\varphi_u^4\Big(
-\,2\,\varphi_z\varphi_{z\overline{z}u}
-
3\,\varphi_{z\overline{z}}\varphi_{zu}
-
\varphi_{zz}\varphi_{\overline{z}u}
-
\varphi_{\overline{z}}\varphi_{zzu}
+
\isqrt\,\varphi_{zz\overline{z}}
\Big)
+
\ \ \ \ \ \ \ \ \ \ \ \ \ \ \ \ \ \ \ \ \ \ \ \ \ \ \ \ \ \ \ \ \ 
\\
&
+
\varphi_u^5\Big(
\varphi_{zz\overline{z}}
\Big).
\endaligned
\]

\proof
Just a direct computation.
\endproof

Strikingly (but as is known), the so-called
{\sl rigid case} where the graphing function:
\[
v
=
\varphi(x,y)
\]
is assumed\,\,---\,\,as a simplifying assumption\,\,---\,\,to
be {\em independent} of the CR-transversal coordinate $u$, 
comes up with a {\em tremendously spectacular 
collapse of complexity}.

Indeed, then:
\[
\aligned
\mathcal{L}
&
=
\frac{\partial}{\partial z}
+
\isqrt\,\varphi_z\,\frac{\partial}{\partial u},
\\
\overline{\mathcal{L}}
&
=
\frac{\partial}{\partial\overline{z}}
-
\isqrt\,\varphi_{\overline{z}}\,
\frac{\partial}{\partial u},
\endaligned
\]
hence:
\[
\mathcal{T}
\,=\,
2\,\varphi_{z\overline{z}}\,
\frac{\partial}{\partial u},
\]
which inverts as:
\[
\frac{\partial}{\partial u}
\,=\,
\frac{1}{2\,\varphi_{z\overline{z}}}\,
\mathcal{T},
\]
so that:
\[
\aligned
\big[\mathcal{L},\mathcal{T}\big]
&
=
2\,\varphi_{zz\overline{z}}\,\frac{\partial}{\partial u}
\\
&
=
\frac{\varphi_{zz\overline{z}}}{\varphi_{z\overline{z}}}\,
\mathcal{T},
\endaligned
\]
{\em i.e.}:
\[
\boxed{\,
P
\,=\,
\frac{\varphi_{zz\overline{z}}}{\varphi_{z\overline{z}}},
\,}
\]
which is a {\em considerable} collapse, indeed!

\medskip\noindent{\bf Dual coframe and its Darboux structure.}
Now, introduce three differential
$1$-forms, sections of $\C\otimes_\R T^*M$:
\[
\big\{
\rho_0,\,\overline{\zeta}_0,\,\zeta_0
\big\}
\]
that are dual to the frame:
\[
\big\{
\mathcal{T},\,\mathcal{L},\,\overline{\mathcal{L}}
\big\}
\]
namely satisfy:
\[
\aligned
\rho_0\big(\mathcal{T}\big)
&
=
1, 
\ \ \ \ \ \ \ \ \ \
\rho_0\big(\overline{\mathcal{L}}\big)=0, 
\ \ \ \ \ \ \ \ \ \
\rho_0\big(\mathcal{L}\big)=0,
\\
\overline{\zeta}_0\big(\mathcal{T}\big)
&
=
0, 
\ \ \ \ \ \ \ \ \ \
\overline{\zeta}_0\big(\overline{\mathcal{L}}\big)=1, 
\ \ \ \ \ \ \ \ \ \
\zeta_0\big(\mathcal{L}\big)=0,
\\
\zeta_0\big(\mathcal{T}\big)
&
=
0, 
\ \ \ \ \ \ \ \ \ \
\zeta_0\big(\overline{\mathcal{L}}\big)=0, 
\ \ \ \ \ \ \ \ \ \
\zeta_0\big(\mathcal{L}\big)=1,
\endaligned
\]
hence make a {\sl coframe} for:
\[
\C\otimes_\R T^*M.
\]

Explicitly (exercise):
\[
\aligned
\rho_0
&
=
\frac{du-A\,dz-\overline{A}\,d\overline{z}}{\ell},
\\
\overline{\zeta}_0
&
=
d\overline{z},
\\
\zeta_0
&
=
dz.
\endaligned
\]

\medskip\noindent{\bf Cartan formula.}
{\em Given a differential $1$-form:}
\[
\omega
\]
{\em and two vector fields:}
\[
\mathcal{X},
\ \ \ \ \ \ \ \ \ \ \ \ \ \ \ \ \
\mathcal{Y},
\]
{\em one has:}
\[
d\omega(\mathcal{
X},\mathcal{Y})
=
\mathcal{X}\big(\omega(\mathcal{Y})\big)-
\mathcal{Y}\big(\omega(\mathcal{X})\big)-\omega\big([
\mathcal{X},\mathcal{Y}]\big).
\qed
\]

\medskip\noindent
{\bf Lie structure and Darboux structure.}
{\em Given a local frame of vector fields:}
\[
\big\{
\mathcal{X}_1,\dots,\mathcal{X}_n
\big\}
\]
{\em on $\R^n$ or on $\C^n$, and given the dual coframe
of differential $1$-forms:}
\[
\big\{
\omega^1,\dots,\omega^n
\big\},
\]
{\em namely:}
\[
\omega^{\nu_1}\big(\mathcal{X}_{\nu_2}\big)
=
\left\{
\aligned
&
1
\ \ \ \ \
\text{\em when}\ \
\nu_1=\nu_2,
\\
&
0
\ \ \ \ \
\text{\em when}\ \
\nu_1\neq\nu_2,
\endaligned\right.
\]
{\em the Lie structure of the frame is:}
\[
\big[\mathcal{X}_{\nu_1},\,\mathcal{X}_{\nu_2}\big]
\,=\,
\sum_{1\leqslant\nu_3\leqslant n}\,
A_{\nu_1,\nu_2}^{\nu_3}\,
\mathcal{X}_{\nu_3}
\ \ \ \ \ \ \ \ \ \ \ \ \ 
{\scriptstyle{(1\,\leqslant\,\nu_1\,<\,\nu_2\,\leqslant\,n)}},
\]
{\em for certain local functions:}
\[
A_{\nu_1,\nu_2}^{\nu_3},
\]
{\em if and only if the Darboux structure of the coframe is:}
\[
d\omega^{\nu_1}
=
-\,
\sum_{1\leqslant\nu_2<\nu_3\leqslant n}\,
A_{\nu_2,\nu_3}^{\nu_1}\,
\omega^{\nu_2}\wedge\omega^{\nu_3}.
\]

\proof
Use the Cartan formula.
\endproof

To visually
determine the Darboux structure of the dual coframe,
introduce an auxiliary array: 
\[
\footnotesize 
\begin{array}{ccccccccccc} 
& & \mathcal{T} & & \overline{\mathcal{L}} & & \mathcal{L} & & 
\\ 
& & \boxed{d\rho_0} & & \boxed{d\overline{\zeta_0}} & & 
\boxed{d\zeta_0} & & 
\\ 
\big[\mathcal{T},\,\overline{\mathcal{L}}\big] & = & 
-\,\overline{P}\cdot\mathcal{T} & + & 0 & + & 0 & 
\boxed{\rho_0\wedge\overline{\zeta_0}} 
\\ 
\big[\mathcal{T},\,\mathcal{L}\big] & = & -\,P\cdot\mathcal{T} & + & 
0 & + & 0 & \boxed{\rho_0\wedge\zeta_0} 
\\ 
\big[\overline{\mathcal{L}},\,\mathcal{L}\big] & = & 
\isqrt\,\cdot\mathcal{T} & + & 0 & + & 0 & 
\boxed{\overline{\zeta_0}\wedge\zeta_0,}
\end{array} 
\] 
read the three columns {\em vertically}, and put an overall minus
sign:
\[
\label{d0} \boxed{ \aligned d\rho_0 & 
=
P\,\rho_0\wedge\zeta_0
+ 
\overline{P}\,\rho_0\wedge\overline{\zeta_0} 
+
\isqrt\,\zeta_0\wedge\overline{\zeta}_0,
\\ 
d\zeta_0 & = 0, 
\\ 
d\overline{\zeta_0} & = 0. 
\endaligned} 
\]
This is the {\sl initial Darboux structure} for the
problem of biholomorphic equivalence in the
general class $\text{\sf I}$.


\bigskip

\section{\sf $M^4 \subset \C^3$ of general class $\text{\sf II}$: 
\\
initial frame and coframe in local coordinates}
\label{initial-II}
\HEAD{\ref{initial-II}.~$M^4 \subset \C^3$ of general class 
$\text{\sf II}$:
initial frame and coframe in local coordinates}{
Jo\"el {\sc Merker}, D\'epartement de Math\'ematiques d'Orsay}

\medskip

Next, consider:
\[
\Big(
M^4
\,\subset\,
\C^3
\Big)
\,\,\in\,\,
\text{\sf General Class $\text{\sf II}$}.
\]

Representing as before
(\cite{ Merker-Pocchiola-Sabzevari-5-CR-II, Merker-5-CR-III})
$M$ in coordinates:
\[
(z,w_1,w_2) 
= 
\big(
x+\isqrt\,y,\, 
u_1+\isqrt\,v_1,\,
u_2+\isqrt\,v_2\big),
\]
as a graph:
\[
\aligned
v_1
&
=
\varphi_1(x,y,u_1,u_2),
\\
v_2
&
=
\varphi_2(x,y,u_1,u_2),
\endaligned
\]
an associated explicit frame for $TM$ is:
\[
\aligned
X
&
=
\frac{\partial}{\partial x}
+
\varphi_{1,x}(x,y,u_1,u_2)\,\frac{\partial}{\partial v_1}
+
\varphi_{2,x}(x,y,u_1,u_2)\,\frac{\partial}{\partial v_2},
\\
Y
&
=
\frac{\partial}{\partial y}
+
\varphi_{1,y}(x,y,u_1,u_2)\,\frac{\partial}{\partial v_1}
+
\varphi_{2,y}(x,y,u_1,u_2)\,\frac{\partial}{\partial v_2},
\\
U_1
&
=
\frac{\partial}{\partial u_1}
+
\varphi_{1,u_1}(x,y,u_1,u_2)\,\frac{\partial}{\partial v_1}
+
\varphi_{2,u_1}(x,y,u_1,u_2)\,\frac{\partial}{\partial v_2},
\\
U_2
&
=
\frac{\partial}{\partial u_2}
+
\varphi_{1,u_2}(x,y,u_1,u_2)\,\frac{\partial}{\partial v_1}
+
\varphi_{2,u_2}(x,y,u_1,u_2)\,\frac{\partial}{\partial v_2}.
\endaligned
\]
The projection:
\[
\pi\colon
\ \ \ \ \ \ \ \ \ \ \ \ \ \ \
\aligned
\R^6
&
\,\longrightarrow\,
\R^4
\\
(x,y,u_1,v_1,u_2,v_2)
&
\,\longmapsto\,
(x,y,u_1,u_2),
\endaligned
\]
makes a {\sl chart} on $M$ and does:
\[
\aligned
\pi_*
\bigg(
\frac{\partial}{\partial x}
+
\varphi_{1,x}\,
\frac{\partial}{\partial v_1}
+
\varphi_{2,x}\,
\frac{\partial}{\partial v_2}
\bigg)
&
=
\frac{\partial}{\partial x},
\\
\pi_*
\bigg(
\frac{\partial}{\partial y}
+
\varphi_{1,y}\,
\frac{\partial}{\partial v_1}
+
\varphi_{2,y}\,
\frac{\partial}{\partial v_2}
\bigg)
&
=
\frac{\partial}{\partial y},
\\
\pi_*
\bigg(
\frac{\partial}{\partial u_1}
+
\varphi_{1,u_1}\,
\frac{\partial}{\partial v_1}
+
\varphi_{2,u_1}\,
\frac{\partial}{\partial v_2}
\bigg)
&
=
\frac{\partial}{\partial u_1},
\\
\pi_*
\bigg(
\frac{\partial}{\partial u_2}
+
\varphi_{1,u_2}\,
\frac{\partial}{\partial v_1}
+
\varphi_{2,u_2}\,
\frac{\partial}{\partial v_2}
\bigg)
&
=
\frac{\partial}{\partial u_2}.
\endaligned
\]

Next, as was explained in~\cite{Merker-Pocchiola-Sabzevari-5-CR-I}, 
a natural {\em intrinsic} local vector field generator
for $T^{1, 0} M$ is:
\[ 
\mathcal{L} 
:= 
\frac{\partial}{\partial z} 
+ 
A_1\,\frac{\partial}{\partial u_1}
+ 
A_2\,\frac{\partial}{\partial u_2}, 
\] 
where:
\[
\aligned
A_1
&
:=
\frac{
\left\vert\!
\begin{array}{cc}
-\,\varphi_{1,z} & \varphi_{1,u_2}\\
-\,\varphi_{2,z} & \isqrt+\varphi_{2,u_2}
\end{array}
\!\right\vert
}{
\left\vert\!
\begin{array}{cc}
\isqrt+\varphi_{1,u_1} & \varphi_{1,u_2}\\
\varphi_{2,u_1} & \isqrt+\varphi_{2,u_2}
\end{array}
\!\right\vert
},
\\
A_2
&
:=
\frac{
\left\vert\!
\begin{array}{cc}
\isqrt+\varphi_{1,u_1} & -\,\varphi_{1,z} \\
\varphi_{2,u_1} & -\,\varphi_{2,z} 
\end{array}
\!\right\vert
}{
\left\vert\!
\begin{array}{cc}
\isqrt+\varphi_{1,u_1} & \varphi_{1,u_2}\\
\varphi_{2,u_1} & \isqrt+\varphi_{2,u_2}
\end{array}
\!\right\vert}.
\endaligned
\]

Set:
\[
\Delta
:=
\left\vert\!
\begin{array}{cc}
\isqrt+\varphi_{1,u_1} & \varphi_{1,u_2}\\
\varphi_{2,u_1} & \isqrt+\varphi_{2,u_2}
\end{array}
\!\right\vert,
\]
whence:
\[
\overline{\Delta}
:=
\left\vert\!
\begin{array}{cc}
-\isqrt+\varphi_{1,u_1} & \varphi_{1,u_2}\\
\varphi_{2,u_1} & -\isqrt+\varphi_{2,u_2}
\end{array}
\!\right\vert.
\]

Also, set:
\[
\aligned
\Lambda_1
&
:=
\left\vert\!
\begin{array}{cc}
-\,\varphi_{1,z} & \varphi_{1,u_2}\\
-\,\varphi_{2,z} & \isqrt+\varphi_{2,u_2}
\end{array}
\!\right\vert,
\\
\Lambda_2
&
:=
\left\vert\!
\begin{array}{cc}
\isqrt+\varphi_{1,u_1} & -\,\varphi_{1,z} \\
\varphi_{2,u_1} & -\,\varphi_{2,z} 
\end{array}
\!\right\vert,
\endaligned
\]
whence:
\[
\aligned
\overline{\Lambda}_1
&
:=
\left\vert\!
\begin{array}{cc}
-\,\varphi_{1,\overline{z}} & \varphi_{1,u_2}\\
-\,\varphi_{2,\overline{z}} & -\isqrt+\varphi_{2,u_2}
\end{array}
\!\right\vert,
\\
\overline{\Lambda}_2
&
:=
\left\vert\!
\begin{array}{cc}
-\isqrt+\varphi_{1,u_1} & -\,\varphi_{1,\overline{z}} \\
\varphi_{2,u_1} & -\,\varphi_{2,\overline{z}} 
\end{array}
\!\right\vert.
\endaligned
\]

In these notations:
\[
\aligned
\mathcal{L}
&
=
\frac{\partial}{\partial z}
+
\frac{\Lambda_1}{\Delta}\,\frac{\partial}{\partial u_1}
+
\frac{\Lambda_2}{\Delta}\,\frac{\partial}{\partial u_2},
\\
\overline{\mathcal{L}}
&
=
\frac{\partial}{\partial\overline{z}}
+
\frac{\overline{\Lambda}_1}{\overline{\Delta}}\,
\frac{\partial}{\partial u_1}
+
\frac{\overline{\Lambda}_2}{\overline{\Delta}}\,
\frac{\partial}{\partial u_2}.
\endaligned
\]

\medskip

Set:
\[
\mathcal{T}
:=
\isqrt\,\big[\mathcal{L},\overline{\mathcal{L}}\big],
\]
and set:
\[
\mathcal{S}
:=
\big[\mathcal{L},\mathcal{T}\big].
\]

By hypothesis, the CR-geometric invariant condition:
\[
\aligned
\C\otimes_\R TM
&
=
T^{1,0}M+T^{0,1}M
+
[T^{1,0}M,\,T^{0,1}M]
+
\big[T^{1,0}M,\,[T^{1,0}M,\,T^{0,1}M]\big]
\endaligned
\]
holds, which means as is known that the $4$ fields:
\[
\big\{
\mathcal{L},\,\overline{\mathcal{L}},\,
\mathcal{T},\,\mathcal{S}
\big\}
\]
constitute a frame for:
\[
{\bf 4}
=
\rank_\C\big(
\C\otimes_\R TM
\big).
\]

\medskip\noindent{\bf Lemma.}
{\em There are certain uniquely defined 
coefficient-functions that are polynomials:}
\[
\aligned
\Upsilon_1
&
=
\Upsilon_1
\Big(
\varphi_{1,x^jy^ku_1^{l_1}u_2^{l_2}},\,\,
\varphi_{2,x^jy^ku_1^{l_1}u_2^{l_2}}
\Big)_{1\leqslant j+k+l_1+l_2\leqslant 2},
\\
\Upsilon_2
&
=
\Upsilon_2
\Big(
\varphi_{1,x^jy^ku_1^{l_1}u_2^{l_2}},\,\,
\varphi_{2,x^jy^ku_1^{l_1}u_2^{l_2}}
\Big)_{1\leqslant j+k+l_1+l_2\leqslant 2},
\endaligned
\]
\[
\aligned
\Pi_1
&
=
\Pi_1
\Big(
\varphi_{1,x^jy^ku_1^{l_1}u_2^{l_2}},\,\,
\varphi_{2,x^jy^ku_1^{l_1}u_2^{l_2}}
\Big)_{1\leqslant j+k+l_1+l_2\leqslant 3},
\\
\Pi_2
&
=
\Pi_2
\Big(
\varphi_{1,x^jy^ku_1^{l_1}u_2^{l_2}},\,\,
\varphi_{2,x^jy^ku_1^{l_1}u_2^{l_2}}
\Big)_{1\leqslant j+k+l_1+l_2\leqslant 3},
\endaligned
\] 
{\em such that:}
\[
\aligned
\mathcal{T}
&
=
\frac{\Upsilon_1}{\Delta^2\,\overline{\Delta}^2}\,
\frac{\partial}{\partial u_1}
+
\frac{\Upsilon_2}{\Delta^2\,\overline{\Delta}^2}\,
\frac{\partial}{\partial u_2}
\\
&
=:
Y_1\,
\frac{\partial}{\partial u_1}
+
Y_2\,
\frac{\partial}{\partial u_2},
\\
\mathcal{S}
&
=
\frac{\Pi_1}{\Delta^4\,\overline{\Delta}^3}\,
\frac{\partial}{\partial u_1}
+
\frac{\Pi_2}{\Delta^4\,\overline{\Delta}^3}\,
\frac{\partial}{\partial u_2}
\\
&
=:
H_1\,
\frac{\partial}{\partial u_1}
+
H_2\,
\frac{\partial}{\partial u_2}.
\endaligned
\]

\proof
In fact firstly:
\[
\aligned
\isqrt\,
\big[
\mathcal{L},\overline{\mathcal{L}}
\big]
&
=
\isqrt\,
\bigg[
\frac{\partial}{\partial z}
+
\frac{\Lambda_1}{\Delta}\,
\frac{\partial}{\partial u_1}
+
\frac{\Lambda_2}{\Delta}\,
\frac{\partial}{\partial u_2},
\,\,\,\,
\frac{\partial}{\partial\overline{z}}
+
\frac{\overline{\Lambda}_1}{\overline{\Delta}}\,
\frac{\partial}{\partial u_1}
+
\frac{\overline{\Lambda}_2}{\overline{\Delta}}\
\frac{\partial}{\partial u_2}
\bigg]
\\
&
=
\isqrt\,
\bigg\{
\mathcal{L}
\bigg(
\frac{\overline{\Lambda}_1}{\overline{\Delta}}
\bigg)
-
\overline{\mathcal{L}}
\bigg(
\frac{\Lambda_1}{\Delta}
\bigg)
\bigg\}
\frac{\partial}{\partial u_1}
-
\isqrt\,
\bigg\{
\mathcal{L}
\bigg(
\frac{\overline{\Lambda}_2}{\overline{\Delta}}
\bigg)
-
\overline{\mathcal{L}}
\bigg(
\frac{\Lambda_2}{\Delta}
\bigg)
\bigg\}
\frac{\partial}{\partial u_2}
\endaligned
\]
Compute:
\[
\aligned
\mathcal{L}
\bigg(
\frac{\overline{\Lambda}_1}{\overline{\Delta}}
\bigg)
&
=
\frac{\mathcal{L}\big(\overline{\Lambda}_1\big)}{
\overline{\Delta}}
-
\frac{\overline{\Lambda}_1\,\mathcal{L}\big(\overline{\Delta}\big)}{
\overline{\Delta}^2}
\\
&
=
\frac{\Delta\,\overline{\Lambda}_{1,z}
+
\Lambda_1\,\overline{\Lambda}_{1,u_1}
+
\Lambda_2\,\overline{\Lambda}_{1,u_2}}{
\Delta\,\overline{\Delta}}
\,-
\\
&
-\,
\frac{
\overline{\Lambda}_1
\big(
\Delta\,\overline{\Delta}_z
+
\Lambda_1\,\overline{\Delta}_{u_1}
+
\Lambda_2\,\overline{\Delta}_{u_2}
\big)}{
\Delta\,\overline{\Delta}^2},
\endaligned
\]
conjugate this:
\[
\aligned
\overline{\mathcal{L}}
\bigg(
\frac{\Lambda_1}{\Delta}
\bigg)
&
=
\frac{\overline{\Delta}\,\Lambda_{1,\overline{z}}
+
\overline{\Lambda}_1\,\Lambda_{1,u_1}
+
\overline{\Lambda}_2\,\Lambda_{1,u_2}}{
\Delta\,\overline{\Delta}}
\,-
\\
&
-
\frac{\Lambda_1\big(
\overline{\Delta}\,\Delta_{\overline{z}}
+
\overline{\Lambda}_1\,\Delta_{u_1}
+
\overline{\Lambda}_2\,\Delta_{u_2}
\big)}{
\Delta^2\,\overline{\Delta}}
\endaligned
\]
whence after reduction to the common denominator
$\Delta^2\, \overline{\Delta}^2$:
\[
\aligned
\isqrt\,
\bigg\{
\mathcal{L}
\bigg(
\frac{\overline{\Lambda}_1}{\overline{\Delta}}
\bigg)
-
\overline{\mathcal{L}}
\bigg(
\frac{\Lambda_1}{\Delta}
\bigg)
\bigg\}
&
=
\frac{\isqrt}{\Delta^2\,\overline{\Delta}^2}\,
\bigg\{
\Delta\,\Delta\,\overline{\Delta}\,
\overline{\Lambda}_{1,z}
+
\Delta\,\overline{\Delta}\,\Lambda_1\,
\overline{\Lambda}_{1,u_1}
+
\Delta\,\overline{\Delta}\,\Lambda_2\,
\overline{\Lambda}_{1,u_2}
\,-
\\
&
\ \ \ \ \ \ \ \ \ \ \ \ \ \ \ 
-\,
\Delta\,\Delta\,\overline{\Delta}_z\,
\overline{\Lambda}_1
-
\Delta\,\Lambda_1\,\overline{\Delta}_{u_1}\,
\overline{\Lambda}_1
-
\Delta\,\Lambda_2\,\overline{\Delta}_{u_2}\,
\overline{\Lambda}_1
\,-
\\
&
\ \ \ \ \ \ \ \ \ \ \ \ \ \ \ 
-\,
\Delta\,\overline{\Delta}\,\overline{\Delta}\,
\Lambda_{1,\overline{z}}
-
\Delta\,\overline{\Delta}\,\overline{\Lambda}_1\,
\Lambda_{1,u_1}
-
\Delta\,\overline{\Delta}\,\overline{\Lambda}_z\,
\Lambda_{1,u_2}
+
\\
&
\ \ \ \ \ \ \ \ \ \ \ \ \ \ \ 
+
\overline{\Delta}\,\overline{\Delta}\,\Delta_{\overline{z}}\,
\Lambda_1
+
\overline{\Delta}\,\overline{\Lambda}_1\,\Delta_{u_1}\,
\Lambda_1
+
\overline{\Delta}\,\overline{\Lambda}_2\,\Delta_{u_2}\,
\Lambda_1
\bigg\},
\endaligned
\]
which provides:
\[
\boxed{\,
\aligned
\Upsilon_1
&
:=
\isqrt\,\Big(
\Delta\,\Delta\,\overline{\Delta}\,
\overline{\Lambda}_{1,z}
+
\Delta\,\overline{\Delta}\,\Lambda_1\,
\overline{\Lambda}_{1,u_1}
+
\Delta\,\overline{\Delta}\,\Lambda_2\,
\overline{\Lambda}_{1,u_2}
\,-\,\,
\\
&
\ \ \ \ \ \ \ \ \ \ \ \ \ \ \ 
-\,
\Delta\,\Delta\,\overline{\Delta}_z\,
\overline{\Lambda}_1
-
\Delta\,\Lambda_1\,\overline{\Delta}_{u_1}\,
\overline{\Lambda}_1
-
\Delta\,\Lambda_2\,\overline{\Delta}_{u_2}\,
\overline{\Lambda}_1
\,-
\\
&
\ \ \ \ \ \ \ \ \ \ \ \ \ \ \ 
-\,
\Delta\,\overline{\Delta}\,\overline{\Delta}\,
\Lambda_{1,\overline{z}}
-
\Delta\,\overline{\Delta}\,\overline{\Lambda}_1\,
\Lambda_{1,u_1}
-
\Delta\,\overline{\Delta}\,\overline{\Lambda}_z\,
\Lambda_{1,u_2}
+
\\
&
\ \ \ \ \ \ \ \ \ \ \ \ \ \ \ 
+
\overline{\Delta}\,\overline{\Delta}\,\Delta_{\overline{z}}\,
\Lambda_1
+
\overline{\Delta}\,\overline{\Lambda}_1\,\Delta_{u_1}\,
\Lambda_1
+
\overline{\Delta}\,\overline{\Lambda}_2\,\Delta_{u_2}\,
\Lambda_1
\Big),
\endaligned
\,}
\]

Similarly:
\[
\isqrt\,
\bigg\{
\mathcal{L}
\bigg(
\frac{\overline{\Lambda}_2}{\overline{\Delta}}
\bigg)
-
\overline{\mathcal{L}}
\bigg(
\frac{\Lambda_2}{\Delta}
\bigg)
\bigg\}
=
\frac{\Upsilon_2}{\Delta^2\,\overline{\Delta}^2},
\]
with:
\[
\boxed{\,
\aligned
\Upsilon_2
&
:=
\isqrt\,\Big(
\Delta\,\Delta\,\overline{\Delta}\,
\overline{\Lambda}_{2,z}
+
\Delta\,\overline{\Delta}\,\Lambda_1\,
\overline{\Lambda}_{2,u_1}
+
\Delta\,\overline{\Delta}\,\Lambda_2\,
\overline{\Lambda}_{2,u_2}
\,-\,\,
\\
&
\ \ \ \ \ \ \ \ \ \ \ \ \ \ \ 
-\,
\Delta\,\Delta\,\overline{\Delta}_z\,
\overline{\Lambda}_2
-
\Delta\,\Lambda_1\,\overline{\Delta}_{u_1}\,
\overline{\Lambda}_2
-
\Delta\,\Lambda_2\,\overline{\Delta}_{u_2}\,
\overline{\Lambda}_2
\,-
\\
&
\ \ \ \ \ \ \ \ \ \ \ \ \ \ \ 
-\,
\Delta\,\overline{\Delta}\,\overline{\Delta}\,
\Lambda_{2,\overline{z}}
-
\Delta\,\overline{\Delta}\,\overline{\Lambda}_1\,
\Lambda_{2,u_1}
-
\Delta\,\overline{\Delta}\,\overline{\Lambda}_z\,
\Lambda_{2,u_2}
+
\\
&
\ \ \ \ \ \ \ \ \ \ \ \ \ \ \ 
+
\overline{\Delta}\,\overline{\Delta}\,\Delta_{\overline{z}}\,
\Lambda_2
+
\overline{\Delta}\,\overline{\Lambda}_1\,\Delta_{u_1}\,
\Lambda_2
+
\overline{\Delta}\,\overline{\Lambda}_2\,\Delta_{u_2}\,
\Lambda_2
\Big).
\endaligned
\,}
\]

\medskip

Secondly:
\[
\aligned
\big[\mathcal{L},\mathcal{T}\big]
&
=
\bigg[
\frac{\partial}{\partial z}
+
\frac{\Lambda_1}{\Delta}\,
\frac{\partial}{\partial u_1}
+
\frac{\Lambda_2}{\Delta}\,
\frac{\partial}{\partial u_2},
\,\,\,\,
\frac{\Upsilon_1}{\Delta^2\,\overline{\Delta}^2}\,
\frac{\partial}{\partial u_1}
+
\frac{\Upsilon_2}{\Delta^2\,\overline{\Delta}^2}\,
\frac{\partial}{\partial u_2}
\bigg]
\\
&
=
\bigg\{
\mathcal{L}
\bigg(
\frac{\Upsilon_1}{\Delta^2\,\overline{\Delta}^2}
\bigg)
-
\mathcal{T}
\bigg(
\frac{\Lambda_1}{\Delta}
\bigg)
\bigg\}\,
\frac{\partial}{\partial u_1}
+
\bigg\{
\mathcal{L}
\bigg(
\frac{\Upsilon_2}{\Delta^2\,\overline{\Delta}^2}
\bigg)
-
\mathcal{T}
\bigg(
\frac{\Lambda_2}{\Delta}
\bigg)
\bigg\}\,
\frac{\partial}{\partial u_2}.
\endaligned
\]

Compute:
\[
\aligned
\mathcal{L}
\bigg(
\frac{\Upsilon_1}{\Delta^2\,\overline{\Delta}^2}
\bigg)
&
=
\frac{\mathcal{L}\big(\Upsilon_1\big)}{
\Delta^2\,\overline{\Delta}^2}
-
\frac{\Upsilon_1\,2\,\mathcal{L}(\Delta)}{
\Delta^3\,\overline{\Delta}^2}
-
\frac{\Upsilon_1\,2\,\mathcal{L}\big(\overline{\Delta}\big)}{
\Delta^2\,\overline{\Delta}^3}
\\
&
=
\frac{
\Delta\,\Upsilon_{1,z}
+
\Lambda_1\,\Upsilon_{1,u_1}
+
\Lambda_2\,\Upsilon_{1,u_2}}{
\Delta^3\,\overline{\Delta}^2}
\,-
\\
&
-\,
\frac{\Upsilon_1\big(
2\,\Delta_z\,\Delta
+
2\,\Lambda_1\,\Delta_{u_1}
+
2\,\Lambda_2\,\Delta_{u_2}
\big)}{
\Delta^4\,\overline{\Delta}^2}
\,-
\\
&
-\,
\frac{\Upsilon_1\big(
2\,\overline{\Delta}_z\,\Delta
+
2\,\Lambda_1\,\overline{\Delta}_{u_1}
+
2\,\Lambda_2\,\overline{\Delta}_{u_2}
\big)}{
\Delta^3\,\overline{\Delta}^3},
\endaligned
\]
whence after reduction to the common denominator
$\Delta^4\, \overline{\Delta}^3$:
\[
\aligned
\mathcal{L}
\bigg(
\frac{\Upsilon_1}{\Delta^2\,\overline{\Delta}^2}
\bigg)
&
=
\frac{1}{\Delta^4\,\overline{\Delta}^3}\,
\bigg\{
\Delta\,\Delta\,\overline{\Delta}\,
\Upsilon_{1,z}
+
\Delta\,\overline{\Delta}\,\Lambda_1\,
\Upsilon_{1,u_1}
+
\Delta\,\overline{\Delta}\,\Lambda_2\,
\Upsilon_{1,u_2}
\,-
\\
&
\ \ \ \ \ \ \ \ \ \ \ \ \ \ \
\,-
2\,\Delta\,\overline{\Delta}\,\Delta_z\,
\Upsilon_1
-
2\,\overline{\Delta}\,\Lambda_1\,\Delta_{u_1}\,
\Upsilon_1
-
2\,\overline{\Delta}\,\Lambda_2\,\Delta_{u_2}\,
\Upsilon_1
\,-
\\
&
\ \ \ \ \ \ \ \ \ \ \ \ \ \ \
\,-
2\,\Delta\,\Delta\,\overline{\Delta}_z\,
\Upsilon_1
-
2\,\Delta\,\Lambda_1\,\overline{\Delta}_{u_1}\,
\Upsilon_1
-
2\,\Delta\,\Lambda_2\,\overline{\Delta}_{u_2}\,
\Upsilon_1
\bigg\}.
\endaligned
\]

Compute also:
\[
\aligned
\mathcal{T}
\bigg(
\frac{\Lambda_1}{\Delta}
\bigg)
&
=
\frac{\mathcal{T}\big(\Lambda_1\big)}{\Delta}
-
\frac{\Lambda_1\,\mathcal{T}\big(\Delta\big)}{\Delta^2}
\\
&
=
\frac{\Upsilon_1\,\Lambda_{1,u_1}
+
\Upsilon_2\,\Lambda_{1,u_2}}{
\Delta^3\,\overline{\Delta}^2}
\,-
\\
&
-\,
\frac{\Lambda_1\big(
\Upsilon_1\,\Delta_{u_1}
+
\Upsilon_2\,\Delta_{u_2}\big)}{
\Delta^4\,\overline{\Delta}^2}
\\
&
=
\frac{1}{\Delta^4\,\overline{\Delta}^3}\,
\Big\{
\Delta\,\overline{\Delta}\,\Upsilon_1\,
\Lambda_{1,u_1}
+
\Delta\,\overline{\Delta}\,\Upsilon_2\,
\Lambda_{1,u_2}
\,-
\\
&
\ \ \ \ \ \ \ \ \ \ \ \ \ \ \
\,-
\overline{\Delta}\,\Upsilon_1\,\Delta_{u_1}\,
\Lambda_1
-
\overline{\Delta}\,\Upsilon_2\,\Delta_{u_2}\,
\Lambda_1
\Big\}.
\endaligned
\]

In sum:
\[
\aligned
\mathcal{L}
\bigg(
\frac{\Upsilon_1}{\Delta^2\,\overline{\Delta}^2}
\bigg)
-
\mathcal{T}
\bigg(
\frac{\Lambda_1}{\Delta}
\bigg)
=
\frac{\Pi_1}{\Delta^4\,\overline{\Delta}^3},
\endaligned
\]
with:
\[
\boxed{\,
\aligned
\Pi_1
&
=
\Delta\,\Delta\,\overline{\Delta}\,
\Upsilon_{1,z}
+
\Delta\,\overline{\Delta}\,\Lambda_1\,
\Upsilon_{1,u_1}
+
\Delta\,\overline{\Delta}\,\Lambda_2\,
\Upsilon_{1,u_2}
\,-
\\
&
\,-
2\,\Delta\,\overline{\Delta}\,\Delta_z\,
\Upsilon_1
-
2\,\overline{\Delta}\,\Lambda_1\,\Delta_{u_1}\,
\Upsilon_1
-
2\,\overline{\Delta}\,\Lambda_2\,\Delta_{u_2}\,
\Upsilon_1
\,-\,\,
\\
&
\,-
2\,\Delta\,\Delta\,\overline{\Delta}_z\,
\Upsilon_1
-
2\,\Delta\,\Lambda_1\,\overline{\Delta}_{u_1}\,
\Upsilon_1
-
2\,\Delta\,\Lambda_2\,\overline{\Delta}_{u_2}\,
\Upsilon_1
\,-
\\
&
\ \ \ \ \ \ \ \ \ \ \ \ \ \ \ \ \ \ \ \ \ \ \ \ \ \ \
-\,
\Delta\,\overline{\Delta}\,\Upsilon_1\,
\Lambda_{1,u_1}
-
\Delta\,\overline{\Delta}\,\Upsilon_2\,
\Lambda_{1,u_2}
+
\\
&
\ \ \ \ \ \ \ \ \ \ \ \ \ \ \ \ \ \ \ \ \ \ \ \ \ \ \
+
\overline{\Delta}\,\Upsilon_1\,\Delta_{u_1}\,
\Lambda_1
+
\overline{\Delta}\,\Upsilon_2\,\Delta_{u_2}\,
\Lambda_1.
\endaligned
}
\]

Similarly:
\[
\aligned
\mathcal{L}
\bigg(
\frac{\Upsilon_2}{\Delta^2\,\overline{\Delta}^2}
\bigg)
-
\mathcal{T}
\bigg(
\frac{\Lambda_2}{\Delta}
\bigg)
=
\frac{\Pi_2}{\Delta^4\,\overline{\Delta}^3},
\endaligned
\]
with:
\[
\boxed{\,
\aligned
\Pi_2
&
=
\Delta\,\Delta\,\overline{\Delta}\,
\Upsilon_{2,z}
+
\Delta\,\overline{\Delta}\,\Lambda_1\,
\Upsilon_{2,u_1}
+
\Delta\,\overline{\Delta}\,\Lambda_2\,
\Upsilon_{2,u_2}
\,-
\\
&
\,-
2\,\Delta\,\overline{\Delta}\,\Delta_z\,
\Upsilon_2
-
2\,\overline{\Delta}\,\Lambda_1\,\Delta_{u_1}\,
\Upsilon_2
-
2\,\overline{\Delta}\,\Lambda_2\,\Delta_{u_2}\,
\Upsilon_2
\,-\,\,
\\
&
\,-
2\,\Delta\,\Delta\,\overline{\Delta}_z\,
\Upsilon_2
-
2\,\Delta\,\Lambda_1\,\overline{\Delta}_{u_1}\,
\Upsilon_2
-
2\,\Delta\,\Lambda_2\,\overline{\Delta}_{u_2}\,
\Upsilon_2
\,-
\\
&
\ \ \ \ \ \ \ \ \ \ \ \ \ \ \ \ \ \ \ \ \ \ \ \ \ \ \
-\,
\Delta\,\overline{\Delta}\,\Upsilon_1\,
\Lambda_{2,u_1}
-
\Delta\,\overline{\Delta}\,\Upsilon_2\,
\Lambda_{2,u_2}
+
\\
&
\ \ \ \ \ \ \ \ \ \ \ \ \ \ \ \ \ \ \ \ \ \ \ \ \ \ \
+
\overline{\Delta}\,\Upsilon_1\,\Delta_{u_1}\,
\Lambda_2
+
\overline{\Delta}\,\Upsilon_2\,\Delta_{u_2}\,
\Lambda_2,
\endaligned
}
\]
which concludes.
\endproof

\noindent{\bf Explicitness obstacle.}
{\em After reduction to a common minimal denominator, the two numerators
in:}
\[
Y_1
=
\frac{\Upsilon_1}{
\Delta^2\,\overline{\Delta}^2},
\ \ \ \ \ \ \ \ \ \ \ \ \ \ \ \ \ \ \ \ \ \ \ \ \ 
Y_2
=
\frac{\Upsilon_2}{
\Delta^2\,\overline{\Delta}^2},
\]
{\em are both polynomials in the ${\bf 2}\cdot{\bf 14}$ partial derivatives:}
\[
\Big(
\varphi_{1,x^jy^ku_1^{l_1}u_2^{l_2}},\,\,
\varphi_{2,x^jy^ku_1^{l_1}u_2^{l_2}}
\Big)_{1\leqslant j+k+l_1+l_2\leqslant 2}
\]
{\em incorporating:}
\[
{\bf 355}
\]
{\em monomials, while the next two numerators in:}
\[
H_1
=
\frac{\Pi_1}{
\Delta^4\,\overline{\Delta}^3},
\ \ \ \ \ \ \ \ \ \ \ \ \ \ \ \ \ \ \ \ \ \ \ \ \ 
H_2
=
\frac{\Pi_2}{
\Delta^4\,\overline{\Delta}^3}
\]
{\em incorporate both:}
\[
{\bf 24\,437}
\]
{\em monomials in the ${\bf 2} \cdot {\bf 34}$ partial derivatives:}
\[
\Big(
\varphi_{1,x^jy^ku_1^{l_1}u_2^{l_2}},\,\,
\varphi_{2,x^jy^ku_1^{l_1}u_2^{l_2}}
\Big)_{1\leqslant j+k+l_1+l_2\leqslant 3}.
\]

\proof
Use a computer software.
\endproof

The Class $\text{\sf II}$ hypothesis now reads as the nonzeroness:
\[
0
\,\neq\,
\det\,
\left(\!
\begin{array}{cc}
\frac{\Upsilon_1}{\Delta^2\overline{\Delta}^2} &
\frac{\Upsilon_2}{\Delta^2\overline{\Delta}^2}
\medskip
\\
\frac{\Pi_1}{\Delta^4\overline{\Delta}^3} &
\frac{\Pi_2}{\Delta^4\overline{\Delta}^3}
\end{array}
\!\right)
(x,y,u_1,u_2),
\]
at every point.

More precisely, if one preliminarily normalizes
coordinates as in~\cite{ Merker-5-CR-III}:
\[
\aligned
v_1
&
=
z\overline{z}
\ \ \ \ \ \ \ \ \ \ \ \ 
+
z\overline{z}\,{\rm O}_2\big(z,\overline{z}\big)
+
z\overline{z}\,{\rm O}_1(u_1)
+
z\overline{z}\,{\rm O}_1(u_2),\,\,
\\
v_2
&
=
z^2\overline{z}
+
z\overline{z}^2
+
z\overline{z}\,{\rm O}_2\big(z,\overline{z}\big)
+
z\overline{z}\,{\rm O}_1(u_1)
+
z\overline{z}\,{\rm O}_1(u_2),\,\,
\endaligned
\]
so that at the origin:
\[
\aligned
\mathcal{L}\big\vert_0
&
=
\frac{\partial}{\partial z}
\bigg\vert_0,
\\
\overline{\mathcal{L}}\big\vert_0
&
=
\frac{\partial}{\partial\overline{z}}
\bigg\vert_0,
\\
\big[\mathcal{L},\overline{\mathcal{L}}\big]
\Big\vert_0
&
=
-\,2\,\isqrt\,
\frac{\partial}{\partial u_1}
\bigg\vert_0,
\\
\big[\mathcal{L},\,\big[\mathcal{L},\overline{\mathcal{L}}\big]\big]
\Big\vert_0
&
=
-\,4\,\isqrt\,
\frac{\partial}{\partial u_2}
\bigg\vert_0,
\endaligned
\]
the determinant in question:
\[
\det\,
\left(\!
\begin{array}{cc}
\frac{\Upsilon_1}{\Delta^2\overline{\Delta}^2} &
\frac{\Upsilon_2}{\Delta^2\overline{\Delta}^2}
\medskip
\\
\frac{\Pi_1}{\Delta^4\overline{\Delta}^3} &
\frac{\Pi_2}{\Delta^4\overline{\Delta}^3}
\end{array}
\!\right)
(0)
\,=\,
\left\vert\!
\begin{array}{cc}
2 & 0
\\
0 & 4
\end{array}
\!\right\vert
\]
becomes quite visibly nonzero, hence also near the origin.

\medskip

But generally, the disease is, that when (necessarily) re-expressing:
\[
\aligned
\frac{\partial}{\partial u_1}
&
=
\frac{\Delta^2\,\overline{\Delta}^2\,\Pi_2}{
\left\vert\!
\begin{array}{cc}
\Upsilon_1 & \Upsilon_2\\
\Pi_1 & \Pi_2
\end{array}
\!\right\vert}
\cdot
\mathcal{T}
-
\frac{\Delta^4\,\overline{\Delta}^3\,\Upsilon_2}{
\left\vert\!
\begin{array}{cc}
\Upsilon_1 & \Upsilon_2\\
\Pi_1 & \Pi_2
\end{array}
\!\right\vert}
\cdot
\mathcal{S},
\\
\frac{\partial}{\partial u_2}
&
=
\frac{\Delta^2\,\overline{\Delta}^2\,\Pi_1}{
\left\vert\!
\begin{array}{cc}
\Upsilon_1 & \Upsilon_2\\
\Pi_1 & \Pi_2
\end{array}
\!\right\vert}
\cdot
\mathcal{T}
-
\frac{\Delta^4\,\overline{\Delta}^3\,\Upsilon_1}{
\left\vert\!
\begin{array}{cc}
\Upsilon_1 & \Upsilon_2\\
\Pi_1 & \Pi_2
\end{array}
\!\right\vert}
\cdot
\mathcal{S},
\endaligned
\]
some quite huge fractions appear, and it becomes unwieldy,
even for a powerful computer, to make explicit
the coefficients in the next $4$ brackets completing
the Lie structure of the frame:
\[
\aligned
\big[\overline{\mathcal{L}},\mathcal{T}\big]
&
=
\function
\cdot
\mathcal{T}
+
\function
\cdot
\mathcal{S},
\\
\big[\mathcal{L},\mathcal{S}\big]
&
=
\function
\cdot
\mathcal{T}
+
\function
\cdot
\mathcal{S},
\\
\big[\overline{\mathcal{L}},\mathcal{S}\big]
&
=
\function
\cdot
\mathcal{T}
+
\function
\cdot
\mathcal{S},
\\
\big[\mathcal{T},\mathcal{S}\big]
&
=
\function
\cdot
\mathcal{T}
+
\function
\cdot
\mathcal{S}.
\endaligned
\]

\medskip\noindent{\bf About a classical simplified context.}
In the (special) so-called {\sl rigid} case:
\[
\aligned
v_1
&
=
\varphi_1(x,y),
\\
v_2
&
=
\varphi_2(x,y),
\endaligned
\]
where the graphing functions do not depend upon
$u_1$, $u_2$, completely explicit formulas can be typesetted.

\medskip

Indeed:
\[
\aligned
\mathcal{L}
&
=
\frac{\partial}{\partial z}
+
\isqrt\,\varphi_{1,z}\,
\frac{\partial}{\partial u_1}
+
\isqrt\,\varphi_{2,z}\,
\frac{\partial}{\partial u_2},
\\
\overline{\mathcal{L}}
&
=
\frac{\partial}{\partial\overline{z}}
-
\isqrt\,\varphi_{1,\overline{z}}\,
\frac{\partial}{\partial u_1}
-
\isqrt\,\varphi_{2,\overline{z}}\,
\frac{\partial}{\partial u_2},
\\
\mathcal{T}
&
=
\ \ \ \ \ \ \ \ \ \ \ \ \,
2\,\varphi_{1,z\overline{z}}\,
\frac{\partial}{\partial u_1}
+
2\,\varphi_{2,z\overline{z}}\,
\frac{\partial}{\partial u_2},
\\
\mathcal{S}
&
=
\ \ \ \ \ \ \ \ \ \ \ \ \,
2\,\varphi_{1,zz\overline{z}}\,
\frac{\partial}{\partial u_1}
+
2\,\varphi_{2,zz\overline{z}}\,
\frac{\partial}{\partial u_2}.
\endaligned
\]

The determinant in question is:
\[
4\,\big(
\varphi_{1,z\overline{z}}\,\varphi_{2,zz\overline{z}}
-
\varphi_{2,z\overline{z}}\,\varphi_{1,zz\overline{z}}
\big),
\]
and one has:
\[
\aligned
\frac{\partial}{\partial u_1}
&
=
\frac{\varphi_{2,zz\overline{z}}}{
2\,(
\varphi_{1,z\overline{z}}\,\varphi_{2,zz\overline{z}}
-
\varphi_{2,z\overline{z}}\,\varphi_{1,zz\overline{z}}
)}
\cdot
\mathcal{T}
-
\frac{\varphi_{2,z\overline{z}}}{
2\,(
\varphi_{1,z\overline{z}}\,\varphi_{2,zz\overline{z}}
-
\varphi_{2,z\overline{z}}\,\varphi_{1,zz\overline{z}}
)}
\cdot
\mathcal{S},
\\
\frac{\partial}{\partial u_2}
&
=
-\,
\frac{\varphi_{1,zz\overline{z}}}{
2\,(
\varphi_{1,z\overline{z}}\,\varphi_{2,zz\overline{z}}
-
\varphi_{2,z\overline{z}}\,\varphi_{1,zz\overline{z}}
)}
\cdot
\mathcal{T}
+
\frac{\varphi_{1,z\overline{z}}}{
2\,(
\varphi_{1,z\overline{z}}\,\varphi_{2,zz\overline{z}}
-
\varphi_{2,z\overline{z}}\,\varphi_{1,zz\overline{z}}
)}
\cdot
\mathcal{S}.
\endaligned
\]

Next:
\[
\aligned
\big[\overline{\mathcal{L}},\mathcal{T}\big]
&
=
2\,\varphi_{1,z\overline{z}\overline{z}}\,
\frac{\partial}{\partial u_1}
+
2\,\varphi_{2,z\overline{z}\overline{z}}\,
\frac{\partial}{\partial u_2}
\\
&
=
\bigg(
\frac{\varphi_{1,z\overline{z}\overline{z}}\,
\varphi_{2,zz\overline{z}}
-
\varphi_{2,z\overline{z}\overline{z}}\,\varphi_{1,zz\overline{z}}}{
\varphi_{1,z\overline{z}}\,\varphi_{2,zz\overline{z}}
-
\varphi_{2,z\overline{z}}\,\varphi_{1,zz\overline{z}}}
\bigg)\cdot
\mathcal{T}
+
\bigg(
\frac{
-\,\varphi_{1,z\overline{z}\overline{z}}\,
\varphi_{2,z\overline{z}}
+
\varphi_{2,z\overline{z}\overline{z}}\,\varphi_{1,z\overline{z}}}{
\varphi_{1,z\overline{z}}\,\varphi_{2,zz\overline{z}}
-
\varphi_{2,z\overline{z}}\,\varphi_{1,zz\overline{z}}}
\bigg)\cdot
\mathcal{S},
\endaligned
\]
and:
\[
\aligned
\big[\mathcal{L},\mathcal{S}\big]
&
=
2\,\varphi_{1,zzz\overline{z}}\,
\frac{\partial}{\partial u_1}
+
2\,\varphi_{2,zzz\overline{z}}\,
\frac{\partial}{\partial u_2}
\\
&
=
\bigg(
\frac{\varphi_{1,zzz\overline{z}}\,
\varphi_{2,zz\overline{z}}
-
\varphi_{2,zzz\overline{z}}\,\varphi_{1,zz\overline{z}}}{
\varphi_{1,z\overline{z}}\,\varphi_{2,zz\overline{z}}
-
\varphi_{2,z\overline{z}}\,\varphi_{1,zz\overline{z}}}
\bigg)\cdot
\mathcal{T}
+
\bigg(
\frac{
-\,\varphi_{1,zzz\overline{z}}\,
\varphi_{2,z\overline{z}}
+
\varphi_{2,zzz\overline{z}}\,\varphi_{1,z\overline{z}}}{
\varphi_{1,z\overline{z}}\,\varphi_{2,zz\overline{z}}
-
\varphi_{2,z\overline{z}}\,\varphi_{1,zz\overline{z}}}
\bigg)\cdot
\mathcal{S},
\endaligned
\]
and also:
\[
\aligned
\big[\overline{\mathcal{L}},\mathcal{S}\big]
&
=
2\,\varphi_{1,zz\overline{z}\overline{z}}\,
\frac{\partial}{\partial u_1}
+
2\,\varphi_{2,zz\overline{z}\overline{z}}\,
\frac{\partial}{\partial u_2}
\\
&
=
\bigg(
\frac{\varphi_{1,zz\overline{z}\overline{z}}\,
\varphi_{2,zz\overline{z}}
-
\varphi_{2,zz\overline{z}\overline{z}}\,\varphi_{1,zz\overline{z}}
}{
\varphi_{1,z\overline{z}}\,\varphi_{2,zz\overline{z}}
-
\varphi_{2,z\overline{z}}\,\varphi_{1,zz\overline{z}}}
\bigg)\cdot
\mathcal{T}
+
\bigg(
\frac{
-\,\varphi_{1,zz\overline{z}\overline{z}}\,
\varphi_{2,z\overline{z}}
+
\varphi_{2,zz\overline{z}\overline{z}}\,\varphi_{1,z\overline{z}}
}{
\varphi_{1,z\overline{z}}\,\varphi_{2,zz\overline{z}}
-
\varphi_{2,z\overline{z}}\,\varphi_{1,zz\overline{z}}}
\bigg)\cdot
\mathcal{S},
\endaligned
\]
while trivially (only in the rigid case!):
\[
\big[\mathcal{S},\mathcal{T}\big]
=
0.
\]

\medskip\noindent{\bf Symbolic treatment of the general case.}
Now, come back to the general case where $\varphi_1$, $\varphi_2$ do
depend on all variables $(x, y, u_1, u_2)$, so that, from the form of
the two defining equations $v_1 = \varphi^1$ and $v_2 = \varphi^2$:
\[
\aligned
\mathcal{L}
&
=
\frac{\partial}{\partial z}
+
\frac{\Lambda_1}{\Delta}\,\frac{\partial}{\partial u_1}
+
\frac{\Lambda_2}{\Delta}\,\frac{\partial}{\partial u_2},
\\
\overline{\mathcal{L}}
&
=
\frac{\partial}{\partial\overline{z}}
+
\frac{\overline{\Lambda}_1}{\overline{\Delta}}\,
\frac{\partial}{\partial u_1}
+
\frac{\overline{\Lambda}_2}{\overline{\Delta}}\,
\frac{\partial}{\partial u_2},
\\
\mathcal{T}
&
=
\ \ \ \ \ \ \ \ \ \ 
\frac{\Upsilon_1}{\Delta^2\,\overline{\Delta}^2}\,
\frac{\partial}{\partial u_1}
+
\frac{\Upsilon_2}{\Delta^2\,\overline{\Delta}^2}\,
\frac{\partial}{\partial u_2},
\\
\mathcal{S}
&
=
\ \ \ \ \ \ \ \ \ \ 
\frac{\Pi_1}{\Delta^4\,\overline{\Delta}^3}\,
\frac{\partial}{\partial u_1}
+
\frac{\Pi_2}{\Delta^4\,\overline{\Delta}^3}\,
\frac{\partial}{\partial u_2}.
\endaligned
\]

\medskip\noindent{\bf Lemma-Exercise.}
{\em
There is a uniquely defined function:}
\[
\ell
=
\ell\big(x,y,u_1,u_2\big)
\]
{\em defined on $M$ near $0$ such that the vector field:}
\[
\mathcal{N}
:=
\mathcal{L}
+
\ell\,\overline{\mathcal{L}}
\]
{\em has the property:}
\[
\big[\mathcal{N},\,\mathcal{T}\big]
\,\equiv\,
0
\ \
\mod\,
\big(
\mathcal{L},\,\overline{\mathcal{L}},\,\mathcal{T}
\big).
\]

\proof
The solution appears implicitly below.
\endproof

Since, as always, $\mathcal{ T}$ is real:
\[
\overline{\mathcal{T}}
=
\mathcal{T},
\]
one sees:
\[
\overline{\mathcal{S}}
=
\overline{
\big[\mathcal{L},\mathcal{T}\big]}
=
\big[\overline{\mathcal{L}},\mathcal{T}\big].
\]
But of course:
\[
\aligned
\big[
\overline{\mathcal{L}},\,
\mathcal{T}
\big]
&
=
\function\cdot\frac{\partial}{\partial u_1}
+
\function\cdot\frac{\partial}{\partial u_2}
\\
&
=
\function\cdot\mathcal{T}
+
\function\cdot\mathcal{S}.
\endaligned
\]

Therefore,
by introducing two {\em fundamental}
functions $A$ and $B$ of $(x,y, u_1, u_2)$ so that:
\[
\aligned
\big[
\overline{\mathcal{L}},\,
\mathcal{T}
\big]
&
=
A\,\mathcal{T}
+
B\,\mathcal{S}
\\
&
=
\overline{\mathcal{S}},
\endaligned
\] 
these two functions\,\,---\,\,that will
be very fundamental when applying Cartan's method later 
on\,\,---\,\,enjoy some relations coming from
the comparison with:
\[
\aligned
\big[\mathcal{L},\,\mathcal{T}\big]
&
=
\overline{
\big[\overline{\mathcal{L}},\,\mathcal{T}\big]}
\\
&
=
\overline{A}\,\mathcal{T}
+
\overline{B}\,
\underbrace{\overline{\mathcal{S}}}_{\sf replace}
\\
&
=
\overline{A}\,\mathcal{T}
+
\overline{B}\big(
A\,\mathcal{T}
+
B\,\mathcal{S}\big)
\\
&
=
\big(\overline{A}+\overline{B}\,A\big)\,\mathcal{T}
+
B\overline{B}\,\mathcal{S}
\\
&
=
\mathcal{S}
\ \ \ \ \ \ \ \ \ \ \ \ \ \ \ \ \ \ \ \ \ \ \ \ \ \ \ \ \ \ \ \ \ \ 
\ \ \ \ \ \
\explain{by definition},
\endaligned
\]
whence by identification:
\[
\boxed{
\aligned
B\overline{B}
&
\equiv
1,
\\
\overline{A}
+
\overline{B}A
&
\equiv
0.
\endaligned
}
\]
Later in the Cartan procedure, one should take account of
these two relations.

For now:
\[
\aligned
\big[\mathcal{L}+\ell\overline{\mathcal{L}},\,\mathcal{T}\big]
&
\,\equiv\,
\big[\mathcal{L},\mathcal{T}\big]
+
\ell\,
\big[\overline{\mathcal{L}},\mathcal{T}\big]
\ \ \
\mod\,
\big(
\mathcal{L},\,\overline{\mathcal{L}},\,\mathcal{T}
\big),
\\
&
\,\equiv\,
\mathcal{S}
+
\ell\,\big(A\,\mathcal{T}+B\,\mathcal{S}\big)
\ \
\mod\,
\big(
\mathcal{L},\,\overline{\mathcal{L}},\,\mathcal{T}
\big)
\\
&
\,\equiv\,
\big(1+\ell\,B\big)\,
\mathcal{S},
\endaligned
\]
hence the solution to the lemma-exercise is:
\[
\ell
=
-\,\frac{1}{B}.
\]

\medskip

Next, introduce two new {\em fundamental} functions $P$ and $Q$ of
$\big( x,y, u_1, u_2 \big)$ so that:
\[
\aligned
\big[\mathcal{L},\,\mathcal{S}\big]
&
=
P\,\mathcal{T}
+
Q\,\mathcal{S},
\endaligned
\]
whence {\em passim} after conjugation and replacement:
\[
\aligned
\big[\overline{\mathcal{L}},\,\overline{\mathcal{S}}\big]
&
=
\overline{P}\,\mathcal{T}
+
\overline{Q}\,\overline{\mathcal{S}}
\\
&
=
\overline{P}\,\mathcal{T}
+
\overline{Q}\,\big(
A\,\mathcal{T}+B\,\mathcal{S}
\big)
\\
&
=
\big(\overline{P}+A\,\overline{Q}\big)\,
\mathcal{T}
+
\big(B\,\overline{Q}\big)\,\mathcal{S}.
\endaligned
\]

\medskip\noindent{\bf Lemma.}
{\em One has the reality:}
\[
\big[\overline{\mathcal{L}},\,\mathcal{S}\big]
=
\big[\mathcal{L},\,\overline{\mathcal{S}}\big].
\]

\proof
Indeed, in full intermediate details:
\[
\aligned
\big[\overline{\mathcal{L}},\,\mathcal{S}\big]
&
=
\big[\overline{\mathcal{L}},\,
\big[\mathcal{L},\,\isqrt\,[\mathcal{L},\overline{\mathcal{L}}]\big]\big]
\\
&
\!\!\!\!
\overset{\sf Jacobi}{\,=\,}
-\,\zero{\big[[\mathcal{L},\overline{\mathcal{L}}],\,
\big[\overline{\mathcal{L}},\,\isqrt\,\mathcal{L}\big]\big]}
-
\big[
\mathcal{L},\,\big[[\mathcal{L},\overline{\mathcal{L}}],\,
\isqrt\,\overline{\mathcal{L}}\big]
\big]
\\
&
=
\big[
\mathcal{L},\,\big[\overline{\mathcal{L}},\,
\isqrt\,[\mathcal{L},\overline{\mathcal{L}}]\big]\big]
\\
&
=
\big[\mathcal{L},\,\overline{\mathcal{S}}\big],
\endaligned
\]
considering that the $\isqrt$ factor does not perturb
Jacobi's identity.
\endproof

As an application, determine:
\[
\aligned
\big[\overline{\mathcal{L}},\,\mathcal{S}\big]
&
=
\big[\mathcal{L},\,\overline{\mathcal{S}}\big]
\\
&
=
\big[\mathcal{L},\,A\mathcal{T}+B\mathcal{S}\big]
\\
&
=
\mathcal{L}(A)\cdot\mathcal{T}
+
\mathcal{L}(B)\cdot\mathcal{S}
+
A\,
\underbrace{\big[\mathcal{L},\mathcal{T}\big]}_{
\mathcal{S}}
+
B\,
\underbrace{\big[\mathcal{L},\mathcal{S}\big]}_{
P\mathcal{T}+Q\mathcal{S}}
\\
&
=
\big(
\mathcal{L}(A)+BP\big)\cdot
\mathcal{T}
+
\big(\mathcal{L}(B)+BQ+A\big)\cdot
\mathcal{S}.
\endaligned
\]

\medskip\noindent{\bf Scholium.}
{\em The reality condition:}
\[
\big[\overline{\mathcal{L}},\,\mathcal{S}\big]
=
\big[\mathcal{L},\,\overline{\mathcal{S}}\big]
\]
{\em entails:}
\[
\boxed{
\aligned
\mathcal{L}(A)
+
BP
&
\equiv
\overline{\mathcal{L}}\big(\overline{A}\big)
+
A\,\overline{\mathcal{L}}\big(\overline{B}\big)
+
\overline{B}\,\overline{P}
+
A\,\overline{B}\,\overline{Q}
+
A\,\overline{A},\,\,
\\
\,\,
\mathcal{L}(B)+BQ+A
&
\equiv
B\,\overline{\mathcal{L}}\big(\overline{B}\big)
+
B\,\overline{B}\,\overline{Q}
+
\overline{A}\,B.
\endaligned}
\]

\proof
Starting from what has been obtained at the moment:
\[
\aligned
\big[\overline{\mathcal{L}},\,\mathcal{S}\big]
&
=
\big(
\mathcal{L}(A)+BP\big)\cdot
\mathcal{T}
+
\big(\mathcal{L}(B)+BQ+A\big)\cdot
\mathcal{S},
\endaligned
\]
a plain conjugation gives:
\[
\aligned
\big[\mathcal{L},\,\overline{\mathcal{S}}\big]
&
=
\big(
\overline{\mathcal{L}}(\overline{A})
+
\overline{B}\,\overline{P}
\big)\cdot\mathcal{T}
+
\big(
\overline{\mathcal{L}}\big(\overline{B}\big)
+
\overline{B}\,\overline{Q}
+
\overline{A}
\big)
\underbrace{\overline{\mathcal{S}}}_{
A\mathcal{T}+B\mathcal{S}}
\\
&
=
\Big(
\overline{\mathcal{L}}\big(\overline{A}\big)
+
A\,\overline{\mathcal{L}}\big(\overline{B}\big)
+
\overline{B}\,\overline{P}
+
A\,\overline{B}\,\overline{Q}
+
A\,\overline{A}
\Big)\cdot\mathcal{T}
+
\\
&
\ \ \ \ \ 
+
\Big(
B\,\overline{\mathcal{L}}\big(\overline{B}\big)
+
B\,\overline{B}\,\overline{Q}
+
\overline{A}\,B
\Big)\cdot\mathcal{S},
\endaligned
\]
and an identification gives the result.
\endproof

Lastly, compute the 6\textsuperscript{th} Lie bracket
$\big[ \mathcal{S}, \, \mathcal{ T} \big]$.

\medskip\noindent{\bf Lemma.}
{\em One has the Jacobi-type Lie bracket relation:}
\[
\aligned
\big[\mathcal{S},\,\mathcal{T}\big]
&
=
\isqrt\,\big[\overline{\mathcal{L}},\,
[\mathcal{L},\mathcal{S}]
\big]
-
\isqrt\,
\big[
\mathcal{L},\,
\big[\overline{\mathcal{L}},\,\mathcal{S}\big]
\big].
\endaligned
\]

\proof
Indeed, in full intermediate details:
\[
\aligned
\big[\mathcal{S},\,\mathcal{T}\big]
&
=
\big[\mathcal{S},\,\isqrt\,[\mathcal{L},\overline{\mathcal{L}}\big]\big]
\\
&
\,=\,
\isqrt\,
\big[\mathcal{S},\,[\mathcal{L},\overline{\mathcal{L}}\big]\big]
\\
&\!\!
\overset{\text{\sf Jacobi}}{\,=\,}
-\,\isqrt\,\big[\overline{\mathcal{L}},\,[\mathcal{S},\mathcal{L}]\big]
-
\isqrt\,\big[\mathcal{L},\,\big[\overline{\mathcal{L}},\mathcal{S}]\big]
\\
&
\,=\,
\isqrt\,\big[\overline{\mathcal{L}},\,
[\mathcal{L},\mathcal{S}].
\big]
-
\isqrt\,
\big[
\mathcal{L},\,
\big[\overline{\mathcal{L}},\,\mathcal{S}\big]
\big],
\endaligned
\]
which is so.
\endproof

Consequently:
\[
\!\!\!\!\!\!\!\!\!\!\!\!\!\!\!\!\!\!\!\!
\aligned
\big[\mathcal{S},\mathcal{T}\big]
&
=
\isqrt\,\big[\overline{\mathcal{L}},\,\,P\mathcal{T}+Q\mathcal{S}\big]
-
\isqrt\,\big[\mathcal{L},\,\,\big(\mathcal{L}(A)+BP\big)\mathcal{T}
+
\big(\mathcal{L}(B)+BQ+A\big)\mathcal{S}
\big]
\\
&
=
\isqrt\,\overline{\mathcal{L}}(P)\cdot\mathcal{T}
+
\isqrt\,\overline{\mathcal{L}}(Q)\cdot\mathcal{S}
+
\isqrt\,
\underbrace{P\,\big[\overline{\mathcal{L}},\mathcal{T}\big]}_{
A\mathcal{T}+B\mathcal{S}}
+
\isqrt\,Q\,
\underbrace{\,\,\,\big[\overline{\mathcal{L}},\mathcal{S}\big]\,\,\,}_{
(\mathcal{L}(A)+BP)\mathcal{T}
+
\atop
+
(\mathcal{L}(B)+BQ+A)\mathcal{S}
}
-
\\
&
\ \ \ \ \
-\,
\isqrt\,\big(\mathcal{L}(\mathcal{L}(A))
+
P\,\mathcal{L}(B)
+
B\,\mathcal{L}(P)\big)\cdot\mathcal{T}
-
\\
&
\ \ \ \ \
-
\isqrt\,\big(\mathcal{L}(\mathcal{L}(B))
+
Q\,\mathcal{L}(B)
+
B\,\mathcal{L}(Q)
+
\mathcal{L}(A)
\big)\cdot\mathcal{S}
-
\\
&
-\,
\isqrt\,\big(\mathcal{L}(A)+B\,P\big)\,
\underbrace{\big[\mathcal{L},\mathcal{T}\big]}_{\mathcal{S}}
-
\isqrt\,\big(\mathcal{L}(B)+BQ+A\big)\,
\underbrace{\big[\mathcal{L},\mathcal{S}\big]}_{
P\mathcal{T}+Q\mathcal{S}}
\\
&
=
\Big(
-\,
\isqrt\,\mathcal{L}(\mathcal{L}(A))
-
\isqrt\,P\,\mathcal{L}(B)
+
\isqrt\,\overline{\mathcal{L}}(P)
+
\isqrt\,Q\,\mathcal{L}(A)
-
\\
&
\ \ \ \ \ \ \ \ \ \ \ \ \ \ \ \ \ \ \ \ \ \ \ \ \ \ \ \ \ \
\ \ \ \ \ \ \ \ \ \ \ \ \ \ \ \ \ \ \ \ \ \ \ \ \ \ \ \ \ \
-
\isqrt\,P\,\mathcal{L}(B)
-
\isqrt\,B\,\mathcal{L}(P)
\Big)\cdot\mathcal{T}
+
\\
&
\ \ \ \ \
+
\Big(
-\,\isqrt\,\mathcal{L}(\mathcal{L}(B))
-
\isqrt\,Q\,\mathcal{L}(B)
-
\isqrt\,B\,\mathcal{L}(Q)
-
\\
&
\ \ \ \ \ \ \ \ \ \ \ \ \ \ \ \ \ \ \ \ \ \ \ \ \ \ \ \ \ \
\ \ \ \ \ \ \ \ \ \ \ \ \ \ \ \ \ \ \ \ \ \ \ \ \ \ \ \ \ \
-
2\isqrt\,\mathcal{L}(A)
+
\isqrt\,\overline{\mathcal{L}}(Q)
\Big)\cdot\mathcal{S}.
\endaligned
\]

\medskip

In order to abbreviate, it will be convenient to set:
\[
\!\!\!\!\!\!\!\!\!\!\!\!\!\!\!\!\!\!\!\!
\!\!\!\!\!\!\!\!\!\!\!\!\!\!\!\!\!\!\!\!
\boxed{\,\,
\aligned
{\sf E}_{\sf rpl}
&
:=
\mathcal{L}(A)
+
BP,
\\
{\sf F}_{\sf rpl}
&
:=
\mathcal{L}(B)+BQ+A,
\\
{\sf G}_{\sf rpl}
&
:=
\isqrt\,\mathcal{L}(\mathcal{L}(A))
+
\isqrt\,P\,\mathcal{L}(B)
-
\isqrt\,\overline{\mathcal{L}}(P)
-
\isqrt\,Q\,\mathcal{L}(A)
+
\isqrt\,P\,\mathcal{L}(B)
+
\isqrt\,B\,\mathcal{L}(P),\,\,
\\
{\sf H}_{\sf rpl}
&
:=
\isqrt\,\mathcal{L}(\mathcal{L}(B))
+
\isqrt\,Q\,\mathcal{L}(B)
+
\isqrt\,B\,\mathcal{L}(Q)
+
2\isqrt\,\mathcal{L}(A)
-
\isqrt\,\overline{\mathcal{L}}(Q),
\endaligned}
\]
the lower index "${}_{\sf rpl}$" meaning that one should replace these
letters by their expression in terms of the {\bf 4} fundamental
functions:
\[
\boxed{\,\,
A,\ \ \ \ \
B,\ \ \ \ \
P,\ \ \ \ \
Q.\,\,}
\]
In other words, {\em the four functions
${\sf E}_{\sf rpl}$, ${\sf F}_{\sf rpl}$, ${\sf G}_{\sf rpl}$, 
${\sf H}_{\sf rpl}$ are not fundamental}.

\medskip\noindent{\bf Summary.} One has the 6 Lie bracket relations:
\[
\boxed{\,\,
\aligned
\big[\mathcal{S},\,\mathcal{T}\big]
&
=
-\,{\sf G}_{\sf rpl}\cdot\mathcal{T}
-
{\sf H}_{\sf rpl}\cdot\mathcal{S},
\\
\big[\mathcal{S},\,\overline{\mathcal{L}}\big]
&
=
-\,{\sf E}_{\sf rpl}\cdot\mathcal{T}
-
{\sf F}_{\sf rpl}\cdot\mathcal{S},
\\
\big[\mathcal{S},\,\mathcal{L}\big]
&
=
-\,P\cdot\mathcal{T}
-
Q\,\mathcal{S},
\\
\big[\mathcal{T},\,\overline{\mathcal{L}}\big]
&
=
-\,A\cdot\mathcal{T}
-
B\,\mathcal{S},
\\
\big[\mathcal{T},\,\mathcal{L}\big]
&
=
\,-\,\mathcal{S},
\\
\big[\overline{\mathcal{L}},\,\mathcal{L}\big]
&
=
\isqrt\,\mathcal{T}.\,\,
\endaligned}
\]

\medskip\noindent{\bf Initial Darboux structure of the 
dual coframe.} The coframe:
\[
\big\{du_2,du_1,\,dz,\,d\overline{z}\big\}
\]
is clearly dual to the frame: 
\[
\Big\{ 
{\textstyle{\frac{\partial}{\partial u_2}}},\,
{\textstyle{\frac{\partial}{\partial u_1}}},\,
{\textstyle{\frac{\partial}{\partial z}}},\,
{\textstyle{\frac{\partial}{\partial\overline{z}}}} 
\Big\}.
\]
Introduce then the coframe:
\[
\big\{
\sigma_0,\,
\rho_0,\,
\overline{\zeta_0},\,
\zeta_0\big\}
\]
which is dual to the frame:
\[
\big\{
\mathcal{S},\,
\mathcal{T},\,
\overline{\mathcal{L}},\,
\mathcal{L}\big\},
\]
namely:
\[
\begin{array}{cccc}
\sigma_0(\mathcal{S})=1 \ \ \ & \ \ \ \sigma_0(\mathcal{T})=0 \ \ \ &
\ \ \ \sigma_0(\overline{\mathcal{L}})=0 \ \ \ & \ \ \
\sigma_0\big(\mathcal{L}\big)=0,
\\
\rho_0(\mathcal{S})=0
\ \ \ & \ \ \ \rho_0(\mathcal{T})=1 \ \ \ & \ \ \
\rho_0(\overline{\mathcal{L}})=0 \ \ \ & \ \ \ \rho_0(\mathcal{L})=0,
\\
\overline{\zeta_0}(\mathcal{S})=0 \ \ \ & \ \ \
\overline{\zeta_0}(\mathcal{T})=0 \ \ \ & \ \ \
\overline{\zeta_0}(\overline{\mathcal{L}})=1 \ \ \ & \ \ \
\overline{\zeta_0}\big(\mathcal{L}\big)=0,
\\
\zeta_0(\mathcal{S})=0 \ \ \ & \ \ \ \zeta_0(\mathcal{T})=0 \ \ \ & \
\ \ \zeta_0(\overline{\mathcal{L}})=0 \ \ \ & \ \ \
\zeta_0\big(\mathcal{L}\big)=1.
\end{array}
\]
One has (exercice):
\[
\zeta_0 
= 
dz 
\ \ \ \ \ \ \ \ \ \ \ \ \ 
\text{\rm and} 
\ \ \ \ \ \ \ \ \ \ \ \ \ 
\overline{\zeta_0} 
= 
d\overline{z}.
\]

Proceeding as in the launching file~\cite{ Merker-launching-2011}, 
and as in~\cite{ Merker-Sabzevari-2014b}, one
determines the Darboux structure of this initial coframe by
reading vertically the following convenient
auxiliary array, putting an overall minus
sign:
\[
\footnotesize
\begin{array}{cccccccccccc}
& & \mathcal{S} & & \mathcal{T} & & \overline{\mathcal{\mathcal{L}}} & &
\mathcal{L} 
\\
& & \boxed{d\sigma_0} & & \boxed{d\rho_0} & & \boxed{d\overline{\zeta}_0} & &
\boxed{d\zeta_0} & 
\\
\big[\mathcal{S},\,\mathcal{T}\big]
& = & -\,{\sf H}_{\sf rpl} & + & -\,{\sf G}_{\sf rpl} & + & 0 & + & 
0 &
\boxed{\sigma_0\wedge\rho_0}
\\
\big[\mathcal{S},\,\overline{\mathcal{L}}\big]
& = & -\,{\sf F}_{\sf rpl} & + & -\,{\sf E}_{\sf rpl}
& + & 0 & + & 0 & 
\boxed{\sigma_0\wedge\overline{\zeta}_0}
\\
\big[\mathcal{S},\,\mathcal{L}\big]
& = & -\,Q & + & -\,P & + & 0 & + & 0 & 
\boxed{\sigma_0\wedge\zeta_0}
\\
\big[\mathcal{T},\,\overline{\mathcal{L}}\big]
& = & -\,B & + & -\,A & + & 0 & + & 0 & 
\boxed{\rho_0\wedge\overline{\zeta_0}}
\\
\big[\mathcal{T},\,\mathcal{L}\big]
& = & -\,1 & + & 0 & + & 0 & + & 0 & 
\boxed{\rho_0\wedge\zeta_0}
\\
\big[\overline{\mathcal{L}},\,\mathcal{L}\big]
& = & 0 & + & \isqrt & + & 0 & + & 0 & 
\boxed{\overline{\zeta}_0\wedge\zeta_0}
\end{array}
\]
This gives:
\[
\boxed{\,\,
\aligned
d\sigma_0
&
=
{\sf H}_{\sf rpl}\cdot\sigma_0\wedge\rho_0
+
{\sf F}_{\sf rpl}\cdot\sigma_0\wedge\overline{\zeta}_0
+
Q\cdot\sigma_0\wedge\zeta_0
+
B\cdot\rho_0\wedge\overline{\zeta}_0
+
\rho_0\wedge\zeta_0,\,\,
\\
d\rho_0
&
=
{\sf G}_{\sf rpl}\cdot\sigma_0\wedge\rho_0
+
{\sf E}_{\sf rpl}\cdot\sigma_0\wedge\overline{\zeta}_0
+
P\cdot\sigma_0\wedge\zeta_0
+
A\cdot\rho_0\wedge\overline{\zeta}_0
+
\isqrt\,\zeta_0\wedge\overline{\zeta}_0,\,\,
\\
d\overline{\zeta}_0
&
=
0,
\\
d\zeta_0
&
=
0.\,\,
\endaligned}
\]


\bigskip

\section{\sf $M^5 \subset \C^4$ of general class 
$\text{\sf III}_{\text{\sf 1}}$: 
\\
initial frame and coframe in local coordinates}
\label{initial-III-1}
\HEAD{\ref{initial-III-1}.~$M^5 \subset \C^4$ of general class 
$\text{\sf III}_{\text{\sf 1}}$
initial frame and coframe in local coordinates}{
Jo\"el {\sc Merker}, D\'epartement de Math\'ematiques d'Orsay}

\medskip

Next, consider:
\[
\Big(
M^5
\,\subset\,
\C^4
\Big)
\,\,\in\,\,
\text{\sf General Class $\text{\sf III}_{\text{\sf 1}}$}.
\]

Representing as before
(\cite{ Merker-Pocchiola-Sabzevari-5-CR-II, Merker-5-CR-III})
$M$ in coordinates:
\[
(z,w_1,w_2,w_3) 
= 
\big(
x+\isqrt\,y,\, 
u_1+\isqrt\,v_1,\,
u_2+\isqrt\,v_2,\,
u_3+\isqrt\,v_3\big),
\]
as a graph:
\[
\aligned
v_1
&
=
\varphi_1(x,y,u_1,u_2,u_3),
\\
v_2
&
=
\varphi_2(x,y,u_1,u_2,u_3),
\\
v_3
&
=
\varphi_3(x,y,u_1,u_2,u_3),
\endaligned
\]
an associated explicit frame for $TM$ is:
\[
\!\!\!\!\!\!\!\!\!\!\!\!\!\!\!\!\!\!\!\!
\aligned
X
&
=
\frac{\partial}{\partial x}
+
\varphi_{1,x}(x,y,u_1,u_2,u_3)\,\frac{\partial}{\partial v_1}
+
\varphi_{2,x}(x,y,u_1,u_2,u_3)\,\frac{\partial}{\partial v_2}
+
\varphi_{3,x}(x,y,u_1,u_2,u_3)\,\frac{\partial}{\partial v_3},
\\
Y
&
=
\frac{\partial}{\partial y}
+
\varphi_{1,y}(x,y,u_1,u_2,u_3)\,\frac{\partial}{\partial v_1}
+
\varphi_{2,y}(x,y,u_1,u_2,u_3)\,\frac{\partial}{\partial v_2}
+
\varphi_{3,y}(x,y,u_1,u_2,u_3)\,\frac{\partial}{\partial v_3},
\\
U_1
&
=
\frac{\partial}{\partial u_1}
+
\varphi_{1,u_1}(x,y,u_1,u_2,u_3)\,\frac{\partial}{\partial v_1}
+
\varphi_{2,u_1}(x,y,u_1,u_2,u_3)\,\frac{\partial}{\partial v_2}
+
\varphi_{3,u_1}(x,y,u_1,u_2,u_3)\,\frac{\partial}{\partial v_3},
\\
U_2
&
=
\frac{\partial}{\partial u_2}
+
\varphi_{1,u_2}(x,y,u_1,u_2,u_3)\,\frac{\partial}{\partial v_1}
+
\varphi_{2,u_2}(x,y,u_1,u_2,u_3)\,\frac{\partial}{\partial v_2}
+
\varphi_{3,u_2}(x,y,u_1,u_2,u_3)\,\frac{\partial}{\partial v_3},
\\
U_3
&
=
\frac{\partial}{\partial u_3}
+
\varphi_{1,u_3}(x,y,u_1,u_2,u_3)\,\frac{\partial}{\partial v_1}
+
\varphi_{2,u_3}(x,y,u_1,u_2,u_3)\,\frac{\partial}{\partial v_2}
+
\varphi_{3,u_3}(x,y,u_1,u_2,u_3)\,\frac{\partial}{\partial v_3}.
\endaligned
\]
The projection:
\[
\pi\colon
\ \ \ \ \ \ \ \ \ \ \ \ \ \ \
\aligned
\R^8
&
\,\longrightarrow\,
\R^5
\\
(x,y,u_1,v_1,u_2,v_2,u_3,v_3)
&
\,\longmapsto\,
(x,y,u_1,u_2,u_3),
\endaligned
\]
makes a {\sl chart} on $M$ and does:
\[
\aligned
\pi_*
\bigg(
\frac{\partial}{\partial x}
+
\varphi_{1,x}\,
\frac{\partial}{\partial v_1}
+
\varphi_{2,x}\,
\frac{\partial}{\partial v_2}
+
\varphi_{3,x}\,
\frac{\partial}{\partial v_3}
\bigg)
&
=
\frac{\partial}{\partial x},
\\
\pi_*
\bigg(
\frac{\partial}{\partial y}
+
\varphi_{1,y}\,
\frac{\partial}{\partial v_1}
+
\varphi_{2,y}\,
\frac{\partial}{\partial v_2}
+
\varphi_{3,y}\,
\frac{\partial}{\partial v_3}
\bigg)
&
=
\frac{\partial}{\partial y},
\\
\pi_*
\bigg(
\frac{\partial}{\partial u_1}
+
\varphi_{1,u_1}\,
\frac{\partial}{\partial v_1}
+
\varphi_{2,u_1}\,
\frac{\partial}{\partial v_2}
+
\varphi_{3,u_1}\,
\frac{\partial}{\partial v_3}
\bigg)
&
=
\frac{\partial}{\partial u_1},
\\
\pi_*
\bigg(
\frac{\partial}{\partial u_2}
+
\varphi_{1,u_2}\,
\frac{\partial}{\partial v_1}
+
\varphi_{2,u_2}\,
\frac{\partial}{\partial v_2}
+
\varphi_{3,u_2}\,
\frac{\partial}{\partial v_3}
\bigg)
&
=
\frac{\partial}{\partial u_2},
\\
\pi_*
\bigg(
\frac{\partial}{\partial u_3}
+
\varphi_{1,u_3}\,
\frac{\partial}{\partial v_1}
+
\varphi_{2,u_3}\,
\frac{\partial}{\partial v_2}
+
\varphi_{3,u_3}\,
\frac{\partial}{\partial v_3}
\bigg)
&
=
\frac{\partial}{\partial u_3}.
\endaligned
\]

Next, a natural {\em intrinsic} local vector field generator
for $T^{1, 0} M$ is:
\[ 
\mathcal{L} 
:= 
\frac{\partial}{\partial z} 
+ 
A_1\,\frac{\partial}{\partial u_1}
+ 
A_2\,\frac{\partial}{\partial u_2}
+ 
A_3\,\frac{\partial}{\partial u_3}, 
\] 
where:
\[
\aligned
A_1
&
:=
\frac{
\left\vert\!
\begin{array}{ccc}
-\,\varphi_{1,z} & \varphi_{1,u_2} & \varphi_{1,u_3}\\
-\,\varphi_{2,z} & \isqrt+\varphi_{2,u_2} & \varphi_{2,u_3}\\
-\,\varphi_{3,z} & \varphi_{3,u_2} & \isqrt+\varphi_{3,u_3}
\end{array}
\!\right\vert
}{
\left\vert\!
\begin{array}{ccc}
\isqrt+\varphi_{1,u_1} & \varphi_{1,u_2} & \varphi_{1,u_3}\\
\varphi_{2,u_1} & \isqrt+\varphi_{2,u_2} & \varphi_{2,u_3}\\
\varphi_{3,u_1} & \varphi_{3,u_2} & \isqrt+\varphi_{3,u_3}
\end{array}
\!\right\vert
},
\endaligned
\]
\[
\aligned
A_2
&
:=
\frac{
\left\vert\!
\begin{array}{ccc}
\isqrt+ \varphi_{1,u_1} & -\,\varphi_{1,z} & \varphi_{1,u_3}\\
\varphi_{2,u_1} & -\,\varphi_{2,z} &  \varphi_{2,u_3}\\
\varphi_{3,u_1} & -\,\varphi_{3,z} & \isqrt+\varphi_{3,u_3}
\end{array}
\!\right\vert
}{
\left\vert\!
\begin{array}{ccc}
\isqrt+\varphi_{1,u_1} & \varphi_{1,u_2} & \varphi_{1,u_3}\\
\varphi_{2,u_1} & \isqrt+\varphi_{2,u_2} & \varphi_{2,u_3}\\
\varphi_{3,u_1} & \varphi_{3,u_2} & \isqrt+\varphi_{3,u_3}
\end{array}
\!\right\vert
},
\endaligned
\]

\[
\aligned
A_3
&
:=
\frac{
\left\vert\!
\begin{array}{ccc}
\isqrt+\varphi_{1,u_1} & \varphi_{1,u_2} & -\,\varphi_{1,z}\\
\varphi_{2,u_1} &  \varphi_{2,u_2} & -\,\varphi_{2,z}\\
\varphi_{3,u_1}& \isqrt+\varphi_{3,u_2} & -\,\varphi_{3,z} 
\end{array}
\!\right\vert
}{
\left\vert\!
\begin{array}{ccc}
\isqrt+\varphi_{1,u_1} & \varphi_{1,u_2} & \varphi_{1,u_3}\\
\varphi_{2,u_1} & \isqrt+\varphi_{2,u_2} & \varphi_{2,u_3}\\
\varphi_{3,u_1} & \varphi_{3,u_2} & \isqrt+\varphi_{3,u_3}
\end{array}
\!\right\vert
}.
\endaligned
\]

Set:
\[
\Delta
:=
\left\vert\!
\begin{array}{ccc}
\isqrt+\varphi_{1,u_1} & \varphi_{1,u_2} & \varphi_{1,u_3}\\
\varphi_{2,u_1} & \isqrt+\varphi_{2,u_2} & \varphi_{2,u_3}\\
\varphi_{3,u_1} & \varphi_{3,u_2} & \isqrt+\varphi_{3,u_3}
\end{array}
\!\right\vert,
\]
whence:
\[
\overline{\Delta}
:=
\left\vert\!
\begin{array}{ccc}
-\isqrt+\varphi_{1,u_1} & \varphi_{1,u_2} & \varphi_{1,u_3}\\
\varphi_{2,u_1} & -\isqrt+\varphi_{2,u_2} & \varphi_{2,u_3}\\
\varphi_{3,u_1} & \varphi_{3,u_2} & -\isqrt+\varphi_{3,u_3}
\end{array}
\!\right\vert.
\]

Also, set:
\[
\aligned
\Lambda_1
&
:=
\left\vert\!
\begin{array}{ccc}
-\,\varphi_{1,z} & \varphi_{1,u_2} & \varphi_{1,u_3}\\
-\,\varphi_{2,z} & \isqrt+\varphi_{2,u_2} & \varphi_{2,u_3}\\
-\,\varphi_{3,z} & \varphi_{3,u_2} & \isqrt+\varphi_{3,u_3}
\end{array}
\!\right\vert,
\\
\Lambda_2
&
:=
\left\vert\!
\begin{array}{ccc}
\isqrt+ \varphi_{1,u_1} & -\,\varphi_{1,z} & \varphi_{1,u_3}\\
\varphi_{2,u_1} & -\,\varphi_{2,z} &  \varphi_{2,u_3}\\
\varphi_{3,u_1} & -\,\varphi_{3,z} & \isqrt+\varphi_{3,u_3}
\end{array}
\!\right\vert,
\\
\Lambda_3
&
:=
\left\vert\!
\begin{array}{ccc}
\isqrt+\varphi_{1,u_1} & \varphi_{1,u_2} & -\,\varphi_{1,z}\\
\varphi_{2,u_1} &  \varphi_{2,u_2} & -\,\varphi_{2,z}\\
\varphi_{3,u_1}& \isqrt+\varphi_{3,u_2} & -\,\varphi_{3,z} 
\end{array}
\!\right\vert,
\endaligned
\]
whence:
\[
\aligned
\overline{\Lambda}_1
&
:=
\left\vert\!
\begin{array}{ccc}
-\,\varphi_{1,\overline{z}} & \varphi_{1,u_2} & \varphi_{1,u_3}\\
-\,\varphi_{2,\overline{z}} & -\isqrt+\varphi_{2,u_2} & \varphi_{2,u_3}\\
-\,\varphi_{3,\overline{z}} & \varphi_{3,u_2} & -\isqrt+\varphi_{3,u_3}
\end{array}
\!\right\vert,
\\
\overline{\Lambda}_2
&
:=
\left\vert\!
\begin{array}{ccc}
-\isqrt+ \varphi_{1,u_1} & -\,\varphi_{1,\overline{z}} & \varphi_{1,u_3}\\
\varphi_{2,u_1} & -\,\varphi_{2,\overline{z}} &  \varphi_{2,u_3}\\
\varphi_{3,u_1} & -\,\varphi_{3,\overline{z}} & -\isqrt+\varphi_{3,u_3}
\end{array}
\!\right\vert,
\\
\overline{\Lambda}_3
&
:=
\left\vert\!
\begin{array}{ccc}
-\isqrt+\varphi_{1,u_1} & \varphi_{1,u_2} & -\,\varphi_{1,\overline{z}}\\
\varphi_{2,u_1} &  \varphi_{2,u_2} & -\,\varphi_{2,\overline{z}}\\
\varphi_{3,u_1}& -\isqrt+\varphi_{3,u_2} & -\,\varphi_{3,\overline{z}} 
\end{array}
\!\right\vert.
\endaligned
\]

In these notations:
\[
\aligned
\mathcal{L}
&
=
\frac{\partial}{\partial z}
+
\frac{\Lambda_1}{\Delta}\,\frac{\partial}{\partial u_1}
+
\frac{\Lambda_2}{\Delta}\,\frac{\partial}{\partial u_2}
+
\frac{\Lambda_3}{\Delta}\,\frac{\partial}{\partial u_3},
\\
\overline{\mathcal{L}}
&
=
\frac{\partial}{\partial\overline{z}}
+
\frac{\overline{\Lambda}_1}{\overline{\Delta}}\,
\frac{\partial}{\partial u_1}
+
\frac{\overline{\Lambda}_2}{\overline{\Delta}}\,
\frac{\partial}{\partial u_2}
+
\frac{\overline{\Lambda}_3}{\overline{\Delta}}\,
\frac{\partial}{\partial u_3}.
\endaligned
\]

\medskip

Set:
\[
\mathcal{T}
:=
\isqrt\,\big[\mathcal{L},\overline{\mathcal{L}}\big],
\]
and set:
\[
\mathcal{S}
:=
\big[\mathcal{L},\mathcal{T}\big],
\]
whence:
\[
\overline{\mathcal{S}}
=
\big[\overline{\mathcal{L}},\mathcal{T}\big].
\]

By hypothesis, the CR-geometric invariant condition:
\[
\aligned
\C\otimes_\R TM
=
T^{1,0}M+T^{0,1}M
+
[T^{1,0}M,\,T^{0,1}M]
&
+
\big[T^{1,0}M,\,[T^{1,0}M,\,T^{0,1}M]\big]
\\
&
+
\big[T^{0,1}M,\,[T^{1,0}M,\,T^{0,1}M]\big]
\endaligned
\]
holds, which means as is known that the $5$ fields:
\[
\big\{
\mathcal{L},\,\overline{\mathcal{L}},\,
\mathcal{T},\,
\mathcal{S},\,\overline{\mathcal{S}}
\big\}
\]
constitute a frame for:
\[
{\bf 5}
=
\rank_\C\big(
\C\otimes_\R TM
\big).
\]

\medskip\noindent{\bf Lemma.}
{\em There are certain uniquely defined 
coefficient-functions that are polynomials:}
\[
\aligned
\Upsilon_1
&
=
\Upsilon_1
\Big(
\varphi_{1,x^jy^ku_1^{l_1}u_2^{l_2}u_3^{l_3}},\,\,
\varphi_{2,x^jy^ku_1^{l_1}u_2^{l_2}u_3^{l_3}},\,\,
\varphi_{3,x^jy^ku_1^{l_1}u_2^{l_2}u_3^{l_3}}
\Big)_{1\leqslant j+k+l_1+l_2+l_3\leqslant 2},
\\
\Upsilon_2
&
=
\Upsilon_2
\Big(
\varphi_{1,x^jy^ku_1^{l_1}u_2^{l_2}u_3^{l_3}},\,\,
\varphi_{2,x^jy^ku_1^{l_1}u_2^{l_2}u_3^{l_3}},\,\,
\varphi_{3,x^jy^ku_1^{l_1}u_2^{l_2}u_3^{l_3}}
\Big)_{1\leqslant j+k+l_1+l_2+l_3\leqslant 2},
\\
\Upsilon_3
&
=
\Upsilon_3
\Big(
\varphi_{1,x^jy^ku_1^{l_1}u_2^{l_2}u_3^{l_3}},\,\,
\varphi_{2,x^jy^ku_1^{l_1}u_2^{l_2}u_3^{l_3}},\,\,
\varphi_{3,x^jy^ku_1^{l_1}u_2^{l_2}u_3^{l_3}}
\Big)_{1\leqslant j+k+l_1+l_2+l_3\leqslant 2},
\endaligned
\]
\[
\aligned
\Pi_1
&
=
\Pi_1
\Big(
\varphi_{1,x^jy^ku_1^{l_1}u_2^{l_2}u_3^{l_3}},\,\,
\varphi_{2,x^jy^ku_1^{l_1}u_2^{l_2}u_3^{l_3}},\,\,
\varphi_{3,x^jy^ku_1^{l_1}u_2^{l_2}u_3^{l_3}}
\Big)_{1\leqslant j+k+l_1+l_2+l_3\leqslant 3},
\\
\Pi_2
&
=
\Pi_2
\Big(
\varphi_{1,x^jy^ku_1^{l_1}u_2^{l_2}u_3^{l_3}},\,\,
\varphi_{2,x^jy^ku_1^{l_1}u_2^{l_2}u_3^{l_3}},\,\,
\varphi_{3,x^jy^ku_1^{l_1}u_2^{l_2}u_3^{l_3}}
\Big)_{1\leqslant j+k+l_1+l_2+l_3\leqslant 3},
\\
\Pi_3
&
=
\Pi_3
\Big(
\varphi_{1,x^jy^ku_1^{l_1}u_2^{l_2}u_3^{l_3}},\,\,
\varphi_{2,x^jy^ku_1^{l_1}u_2^{l_2}u_3^{l_3}},\,\,
\varphi_{3,x^jy^ku_1^{l_1}u_2^{l_2}u_3^{l_3}}
\Big)_{1\leqslant j+k+l_1+l_2+l_3\leqslant 3},
\endaligned
\] 
{\em such that:}
\[
\aligned
\mathcal{T}
&
=
\frac{\Upsilon_1}{\Delta^2\,\overline{\Delta}^2}\,
\frac{\partial}{\partial u_1}
+
\frac{\Upsilon_2}{\Delta^2\,\overline{\Delta}^2}\,
\frac{\partial}{\partial u_2}
+
\frac{\Upsilon_3}{\Delta^2\,\overline{\Delta}^2}\,
\frac{\partial}{\partial u_3}
\\
&
=:
Y_1\,
\frac{\partial}{\partial u_1}
+
Y_2\,
\frac{\partial}{\partial u_2}
+
Y_3\,
\frac{\partial}{\partial u_3},
\\
\mathcal{S}
&
=
\frac{\Pi_1}{\Delta^4\,\overline{\Delta}^3}\,
\frac{\partial}{\partial u_1}
+
\frac{\Pi_2}{\Delta^4\,\overline{\Delta}^3}\,
\frac{\partial}{\partial u_2}
+
\frac{\Pi_3}{\Delta^4\,\overline{\Delta}^3}\,
\frac{\partial}{\partial u_3}
\\
&
=:
H_1\,
\frac{\partial}{\partial u_1}
+
H_2\,
\frac{\partial}{\partial u_2}
+
H_3\,
\frac{\partial}{\partial u_3},
\\
\overline{\mathcal{S}}
&
=
\frac{\overline{\Pi}_1}{\Delta^3\,\overline{\Delta}^4}\,
\frac{\partial}{\partial u_1}
+
\frac{\overline{\Pi}_2}{\Delta^3\,\overline{\Delta}^4}\,
\frac{\partial}{\partial u_2}
+
\frac{\overline{\Pi}_3}{\Delta^3\,\overline{\Delta}^4}\,
\frac{\partial}{\partial u_3}
\\
&
=\,
\overline{H}_1\,
\frac{\partial}{\partial u_1}
+
\overline{H}_2\,
\frac{\partial}{\partial u_2}
+
\overline{H}_3\,
\frac{\partial}{\partial u_3}.
\endaligned
\]

\proof
The proof goes as for the preceding General Class $\text{\sf II}$, 
and it provides firstly:
\[
\boxed{\,
\aligned
\Upsilon_1
&
:=
\isqrt\,\Big(
\Delta\,\Delta\,\overline{\Delta}\,
\overline{\Lambda}_{1,z}
+
\Delta\,\overline{\Delta}\,\Lambda_1\,
\overline{\Lambda}_{1,u_1}
+
\Delta\,\overline{\Delta}\,\Lambda_2\,
\overline{\Lambda}_{1,u_2}
+
\Delta\,\overline{\Delta}\,\Lambda_3\,
\overline{\Lambda}_{1,u_3}
\,-\,\,
\\
&
\ \ \ \ \ \ \ \ \ \ \ \ \ \ \ 
-\,
\Delta\,\Delta\,\overline{\Delta}_z\,
\overline{\Lambda}_1
-
\Delta\,\Lambda_1\,\overline{\Delta}_{u_1}\,
\overline{\Lambda}_1
-
\Delta\,\Lambda_2\,\overline{\Delta}_{u_2}\,
\overline{\Lambda}_1
-
\Delta\,\Lambda_3\,\overline{\Delta}_{u_3}\,
\overline{\Lambda}_1
\,-
\\
&
\ \ \ \ \ \ \ \ \ \ \ \ \ \ \ 
-\,
\Delta\,\overline{\Delta}\,\overline{\Delta}\,
\Lambda_{1,\overline{z}}
-
\Delta\,\overline{\Delta}\,\overline{\Lambda}_1\,
\Lambda_{1,u_1}
-
\Delta\,\overline{\Delta}\,\overline{\Lambda}_2\,
\Lambda_{1,u_2}
-
\Delta\,\overline{\Delta}\,\overline{\Lambda}_3\,
\Lambda_{1,u_3}
+
\\
&
\ \ \ \ \ \ \ \ \ \ \ \ \ \ \ 
+
\overline{\Delta}\,\overline{\Delta}\,\Delta_{\overline{z}}\,
\Lambda_1
+
\overline{\Delta}\,\overline{\Lambda}_1\,\Delta_{u_1}\,
\Lambda_1
+
\overline{\Delta}\,\overline{\Lambda}_2\,\Delta_{u_2}\,
\Lambda_1
+
\overline{\Delta}\,\overline{\Lambda}_3\,\Delta_{u_3}\,
\Lambda_1
\Big),
\endaligned
\,}
\]
\[
\boxed{\,
\aligned
\Upsilon_2
&
:=
\isqrt\,\Big(
\Delta\,\Delta\,\overline{\Delta}\,
\overline{\Lambda}_{2,z}
+
\Delta\,\overline{\Delta}\,\Lambda_1\,
\overline{\Lambda}_{2,u_1}
+
\Delta\,\overline{\Delta}\,\Lambda_2\,
\overline{\Lambda}_{2,u_2}
+
\Delta\,\overline{\Delta}\,\Lambda_3\,
\overline{\Lambda}_{2,u_3}
\,-\,\,
\\
&
\ \ \ \ \ \ \ \ \ \ \ \ \ \ \ 
-\,
\Delta\,\Delta\,\overline{\Delta}_z\,
\overline{\Lambda}_2
-
\Delta\,\Lambda_1\,\overline{\Delta}_{u_1}\,
\overline{\Lambda}_2
-
\Delta\,\Lambda_2\,\overline{\Delta}_{u_2}\,
\overline{\Lambda}_2
-
\Delta\,\Lambda_3\,\overline{\Delta}_{u_3}\,
\overline{\Lambda}_2
\,-
\\
&
\ \ \ \ \ \ \ \ \ \ \ \ \ \ \ 
-\,
\Delta\,\overline{\Delta}\,\overline{\Delta}\,
\Lambda_{2,\overline{z}}
-
\Delta\,\overline{\Delta}\,\overline{\Lambda}_1\,
\Lambda_{2,u_1}
-
\Delta\,\overline{\Delta}\,\overline{\Lambda}_2\,
\Lambda_{2,u_2}
-
\Delta\,\overline{\Delta}\,\overline{\Lambda}_3\,
\Lambda_{2,u_3}
+
\\
&
\ \ \ \ \ \ \ \ \ \ \ \ \ \ \ 
+
\overline{\Delta}\,\overline{\Delta}\,\Delta_{\overline{z}}\,
\Lambda_2
+
\overline{\Delta}\,\overline{\Lambda}_1\,\Delta_{u_1}\,
\Lambda_2
+
\overline{\Delta}\,\overline{\Lambda}_2\,\Delta_{u_2}\,
\Lambda_2
+
\overline{\Delta}\,\overline{\Lambda}_3\,\Delta_{u_3}\,
\Lambda_2
\Big),
\endaligned
\,}
\]
\[
\boxed{\,
\aligned
\Upsilon_3
&
:=
\isqrt\,\Big(
\Delta\,\Delta\,\overline{\Delta}\,
\overline{\Lambda}_{3,z}
+
\Delta\,\overline{\Delta}\,\Lambda_1\,
\overline{\Lambda}_{3,u_1}
+
\Delta\,\overline{\Delta}\,\Lambda_2\,
\overline{\Lambda}_{3,u_2}
+
\Delta\,\overline{\Delta}\,\Lambda_3\,
\overline{\Lambda}_{3,u_3}
\,-\,\,
\\
&
\ \ \ \ \ \ \ \ \ \ \ \ \ \ \ 
-\,
\Delta\,\Delta\,\overline{\Delta}_z\,
\overline{\Lambda}_3
-
\Delta\,\Lambda_1\,\overline{\Delta}_{u_1}\,
\overline{\Lambda}_3
-
\Delta\,\Lambda_2\,\overline{\Delta}_{u_2}\,
\overline{\Lambda}_3
-
\Delta\,\Lambda_3\,\overline{\Delta}_{u_3}\,
\overline{\Lambda}_3
\,-
\\
&
\ \ \ \ \ \ \ \ \ \ \ \ \ \ \ 
-\,
\Delta\,\overline{\Delta}\,\overline{\Delta}\,
\Lambda_{3,\overline{z}}
-
\Delta\,\overline{\Delta}\,\overline{\Lambda}_1\,
\Lambda_{3,u_1}
-
\Delta\,\overline{\Delta}\,\overline{\Lambda}_2\,
\Lambda_{3,u_2}
-
\Delta\,\overline{\Delta}\,\overline{\Lambda}_3\,
\Lambda_{3,u_3}
+
\\
&
\ \ \ \ \ \ \ \ \ \ \ \ \ \ \ 
+
\overline{\Delta}\,\overline{\Delta}\,\Delta_{\overline{z}}\,
\Lambda_3
+
\overline{\Delta}\,\overline{\Lambda}_1\,\Delta_{u_1}\,
\Lambda_3
+
\overline{\Delta}\,\overline{\Lambda}_2\,\Delta_{u_2}\,
\Lambda_3
+
\overline{\Delta}\,\overline{\Lambda}_3\,\Delta_{u_3}\,
\Lambda_3
\Big),
\endaligned
\,}
\]
and secondly:
\[
\boxed{\,
\aligned
\Pi_1
&
=
\Delta\,\Delta\,\overline{\Delta}\,
\Upsilon_{1,z}
+
\Delta\,\overline{\Delta}\,\Lambda_1\,
\Upsilon_{1,u_1}
+
\Delta\,\overline{\Delta}\,\Lambda_2\,
\Upsilon_{1,u_2}
+
\Delta\,\overline{\Delta}\,\Lambda_3\,
\Upsilon_{1,u_3}
\,-
\\
&
\,-
2\,\Delta\,\overline{\Delta}\,\Delta_z\,
\Upsilon_1
-
2\,\overline{\Delta}\,\Lambda_1\,\Delta_{u_1}\,
\Upsilon_1
-
2\,\overline{\Delta}\,\Lambda_2\,\Delta_{u_2}\,
\Upsilon_1
-
2\,\overline{\Delta}\,\Lambda_3\,\Delta_{u_3}\,
\Upsilon_1
\,-\,\,
\\
&
\,-
2\,\Delta\,\Delta\,\overline{\Delta}_z\,
\Upsilon_1
-
2\,\Delta\,\Lambda_1\,\overline{\Delta}_{u_1}\,
\Upsilon_1
-
2\,\Delta\,\Lambda_2\,\overline{\Delta}_{u_2}\,
\Upsilon_1
-
2\,\Delta\,\Lambda_3\,\overline{\Delta}_{u_3}\,
\Upsilon_1
\,-
\\
&
\ \ \ \ \ \ \ \ \ \ \ \ \ \ \ \ \ \ \ \ \ \ \ \ \ \ \
-\,
\Delta\,\overline{\Delta}\,\Upsilon_1\,
\Lambda_{1,u_1}
-
\Delta\,\overline{\Delta}\,\Upsilon_2\,
\Lambda_{1,u_2}
-
\Delta\,\overline{\Delta}\,\Upsilon_3\,
\Lambda_{1,u_3}
+
\\
&
\ \ \ \ \ \ \ \ \ \ \ \ \ \ \ \ \ \ \ \ \ \ \ \ \ \ \
+
\overline{\Delta}\,\Upsilon_1\,\Delta_{u_1}\,
\Lambda_1
+
\overline{\Delta}\,\Upsilon_2\,\Delta_{u_2}\,
\Lambda_1
+
\overline{\Delta}\,\Upsilon_3\,\Delta_{u_3}\,
\Lambda_1,
\endaligned
}
\]
\[
\boxed{\,
\aligned
\Pi_2
&
=
\Delta\,\Delta\,\overline{\Delta}\,
\Upsilon_{2,z}
+
\Delta\,\overline{\Delta}\,\Lambda_1\,
\Upsilon_{2,u_1}
+
\Delta\,\overline{\Delta}\,\Lambda_2\,
\Upsilon_{2,u_2}
+
\Delta\,\overline{\Delta}\,\Lambda_3\,
\Upsilon_{2,u_3}
\,-
\\
&
\,-
2\,\Delta\,\overline{\Delta}\,\Delta_z\,
\Upsilon_2
-
2\,\overline{\Delta}\,\Lambda_1\,\Delta_{u_1}\,
\Upsilon_2
-
2\,\overline{\Delta}\,\Lambda_2\,\Delta_{u_2}\,
\Upsilon_2
-
2\,\overline{\Delta}\,\Lambda_3\,\Delta_{u_3}\,
\Upsilon_2
\,-\,\,
\\
&
\,-
2\,\Delta\,\Delta\,\overline{\Delta}_z\,
\Upsilon_2
-
2\,\Delta\,\Lambda_1\,\overline{\Delta}_{u_1}\,
\Upsilon_2
-
2\,\Delta\,\Lambda_2\,\overline{\Delta}_{u_2}\,
\Upsilon_2
-
2\,\Delta\,\Lambda_3\,\overline{\Delta}_{u_3}\,
\Upsilon_2
\,-
\\
&
\ \ \ \ \ \ \ \ \ \ \ \ \ \ \ \ \ \ \ \ \ \ \ \ \ \ \
-\,
\Delta\,\overline{\Delta}\,\Upsilon_1\,
\Lambda_{2,u_1}
-
\Delta\,\overline{\Delta}\,\Upsilon_2\,
\Lambda_{2,u_2}
-
\Delta\,\overline{\Delta}\,\Upsilon_3\,
\Lambda_{2,u_3}
+
\\
&
\ \ \ \ \ \ \ \ \ \ \ \ \ \ \ \ \ \ \ \ \ \ \ \ \ \ \
+
\overline{\Delta}\,\Upsilon_1\,\Delta_{u_1}\,
\Lambda_2
+
\overline{\Delta}\,\Upsilon_2\,\Delta_{u_2}\,
\Lambda_2
+
\overline{\Delta}\,\Upsilon_3\,\Delta_{u_3}\,
\Lambda_2,
\endaligned
}
\]
\[
\boxed{\,
\aligned
\Pi_3
&
=
\Delta\,\Delta\,\overline{\Delta}\,
\Upsilon_{3,z}
+
\Delta\,\overline{\Delta}\,\Lambda_1\,
\Upsilon_{3,u_1}
+
\Delta\,\overline{\Delta}\,\Lambda_2\,
\Upsilon_{3,u_2}
+
\Delta\,\overline{\Delta}\,\Lambda_3\,
\Upsilon_{3,u_3}
\,-
\\
&
\,-
2\,\Delta\,\overline{\Delta}\,\Delta_z\,
\Upsilon_3
-
2\,\overline{\Delta}\,\Lambda_1\,\Delta_{u_1}\,
\Upsilon_3
-
2\,\overline{\Delta}\,\Lambda_2\,\Delta_{u_2}\,
\Upsilon_3
-
2\,\overline{\Delta}\,\Lambda_3\,\Delta_{u_3}\,
\Upsilon_3
\,-\,\,
\\
&
\,-
2\,\Delta\,\Delta\,\overline{\Delta}_z\,
\Upsilon_3
-
2\,\Delta\,\Lambda_1\,\overline{\Delta}_{u_1}\,
\Upsilon_3
-
2\,\Delta\,\Lambda_2\,\overline{\Delta}_{u_2}\,
\Upsilon_2
-
2\,\Delta\,\Lambda_3\,\overline{\Delta}_{u_3}\,
\Upsilon_3
\,-
\\
&
\ \ \ \ \ \ \ \ \ \ \ \ \ \ \ \ \ \ \ \ \ \ \ \ \ \ \
-\,
\Delta\,\overline{\Delta}\,\Upsilon_1\,
\Lambda_{3,u_1}
-
\Delta\,\overline{\Delta}\,\Upsilon_2\,
\Lambda_{3,u_2}
-
\Delta\,\overline{\Delta}\,\Upsilon_3\,
\Lambda_{3,u_3}
+
\\
&
\ \ \ \ \ \ \ \ \ \ \ \ \ \ \ \ \ \ \ \ \ \ \ \ \ \ \
+
\overline{\Delta}\,\Upsilon_1\,\Delta_{u_1}\,
\Lambda_3
+
\overline{\Delta}\,\Upsilon_2\,\Delta_{u_2}\,
\Lambda_3
+
\overline{\Delta}\,\Upsilon_3\,\Delta_{u_3}\,
\Lambda_3,
\endaligned
}
\]
which concludes.
\endproof

\noindent{\bf Explicitness obstacle.}
{\em After reduction to a common minimal denominator, the three 
numerators in:}
\[
Y_1
=
\frac{\Upsilon_1}{
\Delta^2\,\overline{\Delta}^2},
\ \ \ \ \ \ \ \ \ \ \ \ \ \ \ \ \ \ \ \ \ \ \ \ \ 
Y_2
=
\frac{\Upsilon_2}{
\Delta^2\,\overline{\Delta}^2},
\ \ \ \ \ \ \ \ \ \ \ \ \ \ \ \ \ \ \ \ \ \ \ \ \ 
Y_3
=
\frac{\Upsilon_3}{
\Delta^2\,\overline{\Delta}^2},
\]
{\em are both polynomials in the ${\bf 3} \cdot {\bf 20}$ 
partial derivatives:}
\[
\Big(
\varphi_{1,x^jy^ku_1^{l_1}u_2^{l_2}u_3^{l_3}},\,\,
\varphi_{2,x^jy^ku_1^{l_1}u_2^{l_2}u_3^{l_3}},\,\,
\varphi_{3,x^jy^ku_1^{l_1}u_2^{l_2}u_3^{l_3}}
\Big)_{1\leqslant j+k+l_1+l_2+l_3\leqslant 2}
\]
{\em incorporating:}
\[
{\bf 41\,964}
\]
{\em monomials, while the next three numerators in:}
\[
H_1
=
\frac{\Pi_1}{
\Delta^4\,\overline{\Delta}^3},
\ \ \ \ \ \ \ \ \ \ \ \ \ \ \ \ \ \ \ \ \ \ \ \ \ 
H_2
=
\frac{\Pi_2}{
\Delta^4\,\overline{\Delta}^3},
\ \ \ \ \ \ \ \ \ \ \ \ \ \ \ \ \ \ \ \ \ \ \ \ \ 
H_3
=
\frac{\Pi_3}{
\Delta^4\,\overline{\Delta}^3},
\]
{\em would incorporate approximately more than:}
\[
{\bf 100\,000\,000}
\]
{\em monomials in the ${\bf 3} \cdot {\bf 55}$ partial derivatives:}
\[
\Big(
\varphi_{1,x^jy^ku_1^{l_1}u_2^{l_2}u_3^{l_3}},\,\,
\varphi_{2,x^jy^ku_1^{l_1}u_2^{l_2}u_3^{l_3}},\,\,
\varphi_{3,x^jy^ku_1^{l_1}u_2^{l_2}u_3^{l_3}}
\Big)_{1\leqslant j+k+l_1+l_2+l_3\leqslant 3},
\]
{\em since {\em e.g.} $\Pi_1$ includes:}
\[
\Delta\,\Delta\,\overline{\Delta}\,\Upsilon_{1,z},
\]
{\em inside which $\Delta\, \Delta\, \overline{ \Delta}$ incorporates:}
\[
{\bf 526}
\]
{\em monomials while $\Upsilon_{1, z}$ incorporates:}
\[
{\bf 236\,068}
\]
{\em monomials.}

\proof
The product is really unwieldy, even on a powerful computer software,
because the monomials themselves are almost all of an already large
size, looking like {\em e.g.} the five following ones:
\[
\aligned
&
-\,\big(\varphi_{3,u_1}\big)^2\,\varphi_{1,u_3}\big)^2\,
\big(\varphi_{2,u_2}\big)^2\,\varphi_{2,\overline{z}u}\,
\varphi_{3,u_2}\,\varphi_{2,z}\,
\varphi_{1,u_2},
\\
&
-\,\isqrt\,
\big(\varphi_{2,u_2}\big)^2\,
\varphi_{3,u_1}\,\varphi_{1,u_2}\,\varphi_{2,\overline{z}u_3}\,
\varphi_{1,z}\,\varphi_{3,u_3},
\\
&
\isqrt\,\big(\varphi_{1,u_1}\big)^2\,
\big(\varphi_{2,u_2}\big)^3\,
\big(\varphi_{3,u_3}\big)^2\,
\varphi_{2,u_1}\,
\varphi_{1,u_3}\,
\varphi_{3,u_2}\,\varphi_{1,zz},
\\
&
6\,\isqrt\,\big(\varphi_{1,u_1}\big)^2\,
\varphi_{2,u_3}\,\varphi_{3,u_2}\,\varphi_{1,z}\,\varphi_{3,u_3}\,
\varphi_{3,zu_1}\,
\varphi_{1,u_3}\,\varphi_{2,u_2},
\\
&
-\,4\,\isqrt\,
\varphi_{1,u_1}\,\varphi_{2,u_1}\,\varphi_{2,u_1}\,
\big(\varphi_{1,u_2}\big)^2\,
\varphi_{2,u_3}\,\varphi_{3,z}\,\varphi_{2,zu_2}\,
\varphi_{3,u_3},
\endaligned
\]
so that complete explicitnes\,\,---\,\,even at the basic
level of the coefficients of the initial frame!\,\,---\,\,seems 
to be out of reach.
\endproof

The Class $\text{\sf III}_{\text{\sf 1}}$ 
hypothesis now reads as the nonzeroness:
\[
0
\,\neq\,
\det\,
\left(\!
\begin{array}{ccc}
\frac{\Upsilon_1}{\Delta^2\overline{\Delta}^2} &
\frac{\Upsilon_2}{\Delta^2\overline{\Delta}^2} &
\frac{\Upsilon_3}{\Delta^2\overline{\Delta}^2}
\medskip
\\
\frac{\Pi_1}{\Delta^4\overline{\Delta}^3} &
\frac{\Pi_2}{\Delta^4\overline{\Delta}^3}&
\frac{\Pi_3}{\Delta^4\overline{\Delta}^3}\medskip
\\
\frac{\overline{\Pi}_1}{\Delta^3\overline{\Delta}^4} &
\frac{\overline{\Pi}_2}{\Delta^3\overline{\Delta}^4}&
\frac{\overline{\Pi}_3}{\Delta^3\overline{\Delta}^4}
\end{array}
\!\right)
(x,y,u_1,u_2,u_3),
\]
at every point.

More precisely, if one preliminarily normalizes
coordinates as in~\cite{ Merker-5-CR-III}:
\[
\aligned
v_1
&
=
z\overline{z}
\ \ \ \ \ \ \ \ \ \ \ \ \ \ \ \ \ \ \ \ \ \ \
+
z\overline{z}\,{\rm O}_2\big(z,\overline{z}\big)
+
z\overline{z}\,{\rm O}_1(u_1)
+
z\overline{z}\,{\rm O}_1(u_2)
+
z\overline{z}\,{\rm O}_1(u_3),\,\,
\\
v_2
&
=
z^2\overline{z}
+
z\overline{z}^2
\ \ \ \ \ \ \ \ \ \ \,
+
z\overline{z}\,{\rm O}_2\big(z,\overline{z}\big)
+
z\overline{z}\,{\rm O}_1(u_1)
+
z\overline{z}\,{\rm O}_1(u_2)
+
z\overline{z}\,{\rm O}_1(u_3),\,\,
\\
v_3
&
=
\isqrt\,\big(z^2\overline{z}
-
z\overline{z}^2\big)
+
z\overline{z}\,{\rm O}_2\big(z,\overline{z}\big)
+
z\overline{z}\,{\rm O}_1(u_1)
+
z\overline{z}\,{\rm O}_1(u_2)
+
z\overline{z}\,{\rm O}_1(u_3),
\endaligned
\]
so that at the origin:
\[
\aligned
\mathcal{L}\big\vert_0
&
=
\frac{\partial}{\partial z}
\bigg\vert_0,
\\
\overline{\mathcal{L}}\big\vert_0
&
=
\frac{\partial}{\partial\overline{z}}
\bigg\vert_0,
\\
\big[\mathcal{L},\overline{\mathcal{L}}\big]
\Big\vert_0
&
=
-\,2\,\isqrt\,
\frac{\partial}{\partial u_1}
\bigg\vert_0,
\\
\big[\mathcal{L},\,\big[\mathcal{L},\overline{\mathcal{L}}\big]\big]
\Big\vert_0
&
=
-\,4\,\isqrt\,
\frac{\partial}{\partial u_2}
\bigg\vert_0
+
4\,
\frac{\partial}{\partial u_3}
\bigg\vert_0,
\\
\big[\overline{\mathcal{L}},\,
\big[\mathcal{L},\overline{\mathcal{L}}\big]\big]
\Big\vert_0
&
=
-\,4\,\isqrt\,
\frac{\partial}{\partial u_2}
\bigg\vert_0
-
4\,
\frac{\partial}{\partial u_3}
\bigg\vert_0,
\endaligned
\]
the determinant in question:
\[
\det\,
\left(\!
\begin{array}{ccc}
\frac{\Upsilon_1}{\Delta^2\overline{\Delta}^2} &
\frac{\Upsilon_2}{\Delta^2\overline{\Delta}^2} &
\frac{\Upsilon_3}{\Delta^2\overline{\Delta}^2}
\medskip
\\
\frac{\Pi_1}{\Delta^4\overline{\Delta}^3} &
\frac{\Pi_2}{\Delta^4\overline{\Delta}^3}&
\frac{\Pi_3}{\Delta^4\overline{\Delta}^3}\medskip
\\
\frac{\overline{\Pi}_1}{\Delta^3\overline{\Delta}^4} &
\frac{\overline{\Pi}_2}{\Delta^3\overline{\Delta}^4}&
\frac{\overline{\Pi}_3}{\Delta^3\overline{\Delta}^4}
\end{array}
\!\right)
(0)
\,=\,
\left\vert\!
\begin{array}{ccc}
2 & 0 & 0
\\
0 & 4 & 4\isqrt
\\
0 & 4 & -4\isqrt
\end{array}
\!\right\vert
\]
becomes quite visibly nonzero, hence also near the origin.

\medskip

But generally, the disease is, that when (necessarily) re-expressing:
\[
\aligned
\frac{\partial}{\partial u_1}
&
=
\frac{\Delta^2\,\overline{\Delta}^2\,
\left\vert\!
\begin{array}{cc}
\Pi_2 & \Pi_3\\
\overline{\Pi}_2 & \overline{\Pi}_3
\end{array}
\!\right\vert
}{
\left\vert\!
\begin{array}{ccc}
\Upsilon_1 & \Upsilon_2 & \Upsilon_3\\
\Pi_1 & \Pi_2 & \Pi_3\\
\overline{\Pi}_1 & \overline{\Pi}_2 & \overline{\Pi}_3
\end{array}
\!\right\vert}
\cdot
\mathcal{T}
-
\frac{\Delta^4\,\overline{\Delta}^3\,
\left\vert\!
\begin{array}{cc}
\Upsilon_2 & \Upsilon_3\\
\overline{\Pi}_2 & \overline{\Pi}_3
\end{array}
\!\right\vert
}{
\left\vert\!
\begin{array}{ccc}
\Upsilon_1 & \Upsilon_2 & \Upsilon_3\\
\Pi_1 & \Pi_2 & \Pi_3\\
\overline{\Pi}_1 & \overline{\Pi}_2 & \overline{\Pi}_3
\end{array}
\!\right\vert}
\cdot
\mathcal{S}
+
\frac{\Delta^4\,\overline{\Delta}^3\,
\left\vert\!
\begin{array}{cc}
\Upsilon_2 & \Upsilon_3\\
\Pi_2 & \Pi_3
\end{array}
\!\right\vert
}{
\left\vert\!
\begin{array}{ccc}
\Upsilon_1 & \Upsilon_2 & \Upsilon_3\\
\Pi_1 & \Pi_2 & \Pi_3\\
\overline{\Pi}_1 & \overline{\Pi}_2 & \overline{\Pi}_3
\end{array}
\!\right\vert}
\cdot
\overline{\mathcal{S}},
\endaligned
\]
\[
\aligned
\frac{\partial}{\partial u_2}
&
=
\frac{\Delta^2\,\overline{\Delta}^2\,
\left\vert\!
\begin{array}{cc}
\Pi_3 & \Pi_1\\
\overline{\Pi}_3 & \overline{\Pi}_1
\end{array}
\!\right\vert
}{
\left\vert\!
\begin{array}{ccc}
\Upsilon_1 & \Upsilon_2 & \Upsilon_3\\
\Pi_1 & \Pi_2 & \Pi_3\\
\overline{\Pi}_1 & \overline{\Pi}_2 & \overline{\Pi}_3
\end{array}
\!\right\vert}
\cdot
\mathcal{T}
-
\frac{\Delta^4\,\overline{\Delta}^3\,
\left\vert\!
\begin{array}{cc}
\Upsilon_3 & \Upsilon_1\\
\overline{\Pi}_3 & \overline{\Pi}_1
\end{array}
\!\right\vert
}{
\left\vert\!
\begin{array}{ccc}
\Upsilon_1 & \Upsilon_2 & \Upsilon_3\\
\Pi_1 & \Pi_2 & \Pi_3\\
\overline{\Pi}_1 & \overline{\Pi}_2 & \overline{\Pi}_3
\end{array}
\!\right\vert}
\cdot
\mathcal{S}
+
\frac{\Delta^4\,\overline{\Delta}^3\,
\left\vert\!
\begin{array}{cc}
\Upsilon_3 & \Upsilon_1\\
\Pi_3 & \Pi_1
\end{array}
\!\right\vert
}{
\left\vert\!
\begin{array}{ccc}
\Upsilon_1 & \Upsilon_2 & \Upsilon_3\\
\Pi_1 & \Pi_2 & \Pi_3\\
\overline{\Pi}_1 & \overline{\Pi}_2 & \overline{\Pi}_3
\end{array}
\!\right\vert}
\cdot
\overline{\mathcal{S}},
\endaligned
\]
\[
\aligned
\frac{\partial}{\partial u_3}
&
=
\frac{\Delta^2\,\overline{\Delta}^2\,
\left\vert\!
\begin{array}{cc}
\Pi_1 & \Pi_2\\
\overline{\Pi}_1 & \overline{\Pi}_2
\end{array}
\!\right\vert
}{
\left\vert\!
\begin{array}{ccc}
\Upsilon_1 & \Upsilon_2 & \Upsilon_3\\
\Pi_1 & \Pi_2 & \Pi_3\\
\overline{\Pi}_1 & \overline{\Pi}_2 & \overline{\Pi}_3
\end{array}
\!\right\vert}
\cdot
\mathcal{T}
-
\frac{\Delta^4\,\overline{\Delta}^3\,
\left\vert\!
\begin{array}{cc}
\Upsilon_1 & \Upsilon_2\\
\overline{\Pi}_1 & \overline{\Pi}_2
\end{array}
\!\right\vert
}{
\left\vert\!
\begin{array}{ccc}
\Upsilon_1 & \Upsilon_2 & \Upsilon_3\\
\Pi_1 & \Pi_2 & \Pi_3\\
\overline{\Pi}_1 & \overline{\Pi}_2 & \overline{\Pi}_3
\end{array}
\!\right\vert}
\cdot
\mathcal{S}
+
\frac{\Delta^4\,\overline{\Delta}^3\,
\left\vert\!
\begin{array}{cc}
\Upsilon_1 & \Upsilon_2\\
\Pi_1 & \Pi_2
\end{array}
\!\right\vert
}{
\left\vert\!
\begin{array}{ccc}
\Upsilon_1 & \Upsilon_2 & \Upsilon_3\\
\Pi_1 & \Pi_2 & \Pi_3\\
\overline{\Pi}_1 & \overline{\Pi}_2 & \overline{\Pi}_3
\end{array}
\!\right\vert}
\cdot
\overline{\mathcal{S}},
\endaligned
\]
some quite huge fractions appear, and it becomes
(very, very) impossible for a powerful computer, to make explicit
the coefficients in the next $7$ brackets completing
the Lie structure of the frame:
\[
\aligned
\big[\mathcal{L},\mathcal{S}\big]
&
=
\function
\cdot
\mathcal{T}
+
\function
\cdot
\mathcal{S}
+
\function
\cdot
\overline{\mathcal{S}},
\\
\big[\mathcal{L},\overline{\mathcal{S}}\big]
&
=
\function
\cdot
\mathcal{T}
+
\function
\cdot
\mathcal{S}
+
\function
\cdot
\overline{\mathcal{S}},
\\
\big[\overline{\mathcal{L}},\mathcal{S}\big]
&
=
\function
\cdot
\mathcal{T}
+
\function
\cdot
\mathcal{S}
+
\function
\cdot
\overline{\mathcal{S}},
\\
\big[\overline{\mathcal{L}},\overline{\mathcal{S}}\big]
&
=
\function
\cdot
\mathcal{T}
+
\function
\cdot
\mathcal{S}
+
\function
\cdot
\overline{\mathcal{S}},
\\
\big[\mathcal{T},\mathcal{S}\big]
&
=
\function
\cdot
\mathcal{T}
+
\function
\cdot
\mathcal{S}
+
\function
\cdot
\overline{\mathcal{S}},
\\
\big[\mathcal{T},\overline{\mathcal{S}}\big]
&
=
\function
\cdot
\mathcal{T}
+
\function
\cdot
\mathcal{S}
+
\function
\cdot
\overline{\mathcal{S}},
\\
\big[\mathcal{S},\overline{\mathcal{S}}\big]
&
=
\function
\cdot
\mathcal{T}
+
\function
\cdot
\mathcal{S}
+
\function
\cdot
\overline{\mathcal{S}}.
\endaligned
\]

\medskip\noindent{\bf Rigid collapse.}
In the (special) so-called {\sl rigid} case:
\[
\aligned
v_1
&
=
\varphi_1(x,y),
\\
v_2
&
=
\varphi_2(x,y),
\\
v_2
&
=
\varphi_3(x,y),
\endaligned
\]
where the graphing functions do not depend upon
$u_1$, $u_2$, $u_3$, completely explicit formulas can be typesetted.

Indeed:
\[
\aligned
\mathcal{L}
&
=
\frac{\partial}{\partial z}
+
\isqrt\,\varphi_{1,z}\,
\frac{\partial}{\partial u_1}
+
\isqrt\,\varphi_{2,z}\,
\frac{\partial}{\partial u_2}
+
\isqrt\,\varphi_{3,z}\,
\frac{\partial}{\partial u_3},
\\
\overline{\mathcal{L}}
&
=
\frac{\partial}{\partial\overline{z}}
-
\isqrt\,\varphi_{1,\overline{z}}\,
\frac{\partial}{\partial u_1}
-
\isqrt\,\varphi_{2,\overline{z}}\,
\frac{\partial}{\partial u_2}
-
\isqrt\,\varphi_{3,\overline{z}}\,
\frac{\partial}{\partial u_3},
\\
\mathcal{T}
&
=
\ \ \ \ \ \ \ \ \ \ \ \ \,
2\,\varphi_{1,z\overline{z}}\,
\frac{\partial}{\partial u_1}
+
2\,\varphi_{2,z\overline{z}}\,
\frac{\partial}{\partial u_2}
+
2\,\varphi_{3,z\overline{z}}\,
\frac{\partial}{\partial u_3},
\\
\mathcal{S}
&
=
\ \ \ \ \ \ \ \ \ \ \ \ \,
2\,\varphi_{1,zz\overline{z}}\,
\frac{\partial}{\partial u_1}
+
2\,\varphi_{2,zz\overline{z}}\,
\frac{\partial}{\partial u_2}
+
2\,\varphi_{3,zz\overline{z}}\,
\frac{\partial}{\partial u_3},
\\
\overline{\mathcal{S}}
&
=
\ \ \ \ \ \ \ \ \ \ \ \ \,
2\,\varphi_{1,z\overline{z}\overline{z}}\,
\frac{\partial}{\partial u_1}
+
2\,\varphi_{2,z\overline{z}\overline{z}}\,
\frac{\partial}{\partial u_2}
+
2\,\varphi_{3,z\overline{z}\overline{z}}\,
\frac{\partial}{\partial u_3}.
\endaligned
\]

The relevant determinant is nowhere vanishing:
\[
8\,
\left\vert\!
\begin{array}{ccc}
\varphi_{z\overline{z}}^1 & \varphi_{z\overline{z}}^2 &
\varphi_{z\overline{z}}^3
\\
\varphi_{zz\overline{z}}^1 & \varphi_{zz\overline{z}}^2 &
\varphi_{zz\overline{z}}^3
\\
\varphi_{z\overline{z}\overline{z}}^1 & 
\varphi_{z\overline{z}\overline{z}}^2 &
\varphi_{z\overline{z}\overline{z}}^3
\end{array}
\right\vert
\,\neq\,
0,
\]
hence one solves:
\[
\!\!\!\!\!\!\!\!\!\!\!\!\!\!\!\!\!\!\!\!
\aligned
\frac{\partial}{\partial u_1}
&
=
\frac{
\frac{1}{2}\,
\left\vert\!
\begin{array}{cc}
\varphi_{zz\overline{z}}^2 & 
\varphi_{zz\overline{z}}^3\\
\varphi_{z\overline{z}\overline{z}}^2 & 
\varphi_{z\overline{z}\overline{z}}^3
\end{array}
\!\right\vert
}{
\left\vert\!
\begin{array}{ccc}
\varphi_{z\overline{z}}^1 & \varphi_{z\overline{z}}^2 &
\varphi_{z\overline{z}}^3
\\
\varphi_{zz\overline{z}}^1 & \varphi_{zz\overline{z}}^2 &
\varphi_{zz\overline{z}}^3
\\
\varphi_{z\overline{z}\overline{z}}^1 & 
\varphi_{z\overline{z}\overline{z}}^2 &
\varphi_{z\overline{z}\overline{z}}^3
\end{array}
\right\vert
}
\cdot
\mathcal{T}
-
\frac{
\frac{1}{2}\,
\left\vert\!
\begin{array}{cc}
\varphi_{z\overline{z}}^2 & 
\varphi_{z\overline{z}}^3
\\
\varphi_{z\overline{z}\overline{z}}^2 & 
\varphi_{z\overline{z}\overline{z}}^3
\end{array}
\!\right\vert
}{
\left\vert\!
\begin{array}{ccc}
\varphi_{z\overline{z}}^1 & \varphi_{z\overline{z}}^2 &
\varphi_{z\overline{z}}^3
\\
\varphi_{zz\overline{z}}^1 & \varphi_{zz\overline{z}}^2 &
\varphi_{zz\overline{z}}^3
\\
\varphi_{z\overline{z}\overline{z}}^1 & 
\varphi_{z\overline{z}\overline{z}}^2 &
\varphi_{z\overline{z}\overline{z}}^3
\end{array}
\right\vert
}
\cdot
\mathcal{S}
+
\frac{
\frac{1}{2}\,
\left\vert\!
\begin{array}{cc}
\varphi_{z\overline{z}}^2 & 
\varphi_{z\overline{z}}^3
\\
\varphi_{zz\overline{z}}^2 & 
\varphi_{zz\overline{z}}^3
\end{array}
\!\right\vert
}{
\left\vert\!
\begin{array}{ccc}
\varphi_{z\overline{z}}^1 & \varphi_{z\overline{z}}^2 &
\varphi_{z\overline{z}}^3
\\
\varphi_{zz\overline{z}}^1 & \varphi_{zz\overline{z}}^2 &
\varphi_{zz\overline{z}}^3
\\
\varphi_{z\overline{z}\overline{z}}^1 & 
\varphi_{z\overline{z}\overline{z}}^2 &
\varphi_{z\overline{z}\overline{z}}^3
\end{array}
\right\vert
}
\cdot
\overline{\mathcal{S}},
\endaligned
\]
\[
\!\!\!\!\!\!\!\!\!\!\!\!\!\!\!\!\!\!\!\!
\aligned
\frac{\partial}{\partial u_2}
&
=
\frac{
\frac{1}{2}\,
\left\vert\!
\begin{array}{cc}
\varphi_{zz\overline{z}}^3 & 
\varphi_{zz\overline{z}}^1\\
\varphi_{z\overline{z}\overline{z}}^3 & 
\varphi_{z\overline{z}\overline{z}}^1
\end{array}
\!\right\vert
}{
\left\vert\!
\begin{array}{ccc}
\varphi_{z\overline{z}}^1 & \varphi_{z\overline{z}}^2 &
\varphi_{z\overline{z}}^3
\\
\varphi_{zz\overline{z}}^1 & \varphi_{zz\overline{z}}^2 &
\varphi_{zz\overline{z}}^3
\\
\varphi_{z\overline{z}\overline{z}}^1 & 
\varphi_{z\overline{z}\overline{z}}^2 &
\varphi_{z\overline{z}\overline{z}}^3
\end{array}
\right\vert
}
\cdot
\mathcal{T}
-
\frac{
\frac{1}{2}\,
\left\vert\!
\begin{array}{cc}
\varphi_{z\overline{z}}^1 & 
\varphi_{z\overline{z}}^3
\\
\varphi_{z\overline{z}\overline{z}}^1 & 
\varphi_{z\overline{z}\overline{z}}^3
\end{array}
\!\right\vert
}{
\left\vert\!
\begin{array}{ccc}
\varphi_{z\overline{z}}^1 & \varphi_{z\overline{z}}^2 &
\varphi_{z\overline{z}}^3
\\
\varphi_{zz\overline{z}}^1 & \varphi_{zz\overline{z}}^2 &
\varphi_{zz\overline{z}}^3
\\
\varphi_{z\overline{z}\overline{z}}^1 & 
\varphi_{z\overline{z}\overline{z}}^2 &
\varphi_{z\overline{z}\overline{z}}^3
\end{array}
\right\vert
}
\cdot
\mathcal{S}
+
\frac{
\frac{1}{2}\,
\left\vert\!
\begin{array}{cc}
\varphi_{z\overline{z}}^1 & 
\varphi_{z\overline{z}}^2
\\
\varphi_{zz\overline{z}}^1 & 
\varphi_{zz\overline{z}}^2
\end{array}
\!\right\vert
}{
\left\vert\!
\begin{array}{ccc}
\varphi_{z\overline{z}}^1 & \varphi_{z\overline{z}}^2 &
\varphi_{z\overline{z}}^3
\\
\varphi_{zz\overline{z}}^1 & \varphi_{zz\overline{z}}^2 &
\varphi_{zz\overline{z}}^3
\\
\varphi_{z\overline{z}\overline{z}}^1 & 
\varphi_{z\overline{z}\overline{z}}^2 &
\varphi_{z\overline{z}\overline{z}}^3
\end{array}
\right\vert
}
\cdot
\overline{\mathcal{S}},
\endaligned
\]
\[
\!\!\!\!\!\!\!\!\!\!\!\!\!\!\!\!\!\!\!\!
\aligned
\frac{\partial}{\partial u_3}
&
=
\frac{
\frac{1}{2}\,
\left\vert\!
\begin{array}{cc}
\varphi_{zz\overline{z}}^1 & 
\varphi_{zz\overline{z}}^2\\
\varphi_{z\overline{z}\overline{z}}^1 & 
\varphi_{z\overline{z}\overline{z}}^2
\end{array}
\!\right\vert
}{
\left\vert\!
\begin{array}{ccc}
\varphi_{z\overline{z}}^1 & \varphi_{z\overline{z}}^2 &
\varphi_{z\overline{z}}^3
\\
\varphi_{zz\overline{z}}^1 & \varphi_{zz\overline{z}}^2 &
\varphi_{zz\overline{z}}^3
\\
\varphi_{z\overline{z}\overline{z}}^1 & 
\varphi_{z\overline{z}\overline{z}}^2 &
\varphi_{z\overline{z}\overline{z}}^3
\end{array}
\right\vert
}
\cdot
\mathcal{T}
-
\frac{
\frac{1}{2}\,
\left\vert\!
\begin{array}{cc}
\varphi_{z\overline{z}}^1 & 
\varphi_{z\overline{z}}^2
\\
\varphi_{z\overline{z}\overline{z}}^1 & 
\varphi_{z\overline{z}\overline{z}}^2
\end{array}
\!\right\vert
}{
\left\vert\!
\begin{array}{ccc}
\varphi_{z\overline{z}}^1 & \varphi_{z\overline{z}}^2 &
\varphi_{z\overline{z}}^3
\\
\varphi_{zz\overline{z}}^1 & \varphi_{zz\overline{z}}^2 &
\varphi_{zz\overline{z}}^3
\\
\varphi_{z\overline{z}\overline{z}}^1 & 
\varphi_{z\overline{z}\overline{z}}^2 &
\varphi_{z\overline{z}\overline{z}}^3
\end{array}
\right\vert
}
\cdot
\mathcal{S}
+
\frac{
\frac{1}{2}\,
\left\vert\!
\begin{array}{cc}
\varphi_{z\overline{z}}^1 & 
\varphi_{z\overline{z}}^2
\\
\varphi_{zz\overline{z}}^1 & 
\varphi_{zz\overline{z}}^2
\end{array}
\!\right\vert
}{
\left\vert\!
\begin{array}{ccc}
\varphi_{z\overline{z}}^1 & \varphi_{z\overline{z}}^2 &
\varphi_{z\overline{z}}^3
\\
\varphi_{zz\overline{z}}^1 & \varphi_{zz\overline{z}}^2 &
\varphi_{zz\overline{z}}^3
\\
\varphi_{z\overline{z}\overline{z}}^1 & 
\varphi_{z\overline{z}\overline{z}}^2 &
\varphi_{z\overline{z}\overline{z}}^3
\end{array}
\right\vert
}
\cdot
\overline{\mathcal{S}}.
\endaligned
\]

Next:
\[
\!\!\!\!\!\!\!\!\!\!\!\!\!\!\!\!\!\!\!\!
\aligned
\big[\mathcal{L},\mathcal{S}\big]
&
=
2\,\varphi_{1,zzz\overline{z}}\,
\frac{\partial}{\partial u_1}
+
2\,\varphi_{2,zzz\overline{z}}\,
\frac{\partial}{\partial u_2}
+
2\,\varphi_{3,zzz\overline{z}}\,
\frac{\partial}{\partial u_3}
\\
&
=
\left(
\frac{
\varphi_{1,zzz\overline{z}}\,
\left\vert\!
\begin{array}{cc}
\varphi_{zz\overline{z}}^2 & 
\varphi_{zz\overline{z}}^3\\
\varphi_{z\overline{z}\overline{z}}^2 & 
\varphi_{z\overline{z}\overline{z}}^3
\end{array}
\!\right\vert
}{
\left\vert\!
\begin{array}{ccc}
\varphi_{z\overline{z}}^1 & \varphi_{z\overline{z}}^2 &
\varphi_{z\overline{z}}^3
\\
\varphi_{zz\overline{z}}^1 & \varphi_{zz\overline{z}}^2 &
\varphi_{zz\overline{z}}^3
\\
\varphi_{z\overline{z}\overline{z}}^1 & 
\varphi_{z\overline{z}\overline{z}}^2 &
\varphi_{z\overline{z}\overline{z}}^3
\end{array}
\right\vert
}
-
\frac{
\varphi_{2,zzz\overline{z}}\,
\left\vert\!
\begin{array}{cc}
\varphi_{zz\overline{z}}^1 & 
\varphi_{zz\overline{z}}^3\\
\varphi_{z\overline{z}\overline{z}}^1 &
\varphi_{z\overline{z}\overline{z}}^3
\end{array}
\!\right\vert
}{
\left\vert\!
\begin{array}{ccc}
\varphi_{z\overline{z}}^1 & \varphi_{z\overline{z}}^2 &
\varphi_{z\overline{z}}^3
\\
\varphi_{zz\overline{z}}^1 & \varphi_{zz\overline{z}}^2 &
\varphi_{zz\overline{z}}^3
\\
\varphi_{z\overline{z}\overline{z}}^1 & 
\varphi_{z\overline{z}\overline{z}}^2 &
\varphi_{z\overline{z}\overline{z}}^3
\end{array}
\right\vert
}
+
\frac{
\varphi_{3,zzz\overline{z}}\,
\left\vert\!
\begin{array}{cc}
\varphi_{zz\overline{z}}^1 & 
\varphi_{zz\overline{z}}^2\\
\varphi_{z\overline{z}\overline{z}}^1 & 
\varphi_{z\overline{z}\overline{z}}^2
\end{array}
\!\right\vert
}{
\left\vert\!
\begin{array}{ccc}
\varphi_{z\overline{z}}^1 & \varphi_{z\overline{z}}^2 &
\varphi_{z\overline{z}}^3
\\
\varphi_{zz\overline{z}}^1 & \varphi_{zz\overline{z}}^2 &
\varphi_{zz\overline{z}}^3
\\
\varphi_{z\overline{z}\overline{z}}^1 & 
\varphi_{z\overline{z}\overline{z}}^2 &
\varphi_{z\overline{z}\overline{z}}^3
\end{array}
\right\vert
}
\right)
\cdot
\mathcal{T}
+
\endaligned
\]
\[
\aligned
&
+
\left(
\frac{
-\,\varphi_{1,zzz\overline{z}}\,
\left\vert\!
\begin{array}{cc}
\varphi_{z\overline{z}}^2 & 
\varphi_{z\overline{z}}^3\\
\varphi_{z\overline{z}\overline{z}}^2 & 
\varphi_{z\overline{z}\overline{z}}^3
\end{array}
\!\right\vert
}{
\left\vert\!
\begin{array}{ccc}
\varphi_{z\overline{z}}^1 & \varphi_{z\overline{z}}^2 &
\varphi_{z\overline{z}}^3
\\
\varphi_{zz\overline{z}}^1 & \varphi_{zz\overline{z}}^2 &
\varphi_{zz\overline{z}}^3
\\
\varphi_{z\overline{z}\overline{z}}^1 & 
\varphi_{z\overline{z}\overline{z}}^2 &
\varphi_{z\overline{z}\overline{z}}^3
\end{array}
\right\vert
}
+
\frac{
\varphi_{2,zzz\overline{z}}\,
\left\vert\!
\begin{array}{cc}
\varphi_{z\overline{z}}^1 & 
\varphi_{z\overline{z}}^3\\
\varphi_{z\overline{z}\overline{z}}^1 &
\varphi_{z\overline{z}\overline{z}}^3
\end{array}
\!\right\vert
}{
\left\vert\!
\begin{array}{ccc}
\varphi_{z\overline{z}}^1 & \varphi_{z\overline{z}}^2 &
\varphi_{z\overline{z}}^3
\\
\varphi_{zz\overline{z}}^1 & \varphi_{zz\overline{z}}^2 &
\varphi_{zz\overline{z}}^3
\\
\varphi_{z\overline{z}\overline{z}}^1 & 
\varphi_{z\overline{z}\overline{z}}^2 &
\varphi_{z\overline{z}\overline{z}}^3
\end{array}
\right\vert
}
-
\frac{
\varphi_{3,zzz\overline{z}}\,
\left\vert\!
\begin{array}{cc}
\varphi_{z\overline{z}}^1 & 
\varphi_{z\overline{z}}^2\\
\varphi_{z\overline{z}\overline{z}}^1 & 
\varphi_{z\overline{z}\overline{z}}^2
\end{array}
\!\right\vert
}{
\left\vert\!
\begin{array}{ccc}
\varphi_{z\overline{z}}^1 & \varphi_{z\overline{z}}^2 &
\varphi_{z\overline{z}}^3
\\
\varphi_{zz\overline{z}}^1 & \varphi_{zz\overline{z}}^2 &
\varphi_{zz\overline{z}}^3
\\
\varphi_{z\overline{z}\overline{z}}^1 & 
\varphi_{z\overline{z}\overline{z}}^2 &
\varphi_{z\overline{z}\overline{z}}^3
\end{array}
\right\vert
}
\right)
\cdot
\mathcal{S}
+
\endaligned
\]
\[
\aligned
&
+
\left(
\frac{
\varphi_{1,zzz\overline{z}}\,
\left\vert\!
\begin{array}{cc}
\varphi_{z\overline{z}}^2 & 
\varphi_{z\overline{z}}^3\\
\varphi_{zz\overline{z}}^2 & 
\varphi_{zz\overline{z}}^3
\end{array}
\!\right\vert
}{
\left\vert\!
\begin{array}{ccc}
\varphi_{z\overline{z}}^1 & \varphi_{z\overline{z}}^2 &
\varphi_{z\overline{z}}^3
\\
\varphi_{zz\overline{z}}^1 & \varphi_{zz\overline{z}}^2 &
\varphi_{zz\overline{z}}^3
\\
\varphi_{z\overline{z}\overline{z}}^1 & 
\varphi_{z\overline{z}\overline{z}}^2 &
\varphi_{z\overline{z}\overline{z}}^3
\end{array}
\right\vert
}
-
\frac{
\varphi_{2,zzz\overline{z}}\,
\left\vert\!
\begin{array}{cc}
\varphi_{z\overline{z}}^1 & 
\varphi_{z\overline{z}}^3\\
\varphi_{zz\overline{z}}^1 &
\varphi_{zz\overline{z}}^3
\end{array}
\!\right\vert
}{
\left\vert\!
\begin{array}{ccc}
\varphi_{z\overline{z}}^1 & \varphi_{z\overline{z}}^2 &
\varphi_{z\overline{z}}^3
\\
\varphi_{zz\overline{z}}^1 & \varphi_{zz\overline{z}}^2 &
\varphi_{zz\overline{z}}^3
\\
\varphi_{z\overline{z}\overline{z}}^1 & 
\varphi_{z\overline{z}\overline{z}}^2 &
\varphi_{z\overline{z}\overline{z}}^3
\end{array}
\right\vert
}
+
\frac{
\varphi_{3,zzz\overline{z}}\,
\left\vert\!
\begin{array}{cc}
\varphi_{z\overline{z}}^1 & 
\varphi_{z\overline{z}}^2\\
\varphi_{zz\overline{z}}^1 & 
\varphi_{zz\overline{z}}^2
\end{array}
\!\right\vert
}{
\left\vert\!
\begin{array}{ccc}
\varphi_{z\overline{z}}^1 & \varphi_{z\overline{z}}^2 &
\varphi_{z\overline{z}}^3
\\
\varphi_{zz\overline{z}}^1 & \varphi_{zz\overline{z}}^2 &
\varphi_{zz\overline{z}}^3
\\
\varphi_{z\overline{z}\overline{z}}^1 & 
\varphi_{z\overline{z}\overline{z}}^2 &
\varphi_{z\overline{z}\overline{z}}^3
\end{array}
\right\vert
}
\right)
\cdot
\overline{\mathcal{S}},
\endaligned
\]
and:
\[
\!\!\!\!\!\!\!\!\!\!\!\!\!\!\!\!\!\!\!\!
\aligned
\big[\overline{\mathcal{L}},\mathcal{S}\big]
&
=
2\,\varphi_{1,zz\overline{z}\overline{z}}\,
\frac{\partial}{\partial u_1}
+
2\,\varphi_{2,zz\overline{z}\overline{z}}\,
\frac{\partial}{\partial u_2}
+
2\,\varphi_{3,zz\overline{z}\overline{z}}\,
\frac{\partial}{\partial u_3}
\\
&
=
\left(
\frac{
\varphi_{1,zz\overline{z}\overline{z}}\,
\left\vert\!
\begin{array}{cc}
\varphi_{zz\overline{z}}^2 & 
\varphi_{zz\overline{z}}^3\\
\varphi_{z\overline{z}\overline{z}}^2 & 
\varphi_{z\overline{z}\overline{z}}^3
\end{array}
\!\right\vert
}{
\left\vert\!
\begin{array}{ccc}
\varphi_{z\overline{z}}^1 & \varphi_{z\overline{z}}^2 &
\varphi_{z\overline{z}}^3
\\
\varphi_{zz\overline{z}}^1 & \varphi_{zz\overline{z}}^2 &
\varphi_{zz\overline{z}}^3
\\
\varphi_{z\overline{z}\overline{z}}^1 & 
\varphi_{z\overline{z}\overline{z}}^2 &
\varphi_{z\overline{z}\overline{z}}^3
\end{array}
\right\vert
}
-
\frac{
\varphi_{2,zz\overline{z}\overline{z}}\,
\left\vert\!
\begin{array}{cc}
\varphi_{zz\overline{z}}^1 & 
\varphi_{zz\overline{z}}^3\\
\varphi_{z\overline{z}\overline{z}}^1 &
\varphi_{z\overline{z}\overline{z}}^3
\end{array}
\!\right\vert
}{
\left\vert\!
\begin{array}{ccc}
\varphi_{z\overline{z}}^1 & \varphi_{z\overline{z}}^2 &
\varphi_{z\overline{z}}^3
\\
\varphi_{zz\overline{z}}^1 & \varphi_{zz\overline{z}}^2 &
\varphi_{zz\overline{z}}^3
\\
\varphi_{z\overline{z}\overline{z}}^1 & 
\varphi_{z\overline{z}\overline{z}}^2 &
\varphi_{z\overline{z}\overline{z}}^3
\end{array}
\right\vert
}
+
\frac{
\varphi_{3,zz\overline{z}\overline{z}}\,
\left\vert\!
\begin{array}{cc}
\varphi_{zz\overline{z}}^1 & 
\varphi_{zz\overline{z}}^2\\
\varphi_{z\overline{z}\overline{z}}^1 & 
\varphi_{z\overline{z}\overline{z}}^2
\end{array}
\!\right\vert
}{
\left\vert\!
\begin{array}{ccc}
\varphi_{z\overline{z}}^1 & \varphi_{z\overline{z}}^2 &
\varphi_{z\overline{z}}^3
\\
\varphi_{zz\overline{z}}^1 & \varphi_{zz\overline{z}}^2 &
\varphi_{zz\overline{z}}^3
\\
\varphi_{z\overline{z}\overline{z}}^1 & 
\varphi_{z\overline{z}\overline{z}}^2 &
\varphi_{z\overline{z}\overline{z}}^3
\end{array}
\right\vert
}
\right)
\cdot
\mathcal{T}
+
\endaligned
\]
\[
\aligned
&
+
\left(
\frac{
-\,\varphi_{1,zz\overline{z}\overline{z}}\,
\left\vert\!
\begin{array}{cc}
\varphi_{z\overline{z}}^2 & 
\varphi_{z\overline{z}}^3\\
\varphi_{z\overline{z}\overline{z}}^2 & 
\varphi_{z\overline{z}\overline{z}}^3
\end{array}
\!\right\vert
}{
\left\vert\!
\begin{array}{ccc}
\varphi_{z\overline{z}}^1 & \varphi_{z\overline{z}}^2 &
\varphi_{z\overline{z}}^3
\\
\varphi_{zz\overline{z}}^1 & \varphi_{zz\overline{z}}^2 &
\varphi_{zz\overline{z}}^3
\\
\varphi_{z\overline{z}\overline{z}}^1 & 
\varphi_{z\overline{z}\overline{z}}^2 &
\varphi_{z\overline{z}\overline{z}}^3
\end{array}
\right\vert
}
+
\frac{
\varphi_{2,zz\overline{z}\overline{z}}\,
\left\vert\!
\begin{array}{cc}
\varphi_{z\overline{z}}^1 & 
\varphi_{z\overline{z}}^3\\
\varphi_{z\overline{z}\overline{z}}^1 &
\varphi_{z\overline{z}\overline{z}}^3
\end{array}
\!\right\vert
}{
\left\vert\!
\begin{array}{ccc}
\varphi_{z\overline{z}}^1 & \varphi_{z\overline{z}}^2 &
\varphi_{z\overline{z}}^3
\\
\varphi_{zz\overline{z}}^1 & \varphi_{zz\overline{z}}^2 &
\varphi_{zz\overline{z}}^3
\\
\varphi_{z\overline{z}\overline{z}}^1 & 
\varphi_{z\overline{z}\overline{z}}^2 &
\varphi_{z\overline{z}\overline{z}}^3
\end{array}
\right\vert
}
-
\frac{
\varphi_{3,zz\overline{z}\overline{z}}\,
\left\vert\!
\begin{array}{cc}
\varphi_{z\overline{z}}^1 & 
\varphi_{z\overline{z}}^2\\
\varphi_{z\overline{z}\overline{z}}^1 & 
\varphi_{z\overline{z}\overline{z}}^2
\end{array}
\!\right\vert
}{
\left\vert\!
\begin{array}{ccc}
\varphi_{z\overline{z}}^1 & \varphi_{z\overline{z}}^2 &
\varphi_{z\overline{z}}^3
\\
\varphi_{zz\overline{z}}^1 & \varphi_{zz\overline{z}}^2 &
\varphi_{zz\overline{z}}^3
\\
\varphi_{z\overline{z}\overline{z}}^1 & 
\varphi_{z\overline{z}\overline{z}}^2 &
\varphi_{z\overline{z}\overline{z}}^3
\end{array}
\right\vert
}
\right)
\cdot
\mathcal{S}
+
\endaligned
\]
\[
\aligned
&
+
\left(
\frac{
\varphi_{1,zz\overline{z}\overline{z}}\,
\left\vert\!
\begin{array}{cc}
\varphi_{z\overline{z}}^2 & 
\varphi_{z\overline{z}}^3\\
\varphi_{zz\overline{z}}^2 & 
\varphi_{zz\overline{z}}^3
\end{array}
\!\right\vert
}{
\left\vert\!
\begin{array}{ccc}
\varphi_{z\overline{z}}^1 & \varphi_{z\overline{z}}^2 &
\varphi_{z\overline{z}}^3
\\
\varphi_{zz\overline{z}}^1 & \varphi_{zz\overline{z}}^2 &
\varphi_{zz\overline{z}}^3
\\
\varphi_{z\overline{z}\overline{z}}^1 & 
\varphi_{z\overline{z}\overline{z}}^2 &
\varphi_{z\overline{z}\overline{z}}^3
\end{array}
\right\vert
}
-
\frac{
\varphi_{2,zz\overline{z}\overline{z}}\,
\left\vert\!
\begin{array}{cc}
\varphi_{z\overline{z}}^1 & 
\varphi_{z\overline{z}}^3\\
\varphi_{zz\overline{z}}^1 &
\varphi_{zz\overline{z}}^3
\end{array}
\!\right\vert
}{
\left\vert\!
\begin{array}{ccc}
\varphi_{z\overline{z}}^1 & \varphi_{z\overline{z}}^2 &
\varphi_{z\overline{z}}^3
\\
\varphi_{zz\overline{z}}^1 & \varphi_{zz\overline{z}}^2 &
\varphi_{zz\overline{z}}^3
\\
\varphi_{z\overline{z}\overline{z}}^1 & 
\varphi_{z\overline{z}\overline{z}}^2 &
\varphi_{z\overline{z}\overline{z}}^3
\end{array}
\right\vert
}
+
\frac{
\varphi_{3,zz\overline{z}\overline{z}}\,
\left\vert\!
\begin{array}{cc}
\varphi_{z\overline{z}}^1 & 
\varphi_{z\overline{z}}^2\\
\varphi_{zz\overline{z}}^1 & 
\varphi_{zz\overline{z}}^2
\end{array}
\!\right\vert
}{
\left\vert\!
\begin{array}{ccc}
\varphi_{z\overline{z}}^1 & \varphi_{z\overline{z}}^2 &
\varphi_{z\overline{z}}^3
\\
\varphi_{zz\overline{z}}^1 & \varphi_{zz\overline{z}}^2 &
\varphi_{zz\overline{z}}^3
\\
\varphi_{z\overline{z}\overline{z}}^1 & 
\varphi_{z\overline{z}\overline{z}}^2 &
\varphi_{z\overline{z}\overline{z}}^3
\end{array}
\right\vert
}
\right)
\cdot
\overline{\mathcal{S}},
\endaligned
\]

while trivially (only in the rigid case!):
\[
\aligned
\big[\mathcal{T},\mathcal{S}\big]
&
=
0,
\\
\big[\mathcal{T},\overline{\mathcal{S}}\big]
&
=
0,
\\
\big[\mathcal{S},\overline{\mathcal{S}}\big]
&
=
0.
\endaligned
\]

\medskip\noindent{\bf Symbolic treatment of the general case.}
Now, come back to the general case where $\varphi_1$, $\varphi_2$,
$\varphi_3$ do depend on all variables $(x, y, u_1, u_2, u_3)$, so
that:
\[
\aligned
\mathcal{L}
&
=
\frac{\partial}{\partial z}
+
\frac{\Lambda_1}{\Delta}\,\frac{\partial}{\partial u_1}
+
\frac{\Lambda_2}{\Delta}\,\frac{\partial}{\partial u_2}
+
\frac{\Lambda_3}{\Delta}\,\frac{\partial}{\partial u_3},
\\
\overline{\mathcal{L}}
&
=
\frac{\partial}{\partial\overline{z}}
+
\frac{\overline{\Lambda}_1}{\overline{\Delta}}\,
\frac{\partial}{\partial u_1}
+
\frac{\overline{\Lambda}_2}{\overline{\Delta}}\,
\frac{\partial}{\partial u_2}
+
\frac{\overline{\Lambda}_3}{\overline{\Delta}}\,
\frac{\partial}{\partial u_3},
\\
\mathcal{T}
&
=
\ \ \ \ \ \ \ \ \ \ 
\frac{\Upsilon_1}{\Delta^2\,\overline{\Delta}^2}\,
\frac{\partial}{\partial u_1}
+
\frac{\Upsilon_2}{\Delta^2\,\overline{\Delta}^2}\,
\frac{\partial}{\partial u_2}
+
\frac{\Upsilon_3}{\Delta^2\,\overline{\Delta}^2}\,
\frac{\partial}{\partial u_3},
\\
\mathcal{S}
&
=
\ \ \ \ \ \ \ \ \ \ 
\frac{\Pi_1}{\Delta^4\,\overline{\Delta}^3}\,
\frac{\partial}{\partial u_1}
+
\frac{\Pi_2}{\Delta^4\,\overline{\Delta}^3}\,
\frac{\partial}{\partial u_2}
+
\frac{\Pi_3}{\Delta^4\,\overline{\Delta}^3}\,
\frac{\partial}{\partial u_3},
\\
\overline{\mathcal{S}}
&
=
\ \ \ \ \ \ \ \ \ \ 
\frac{\overline{\Pi}_1}{\Delta^3\,\overline{\Delta}^4}\,
\frac{\partial}{\partial u_1}
+
\frac{\overline{\Pi}_2}{\Delta^3\,\overline{\Delta}^4}\,
\frac{\partial}{\partial u_2}
+
\frac{\overline{\Pi}_3}{\Delta^3\,\overline{\Delta}^4}\,
\frac{\partial}{\partial u_3}.
\endaligned
\]

Since, as always, $\mathcal{ T}$ is real, one sees:
\[
\overline{\mathcal{S}}
=
\overline{
\big[\mathcal{L},\mathcal{T}\big]}
=
\big[\overline{\mathcal{L}},\mathcal{T}\big].
\]

Now, introduce named functions so that:
\[
\aligned
\big[\mathcal{L},\mathcal{S}\big]
&
=
P\cdot\mathcal{T}
+
Q\cdot\mathcal{S}
+
R\cdot\overline{\mathcal{S}},
\\
\big[\overline{\mathcal{S}},\mathcal{S}\big]
&
=
A\cdot\mathcal{T}
+
B\cdot\mathcal{S}
+
C\cdot\overline{\mathcal{S}}.
\endaligned
\]

\medskip\noindent{\bf Lemma.}
{\em One has the reality condition:}
\[
\big[\overline{\mathcal{L}},\mathcal{S}\big]
=
\big[\mathcal{L},\overline{\mathcal{S}}\big].
\]

\proof
Indeed, a consequence of the Jacobi identity shows that for
any two words $h_1$ and $h_2$ in a free Lie algebra:
\[
\big[h_2,\,\big[h_1,\,[h_1,h_2]\big]\big]
=
\big[h_1,\,\big[h_2,\,[h_1,h_2]\big]\big].
\]
Apply:
\[
\aligned
\big[\overline{\mathcal{L}},\mathcal{S}\big]
&
=
\big[\overline{\mathcal{L}},\,
\big[\mathcal{L},\,
\isqrt\,\big[\mathcal{L},\overline{\mathcal{L}}\big]
\big]\big]
\\
&
=
\isqrt\,
\big[\overline{\mathcal{L}},\,
\big[\mathcal{L},\,
\big[\mathcal{L},\overline{\mathcal{L}}\big]
\big]\big]
\\
&
=
\isqrt\,
\big[\mathcal{L},\,
\big[\overline{\mathcal{L}},\,
\big[\mathcal{L},\overline{\mathcal{L}}\big]
\big]\big]
\\
&
=
\big[\mathcal{L},\,
\big[\overline{\mathcal{L}},\,
\isqrt\,\big[\mathcal{L},\overline{\mathcal{L}}\big]
\big]\big]
\\
&
=
\big[\mathcal{L},\overline{\mathcal{S}}\big]
\endaligned
\]
which is so.
\endproof

Hence necessarily:
\[
\aligned
A
&
=
\overline{A},
\\
C
&
=
\overline{B},
\endaligned
\]
and gathering:
\[
\boxed{
\aligned
\big[\mathcal{L},\mathcal{S}\big]
&
=
P\cdot\mathcal{T}
+
Q\cdot\mathcal{S}
+
R\cdot\overline{\mathcal{S}},
\\
\big[\overline{\mathcal{L}},\mathcal{S}\big]
&
=
A\cdot\mathcal{T}
+
B\cdot\mathcal{S}
+
\overline{B}\cdot\overline{\mathcal{S}},
\\
\big[\mathcal{L},\overline{\mathcal{S}}\big]
&
=
A\cdot\mathcal{T}
+
B\cdot\mathcal{S}
+
\overline{B}\cdot\overline{\mathcal{S}},
\\
\big[\overline{\mathcal{L}},\overline{\mathcal{S}}\big]
&
=
\overline{P}\cdot\mathcal{T}
+
\overline{R}\cdot\mathcal{S}
+
\overline{Q}\cdot\overline{\mathcal{S}}.\,\,
\endaligned}
\]
These five functions:
\[
A,\ \ \ \ \ \ \ \ \
B,\ \ \ \ \ \ \ \ \
P,\ \ \ \ \ \ \ \ \
Q,\ \ \ \ \ \ \ \ \
R,
\]
will appear to be the only fundamental ones.

\medskip

Indeed, it yet remains to compute the three brackets:
\[
\aligned
&
\big[\mathcal{T},\mathcal{S}\big],
\\
&
\big[\mathcal{T},\overline{\mathcal{S}}\big],
\\
&
\big[\mathcal{S},\overline{\mathcal{S}}\big],
\endaligned
\]
using the Jacobi identity.

\medskip\noindent{\bf Lemma.}
{\em The coefficients of the two Lie brackets:}
\[
\aligned
\big[\mathcal{T},\mathcal{S}\big]
&
=
E_{\sf rpl}\cdot\mathcal{T}
+
F_{\sf rpl}\cdot\mathcal{S}
+
G_{\sf rpl}\cdot\overline{\mathcal{S}},
\\
\big[\mathcal{T},\overline{\mathcal{S}}\big]
&
=
\overline{E}_{\sf rpl}\cdot\mathcal{T}
+
\overline{G}_{\sf rpl}\cdot\mathcal{S}
+
\overline{F}_{\sf rpl}\cdot\overline{\mathcal{S}}
\endaligned
\]
{\em express in terms of $A$, $B$, $P$, $Q$, $R$ as:}
\[
\aligned
E_{\sf rpl}
&
=
\isqrt\,
\Big(
\mathcal{L}(A)
-
\overline{\mathcal{L}}(P)
+
A\overline{B}
+
BP
-
AQ
-
\overline{P}R
\Big),
\\
F_{\sf rpl}
&
=
\isqrt\,
\Big(
\mathcal{L}(B)
-
\overline{\mathcal{L}}(Q)
+
A
+
B\overline{B}
-
R\overline{R}
\Big),
\\
G_{\sf rpl}
&
=
\isqrt\,
\Big(
\mathcal{L}(\overline{B})
-
\overline{\mathcal{L}}(R)
+
\overline{B}\overline{B}
+
BR
-
P
-
\overline{B}Q
-
R\overline{Q}
\Big).
\endaligned
\]

\proof
Indeed, for $h_1$ and $h_2$ any two words in a free Lie algebra,
the Jacobi identity gives:
\[
\aligned
\Big[
[h_1,h_2],\,\,
\big[h_1,\,[h_1,h_2]\big]
\Big]
&
=
-\,
\Big[
\Big[
\big[h_1,\,[h_1,h_2]\big],\,h_1
\Big],\,\,
h_2
\Big]
-
\Big[
\Big[
h_2,\,\big[h_1,\,[h_1,h_2]\big]
\Big],\,\,
h_1
\Big]
\\
&
=
-\,
\Big[
h_2,\,\,
\Big[
h_1,\,
\big[h_1,\,[h_1,h_2]\big]
\Big]
\Big]
+
\Big[
h_1,\,\,
\Big[
h_2,\,
\big[h_1,\,[h_1,h_2]\big]
\Big]
\Big].
\endaligned
\]
Apply this to $h_1 := \mathcal{L}$ and $h_2 := \overline{\mathcal{L}}$
and get:
\[
\aligned
-\,\isqrt\,\big[\mathcal{T},\mathcal{S}\big]
&
=
\Big[
\big[\mathcal{L},\overline{\mathcal{L}}\big],\,\,
\big[\mathcal{L},\,\big[\mathcal{L},\overline{\mathcal{L}}\big]\big]
\Big]
\\
&
=
-\,
\Big[
\overline{\mathcal{L}},\,
\underbrace{\Big[\mathcal{L},\,
\big[\mathcal{L},\big[\mathcal{L},\overline{\mathcal{L}}\big]\big]
\Big]}_{
[\mathcal{L},\mathcal{S}]}
\Big]
+
\Big[
\mathcal{L},\,
\underbrace{
\Big[\overline{\mathcal{L}},\,
\big[\mathcal{L},\big[\mathcal{L},\overline{\mathcal{L}}\big]\big]
\Big]}_{[\overline{\mathcal{L}},\mathcal{S}]}
\Big],
\endaligned
\]
that is to say after replacement:
\[
\aligned
-\,\isqrt\,\big[\mathcal{T},\mathcal{S}\big]
&
=
-\,
\Big[
\overline{\mathcal{L}},\,
P\cdot\mathcal{T}
+
Q\cdot\mathcal{S}
+
R\cdot\overline{\mathcal{S}}
\Big]
+
\Big[
\mathcal{L},\,
A\cdot\mathcal{T}
+
B\cdot\mathcal{S}
+
\overline{B}\cdot\overline{\mathcal{S}}
\Big]
\\
&
=
-\,\overline{\mathcal{L}}(P)
\cdot
\mathcal{T}
-
\overline{\mathcal{L}}(Q)
\cdot
\mathcal{S}
-
\overline{\mathcal{L}}(R)
\cdot
\overline{\mathcal{S}}
+
\\
&
\ \ \ \ \
+
\mathcal{L}(A)\cdot\mathcal{T}
+
\mathcal{L}(B)\cdot\mathcal{S}
+
\mathcal{L}(C)\cdot\overline{\mathcal{S}}
-
\\
&
\ \ \ \ \
-\,
P\,
\underbrace{
\big[\overline{\mathcal{L}},\mathcal{T}\big]}_{
\overline{\mathcal{S}}}
-
Q\,
\underbrace{\big[\overline{\mathcal{L}},\mathcal{S}\big]}_{
A\mathcal{T}+B\mathcal{S}+
\overline{B}\overline{\mathcal{S}}}
-
R\,
\underbrace{
\big[\overline{\mathcal{L}},\overline{\mathcal{S}}\big]}_{
\overline{P}\mathcal{T}+\overline{R}\mathcal{S}
+\overline{Q}\overline{\mathcal{S}}}
+
\\
&
\ \ \ \ \
+
A\,
\underbrace{\big[\mathcal{L},\mathcal{T}\big]}_{
\mathcal{S}}
+
B\,
\underbrace{\big[\mathcal{L},\mathcal{S}\big]}_{
P\mathcal{T}+Q\mathcal{S}+R\overline{\mathcal{S}}}
+
\overline{B}\,
\underbrace{
\big[\mathcal{L},\overline{\mathcal{S}}\big]}_{
\overline{A}\mathcal{T}+B\mathcal{S}+\overline{B}\overline{\mathcal{S}}},
\endaligned
\]
which collects as:
\[
\aligned
-\,\isqrt\,\big[\mathcal{T},\mathcal{S}\big]
&
=
\Big(
\mathcal{L}(A)
-
\overline{\mathcal{L}}(P)
-
AQ
-
\overline{P}R
+
BP
+
A\overline{B}
\Big)
\cdot
\mathcal{T}
+
\\
&
+
\Big(
\mathcal{L}(B)
-
\overline{\mathcal{L}}(Q)
+
A
+
B\overline{B}
-
R\overline{R}
\Big)
\cdot
\mathcal{S}
+
\\
&
+
\Big(
\mathcal{L}(C)
-
\overline{\mathcal{L}}(R)
-
P
+
BR
-
\overline{B}Q
-
R\overline{Q}
+
\overline{B}\overline{B}
\Big)
\cdot
\overline{\mathcal{S}},
\endaligned
\]
as stated.
\endproof

Notice that the last remaining Lie bracket is
purely imaginary:
\[
\overline{\big[\mathcal{S},\overline{\mathcal{S}}\big]}
=
\big[\overline{\mathcal{S}},\mathcal{S}\big]
=
-\,
\big[\mathcal{S},\overline{\mathcal{S}}\big].
\]

\medskip\noindent{\bf Lemma.}
{\em The coefficients of:}
\[
\big[\mathcal{S},\overline{\mathcal{S}}\big]
=
\isqrt\,J_{\sf rpl}\cdot\mathcal{T}
+
K_{\sf rpl}\cdot\mathcal{S}
-
\overline{K}_{\sf rpl}\cdot\overline{\mathcal{S}}
\]
{\em express in terms of $A$, $B$, $P$, $Q$, $R$ as:}
\[
\!\!\!\!\!\!\!\!\!\!\!\!\!\!\!\!\!\!\!\!
\aligned
J_{\sf rpl}
=
{\textstyle{\frac{1}{2}}}\,
\Big(
&
\overline{\mathcal{L}}\big(\overline{\mathcal{L}}(P)\big)
+
\mathcal{L}\big(\mathcal{L}(\overline{P})\big)
-
\overline{\mathcal{L}}\big(\mathcal{L}(A)\big)
-
\mathcal{L}\big(\overline{\mathcal{L}}(A)\big)
+
\\
&
+
Q\,\overline{\mathcal{L}}(A)
+
\overline{Q}\,\mathcal{L}(A)
+
2\,A\,\overline{\mathcal{L}}(Q)
+
2\,A\,\mathcal{L}\big(\overline{Q}\big)
+
R\,\overline{\mathcal{L}}\big(\overline{P}\big)
+
\overline{R}\,\mathcal{L}(P)
+
\\
&
+
2\,\overline{P}\,\overline{\mathcal{L}}(R)
+
2\,P\,\mathcal{L}\big(\overline{R}\big)
-
2\,P\,\overline{\mathcal{L}}(B)
-
2\,\overline{P}\,\mathcal{L}\big(\overline{B}\big)
-
B\,\overline{\mathcal{L}}(P)
-
\overline{B}\,\mathcal{L}\big(\overline{P}\big)
\,-
\\
&
-\,
\overline{B}\,\overline{\mathcal{L}}(A)
-
B\,\mathcal{L}(A)
-
2\,A\,\overline{\mathcal{L}}\big(\overline{B}\big)
-
2\,A\,\mathcal{L}(B)
\,-
\\
&
-\,
2\,AA
-
2\,AB\overline{B}
-
BBP
-
\overline{B}\overline{B}\overline{P}
-
B\overline{P}R
-
\overline{B}P\overline{R}
+
PQ\overline{R}
+
\overline{P}\overline{Q}R
+
\\
&
+
2\,AR\overline{R}
+
2\,P\overline{P}
+
\overline{B}\overline{P}Q
+
BP\overline{Q}
\Big),
\endaligned
\]
which is real, and:
\[
\!\!\!\!\!\!\!\!\!\!\!\!\!\!\!\!\!\!\!\!
\!\!\!\!\!\!\!\!\!\!\!\!\!\!\!\!\!\!\!\!
\aligned
K_{\sf rpl}
=
{\textstyle{\frac{1}{2}}}\,
\Big(
&
\overline{\mathcal{L}}\big(\overline{\mathcal{L}}(Q)\big)
-
\overline{\mathcal{L}}\big(\mathcal{L}(B)\big)
-
\mathcal{L}\big(\overline{\mathcal{L}}(B)\big)
+
\mathcal{L}\big(\mathcal{L}(\overline{R})\big)
+
\\
&
+
2\,\overline{R}\,\overline{\mathcal{L}}(R)
+
R\,\overline{\mathcal{L}}(\overline{R})
+
B\,\overline{\mathcal{L}}(Q)
+
2\,\mathcal{L}(\overline{P})
+
\overline{R}\,\mathcal{L}(Q)
+
2\,Q\,\mathcal{L}(\overline{R})
+
\overline{Q}\,\mathcal{L}(B)
-
Q\,\overline{\mathcal{L}}(B)
+
\\
&
+
2\,B\,\mathcal{L}(\overline{Q})
-
2\,\overline{\mathcal{L}}(A)
-
\overline{B}\,\overline{\mathcal{L}}(B)
-
2\,B\,\overline{\mathcal{L}}(\overline{B})
-
3\,B\,\mathcal{L}(B)
-
2\,\overline{R}\,\mathcal{L}(\overline{B})
-
\overline{B}\,\mathcal{L}(\overline{R})
-
\\
&
-
3\,AB
-
BBQ
-
\overline{B}\overline{B}\overline{R}
-
2\,BB\overline{B}
-
\overline{B}\overline{P}
+
A\overline{Q}
+
\overline{P}Q
+
\\
&
+
QQ\overline{R}
+
BQ\overline{Q}
+
BR\overline{R}
+
2\,P\overline{R}
+
\overline{Q}R\overline{R}
\Big).
\endaligned
\]

\proof
In subsection~7.3 of~\cite{ Aghasi-Merker-Sabzevari-5-2011}), 
in a free Lie containing $h_1$ and $h_2$, one has:
\[
\aligned
0
&
=
\big[\big[h_1,[h_1,h_2]\big],
\big[h_2,[h_1,h_2]\big]\big]
-
\big[[h_1,h_2],
\big[h_1,\big[h_2,[h_1,h_2]\big]\big]\big]
-
\\
&
\ \ \ \ \
-\,
\big[h_2,\big[h_2,\big[h_1,\big[h_1,
[h_1,h_2]\big]\big]\big]\big]
+
\big[h_2,\big[h_1,\big[h_1,\big[h_2,
[h_1,h_2]\big]\big]\big]\big],
\\
0
&
=
\big[\big[h_1,[h_1,h_2]\big],
\big[h_2,[h_1,h_2]\big]\big]
+
\big[[h_1,h_2],
\big[h_1,\big[h_2,[h_1,h_2]\big]\big]\big]
+
\\
&
\ \ \ \ \
+
\big[h_1,\big[h_2,\big[h_1,\big[h_2,
[h_1,h_2]\big]\big]\big]\big]
-
\big[h_1,\big[h_1,\big[h_2,\big[h_2,
[h_1,h_2]\big]\big]\big]\big].
\endaligned
\]
Adding and replacing $h_1 := \mathcal{ L}$ and
$h_2 := \overline{ \mathcal{ L}}$:
\[
\aligned
-\,2\,
\Big[
\big[\mathcal{L},\,[\mathcal{L},\overline{\mathcal{L}}]\big],\,
\big[\overline{\mathcal{L}},\,[\mathcal{L},\overline{\mathcal{L}}]\big]
\Big]
&
=
-\,
\Big[\overline{\mathcal{L}},\,\Big[\overline{\mathcal{L}},\,
\big[\mathcal{L},\,\big[\mathcal{L},\,
[\mathcal{L},\overline{\mathcal{L}}]
\big]\big]\Big]\Big]
+
\\
&
\ \ \ \ \
+
\Big[\overline{\mathcal{L}},\,\Big[\mathcal{L},\,
\big[\mathcal{L},\,\big[\overline{\mathcal{L}},
[\mathcal{L},\overline{\mathcal{L}}]
\big]\big]\Big]\Big]
+
\\
&
\ \ \ \ \
+
\Big[\mathcal{L},\,\Big[\overline{\mathcal{L}},\,
\big[\mathcal{L},\,\big[\overline{\mathcal{L}},
[\mathcal{L},\overline{\mathcal{L}}]
\big]\big]\Big]\Big]
\,-
\\
&
\ \ \ \ \
-\,
\Big[\mathcal{L},\,\Big[\mathcal{L},\,
\big[\overline{\mathcal{L}},\,\big[\overline{\mathcal{L}},
[\mathcal{L},\overline{\mathcal{L}}]
\big]\big]\Big]\Big].
\endaligned
\]
Taking account of some $\isqrt$:
\[
\aligned
\big[\mathcal{S},\overline{\mathcal{S}}\big]
&
=
\Big[
\big[\mathcal{L},\,\isqrt\,
\big[\mathcal{L},\overline{\mathcal{L}}\big]\big],\,\,
\big[\overline{\mathcal{L}},\,\isqrt\,
\big[\mathcal{L},\overline{\mathcal{L}}\big]\big]
\Big]
\\
&
=
-\,
\Big[
\big[\mathcal{L},\,[\mathcal{L},\overline{\mathcal{L}}]\big],\,
\big[\overline{\mathcal{L}},\,[\mathcal{L},\overline{\mathcal{L}}]\big]
\Big]
\\
&
=
\frac{\isqrt}{2}\,
\bigg(
+
\Big[\overline{\mathcal{L}},\,\Big[\overline{\mathcal{L}},\,
\underbrace{
\big[\mathcal{L},\,\big[\mathcal{L},\,
\isqrt\,[\mathcal{L},\overline{\mathcal{L}}]
\big]\big]}_{
[\mathcal{L},\mathcal{S}]}
\Big]\Big]
-
\\
&
\ \ \ \ \ \ \ \ \ \ \ \ \ \ \ \  
-
\Big[\overline{\mathcal{L}},\,\Big[\mathcal{L},\,
\underbrace{
\big[\mathcal{L},\,\big[\overline{\mathcal{L}},
\isqrt\,[\mathcal{L},\overline{\mathcal{L}}]
\big]\big]}_{
[\mathcal{L},\overline{\mathcal{S}}]}
\Big]\Big]
-
\\
&
\ \ \ \ \ \ \ \ \ \ \ \ \ \ \ \  
-
\Big[\mathcal{L},\,\Big[\overline{\mathcal{L}},\,
\underbrace{
\big[\mathcal{L},\,\big[\overline{\mathcal{L}},
\isqrt\,[\mathcal{L},\overline{\mathcal{L}}]
\big]\big]}_{
[\mathcal{L},\overline{\mathcal{S}}]}
\Big]\Big]
+
\\
&
\ \ \ \ \ \ \ \ \ \ \ \ \ \ \ \ 
+
\Big[\mathcal{L},\,\Big[\mathcal{L},\,
\underbrace{
\big[\overline{\mathcal{L}},\,\big[\overline{\mathcal{L}},
\isqrt\,[\mathcal{L},\overline{\mathcal{L}}]
\big]\big]}_{
[\overline{\mathcal{L}},\overline{\mathcal{S}}]}
\Big]\Big].
\bigg),
\endaligned
\]
that is to say:
\[
\aligned
\big[\mathcal{S},\overline{\mathcal{S}}\big]
&
=
\isqrt\,J_{\sf rpl}\cdot\mathcal{T}
+
K_{\sf rpl}\cdot\mathcal{S}
-
\overline{K}_{\sf rpl}\cdot\overline{\mathcal{S}}
\\
&
=
\frac{\isqrt}{2}
\bigg(
\Big[\overline{\mathcal{L}},\,
\big[\overline{\mathcal{L}},\,
P\mathcal{T}+Q\mathcal{S}+R\overline{\mathcal{S}}
\big]\Big]
-
\Big[\overline{\mathcal{L}},\,
\big[\mathcal{L},\,
A\,\mathcal{T}+B\mathcal{S}+\overline{B}\overline{\mathcal{S}}
\big]\Big]
\,-
\\
&
\ \ \ \ \ \ \ \ \ \ 
-\,
\Big[\mathcal{L},\,
\big[\overline{\mathcal{L}},\,
A\,\mathcal{T}+B\mathcal{S}+\overline{B}\overline{\mathcal{S}}
\big]\Big]
+
\Big[\mathcal{L},\,
\big[\mathcal{L},\,
\overline{P}\,\mathcal{T}+\overline{R}\mathcal{S}
+\overline{Q}\overline{\mathcal{S}}
\big]\Big]
\bigg).
\endaligned
\]

On preliminarily computes:
\[
\aligned
\Big[\overline{\mathcal{L}},\,
P\mathcal{T}+Q\mathcal{S}+R\overline{\mathcal{S}}
\Big]
&
=
\overline{\mathcal{L}}(P)\,\mathcal{T}
+
\overline{\mathcal{L}}(Q)\,\mathcal{S}
+
\overline{\mathcal{L}}(R)\,\overline{\mathcal{S}}
+
\\
&
\ \ \ \ \ 
+
P\,\big[\overline{\mathcal{L}},\mathcal{T}\big]
+
Q\,\big[\overline{\mathcal{L}},\mathcal{S}\big]
+
R\,\big[\overline{\mathcal{L}},\overline{\mathcal{S}}\big]
\\
&
=
\big(
\overline{\mathcal{L}}(P)
+
AQ
+
\overline{P}R
\big)
\cdot
\mathcal{T}
+
\\
&
\ \ \ \ \
+
\big(
\overline{\mathcal{L}}(Q)
+
BQ
+
R\overline{R}
\big)
\cdot
\mathcal{S}
+
\\
&
\ \ \ \ \
+
\big(
\overline{\mathcal{L}}(R)
+
P
+
\overline{B}Q
+
\overline{Q}R
\big)
\cdot
\overline{\mathcal{S}}.
\endaligned
\]
Similarly:
\[
\aligned
\Big[
\mathcal{L},\,\,
A\,\mathcal{T}
+
B\,\mathcal{S}
+
\overline{B}\,\overline{\mathcal{S}}
\Big]
&
=
\big(
\mathcal{L}(A)+BP+A\overline{B}
\big)
\cdot
\mathcal{T}
+
\\
&
\ \ \ \ \
\big(
\mathcal{L}(B)
+
A
+
BQ
+
B\overline{B}
\big)
\cdot
\mathcal{S}
+
\\
&
\ \ \ \ \
\big(
\mathcal{L}(\overline{B})
+
BR
+
\overline{B}\overline{B}
\big)
\cdot
\overline{\mathcal{S}},
\endaligned
\]
with conjugate:
\[
\aligned
\Big[
\overline{\mathcal{L}},\,\,
A\,\mathcal{T}
+
B\,\mathcal{S}
+
\overline{B}\,\overline{\mathcal{S}}
\Big]
&
=
\big(
\overline{\mathcal{L}}(A)
+
\overline{B}\overline{P}
+
AB
\big)
\cdot
\mathcal{T}
+
\\
&
\ \ \ \ \
\big(
\overline{\mathcal{L}}(B)
+
\overline{B}\overline{R}
+
BB
\big)
\cdot
\mathcal{S}
+
\\
&
\ \ \ \ \
\big(
\overline{\mathcal{L}}(\overline{B})
+
B\overline{B}
+
A
+
\overline{B}\overline{Q}
\big)
\cdot
\overline{\mathcal{S}},
\endaligned
\]
and lastly:
\[
\aligned
\Big[
\mathcal{L},\,\,
\overline{P}\,\mathcal{T}
+
\overline{R}\,\mathcal{S}
+
\overline{Q}\,\overline{\mathcal{S}}
\Big]
&
=
\big(
\mathcal{L}(\overline{P})
+
P\overline{R}
+
A\overline{Q}
\big)
\cdot
\mathcal{T}
+
\\
&
\ \ \ \ \
\big(
\mathcal{L}(\overline{R})
+
\overline{P}
+
Q\overline{R}
+
B\overline{Q}
\big)
\cdot
\mathcal{S}
+
\\
&
\ \ \ \ \
\big(
\mathcal{L}(\overline{Q})
+
R\overline{R}
+
\overline{B}\overline{Q}
\big)
\cdot
\overline{\mathcal{S}},
\endaligned
\]

One can therefore computes:
\[
\!\!\!\!\!\!\!\!\!\!\!\!\!\!\!\!\!\!\!\!
\!\!\!\!\!\!\!\!\!\!\!\!\!\!\!\!\!\!\!\!
\aligned
\Big[
\overline{\mathcal{L}},\,\,
\Big[
\overline{\mathcal{L}},\,
\big[\mathcal{L},\mathcal{S}\big]
\Big]
\Big]
&
=
\Big[
\overline{\mathcal{L}},\,\,
\Big[
\overline{\mathcal{L}},\,\,
P\,\mathcal{T}
+
Q\,\mathcal{S}
+
R\,\overline{\mathcal{S}}
\Big]
\Big]
\\
&
=
\Big[
\overline{\mathcal{L}},\,\,
\big(
\overline{\mathcal{L}}(P)
+
AQ
+
\overline{P}R
\big)
\cdot\mathcal{T}
+
\\
&
\ \ \ \ \ \ \ \ \ \ \ \ \ 
+
\big(
\overline{\mathcal{L}}(Q)
+
BQ
+
R\overline{R}
\big)
\cdot
\mathcal{S}
+
\\
&
\ \ \ \ \ \ \ \ \ \ \ \ \ 
+
\big(
\overline{\mathcal{L}}(R)
+
P
+
\overline{B}Q
+
\overline{Q}R
\big)
\cdot
\overline{\mathcal{S}}
\Big]
\\
&
=
\Big(
\overline{\mathcal{L}}\big(\overline{\mathcal{L}}(P)\big)
+
Q\,\overline{\mathcal{L}}(A)
+
A\,\overline{\mathcal{L}}(Q)
+
R\,\overline{\mathcal{L}}\big(\overline{P}\big)
+
\overline{P}\,\overline{\mathcal{L}}(R)
\Big)
\cdot
\mathcal{T}
+
\\
&
+
\Big(
\overline{\mathcal{L}}\big(\overline{\mathcal{L}}(Q)\big)
+
Q\,\overline{\mathcal{L}}(B)
+
B\,\overline{\mathcal{L}}(Q)
+
\overline{R}\,\overline{\mathcal{L}}(R)
+
R\,\overline{\mathcal{L}}(\overline{R})
\Big)
\cdot
\mathcal{S}
+
\\
&
+
\Big(
\overline{\mathcal{L}}\big(\overline{\mathcal{L}}(R)\big)
+
\overline{L}(P)
+
Q\,\overline{\mathcal{L}}(B)
+
\overline{B}\,\overline{\mathcal{L}}(Q)
+
R\,\overline{\mathcal{L}}(\overline{Q})
+
\overline{Q}\,\overline{\mathcal{L}}(R)
\Big)
\cdot
\overline{\mathcal{S}}
+
\\
&
+
\big(\overline{\mathcal{L}}(P)+AQ+\overline{P}R\big)\,
\underbrace{
\big[\overline{\mathcal{L}},\mathcal{T}\big]}_{
\overline{\mathcal{S}}}
+
\\
&
+
\big(\overline{\mathcal{L}}(Q)+BQ+R\overline{R}\big)\,
\underbrace{
\big[\overline{\mathcal{L}},\mathcal{S}\big]}_{
A\mathcal{T}+B\mathcal{S}+\overline{B}\overline{\mathcal{S}}}
+
\\
&
+
\big(\overline{\mathcal{L}}(R)+P+\overline{B}Q+\overline{Q}R\big)\,
\underbrace{
\big[\overline{\mathcal{L}},\overline{\mathcal{S}}\big]}_{
\overline{P}\mathcal{T}+\overline{R}\mathcal{S}
+\overline{Q}\overline{\mathcal{S}}},
\endaligned
\]
that is to say after collecting:
\[
\!\!\!\!\!\!\!\!\!\!\!\!\!\!\!\!\!\!\!\!
\!\!\!\!\!\!\!\!\!\!\!\!\!\!\!\!\!\!\!\!
\aligned
\Big[
\overline{\mathcal{L}},\,\,
\Big[
\overline{\mathcal{L}},\,
\big[\mathcal{L},\mathcal{S}\big]
\Big]
\Big]
&
=
\mathcal{T}
\cdot
\left(
\aligned
\overline{\mathcal{L}}\big(\overline{\mathcal{L}}(P)\big)
&
+
Q\,\overline{\mathcal{L}}(A)
+
A\,\overline{\mathcal{L}}(Q)
+
R\,\overline{\mathcal{L}}(\overline{P})
+
\overline{P}\,\overline{\mathcal{L}}(R)
+
\\
&
+
A\,\overline{\mathcal{L}}(Q)
+
\overline{P}\,\overline{\mathcal{L}}(R)
+
\\
&
+
ABQ
+
AR\overline{R}
+
P\overline{P}
+
\overline{B}\overline{P}Q
+
\overline{P}\overline{Q}R
\endaligned
\right)
+
\\
&
+
\mathcal{S}
\cdot
\left(
\aligned
\overline{\mathcal{L}}\big(\overline{\mathcal{L}}(Q)\big)
&
+
Q\,\overline{\mathcal{L}}(B)
+
B\,\overline{\mathcal{L}}(Q)
+
\overline{R}\,\overline{\mathcal{L}}(R)
+
R\,\overline{\mathcal{L}}(\overline{R})
+
\\
&
+
B\,\overline{\mathcal{L}}(Q)
+
\overline{R}\,\overline{\mathcal{L}}(R)
+
\\
&
+
BBQ
+
BR\overline{R}
+
P\overline{R}
+
\overline{B}Q\overline{R}
+
\overline{Q}R\overline{R}
\endaligned
\right)
+
\\
&
+
\overline{\mathcal{S}}
\cdot
\left(
\aligned
\overline{\mathcal{L}}\big(\overline{\mathcal{L}}(R)\big)
&
+
\overline{\mathcal{L}}(P)
+
Q\,\overline{\mathcal{L}}(B)
+
\overline{B}\,\overline{\mathcal{L}}(Q)
+
R\,\overline{\mathcal{L}}(\overline{Q})
+
\overline{Q}\,\overline{\mathcal{L}}(R)
+
\\
&
+
\overline{\mathcal{L}}(P)
+
\overline{B}\,\overline{\mathcal{L}}(Q)
+
\overline{Q}\,\overline{\mathcal{L}}(R)
+
\\
&
+
AQ
+
\overline{P}R
+
B\overline{B}Q
+
\overline{B}R\overline{R}
+
P\overline{Q}
+
\overline{B}Q\overline{Q}
+
\overline{Q}\overline{Q}R
\endaligned
\right),
\endaligned
\]

Quite similar computations the remaining three provide
(mind minus signs):
\[
\!\!\!\!\!\!\!\!\!\!
\!\!\!\!\!\!\!\!\!\!\!\!\!\!\!\!\!\!\!\!
\!\!\!\!\!\!\!\!\!\!\!\!\!\!\!\!\!\!\!\!
\aligned
-\,
\Big[
\overline{\mathcal{L}},\,\,
\Big[
\mathcal{L},\,
\big[\mathcal{L},\overline{\mathcal{S}}\big]
\Big]
\Big]
&
=
\mathcal{T}
\cdot
\left(
\aligned
-\,\overline{\mathcal{L}}\big(\mathcal{L}(A)\big)
&
-
P\,\overline{\mathcal{L}}(B)
-
B\,\overline{\mathcal{L}}(P)
-
\overline{B}\,\overline{\mathcal{L}}(A)
-
A\,\overline{\mathcal{L}}\big(\overline{B}\big)
\,-
\\
&
-\,
A\,\mathcal{L}(B)
-
\overline{P}\,\mathcal{L}\big(\overline{B}\big)
-
\\
&
-\,
AA
-
ABQ
-
AB\overline{B}
-
B\overline{P}R
-
\overline{B}\overline{B}\overline{P}
\endaligned
\right)
+
\\
&
+
\mathcal{S}
\cdot
\left(
\aligned
-\,\overline{\mathcal{L}}\big(\mathcal{L}(B)\big)
&
-
\overline{\mathcal{L}}(A)
-
Q\,\overline{\mathcal{L}}(B)
-
B\,\overline{\mathcal{L}}(Q)
-
\overline{B}\,\overline{\mathcal{L}}(B)
-
B\,\overline{\mathcal{L}}\big(\overline{B}\big)
\,-
\\
&
-\,
B\,\mathcal{L}(B)
-
\overline{R}\,\mathcal{L}\big(\overline{B}\big)
\,-
\\
&
-\,
AB
-
BBQ
-
BB\overline{B}
-
BR\overline{R}
-
\overline{B}\overline{B}\overline{R}
\endaligned
\right)
+
\\
&
+
\overline{\mathcal{S}}
\cdot
\left(
\aligned
-\,\overline{\mathcal{L}}\big(\mathcal{L}(\overline{B})\big)
&
-
R\,\overline{\mathcal{L}}(B)
-
B\,\overline{\mathcal{L}}(R)
-
2\,\overline{B}\,\overline{\mathcal{L}}(\overline{B})
\,-
\\
&
-\,
\mathcal{L}(A)
-
\overline{B}\,\mathcal{L}(B)
-
\overline{Q}\,\mathcal{L}\big(\overline{B}\big)
\,-
\\
&
-\,
BP
-
2\,A\overline{B}
-
B\overline{B}Q
-
B\overline{B}\overline{B}
-
B\overline{Q}R
-
\overline{B}\overline{B}\overline{Q}
\endaligned
\right)
\endaligned
\]
\[
\!\!\!\!\!\!\!\!\!\!
\!\!\!\!\!\!\!\!\!\!\!\!\!\!\!\!\!\!\!\!
\!\!\!\!\!\!\!\!\!\!\!\!\!\!\!\!\!\!\!\!
\aligned
-\,
\Big[
\mathcal{L},\,\,
\Big[
\overline{\mathcal{L}},\,
\big[\mathcal{L},\overline{\mathcal{S}}\big]
\Big]
\Big]
&
=
\mathcal{T}
\cdot
\left(
\aligned
-\,\mathcal{L}\big(\overline{\mathcal{L}}(A)\big)
&
-
\overline{P}\,\mathcal{L}\big(\overline{B}\big)
-
\overline{B}\,\mathcal{L}\big(\overline{P}\big)
-
B\,\mathcal{L}(A)
-
A\,\mathcal{L}(B)
\,-
\\
&
-\,
A\,\overline{\mathcal{L}}\big(\overline{B}\big)
-
P\,\overline{\mathcal{L}}(B)
-
\\
&
-\,
AA
-
A\overline{B}\overline{Q}
-
AB\overline{B}
-
\overline{B}P\overline{R}
-
BBP
\endaligned
\right)
+
\\
&
+
\mathcal{S}
\cdot
\left(
\aligned
-\,
\mathcal{L}\big(\overline{\mathcal{L}}(B)\big)
&
-
\overline{R}\,\mathcal{L}\big(\overline{B}\big)
-
\overline{B}\,\mathcal{L}\big(\overline{R}\big)
-
2\,B\,\mathcal{L}(B)
-
\\
&
-\,
\overline{\mathcal{L}}(A)
-
B\,\overline{\mathcal{L}}\big(\overline{B}\big)
-
Q\,\overline{\mathcal{L}}(B)
\,-
\\
&
-\,
\overline{B}\overline{P}
-
2\,A\,B
-
B\overline{B}\overline{Q}
-
BB\overline{B}
-
\overline{B}Q\overline{R}
-
BBQ
\endaligned
\right)
+
\\
&
+
\overline{\mathcal{S}}
\cdot
\left(
\aligned
-\,\mathcal{L}\big(\overline{\mathcal{L}}\big(
\overline{B}\big)\big)
&
-
\mathcal{L}(A)
-
\overline{Q}\,\mathcal{L}\big(\overline{B}\big)
-
\overline{B}\,\mathcal{L}\big(\overline{Q}\big)
-
B\,\mathcal{L}\big(\overline{B}\big)
-
\overline{B}\,\mathcal{L}(B)
\,-
\\
&
-\,
\overline{B}\,\overline{\mathcal{L}}\big(\overline{B}\big)
-
R\,\overline{\mathcal{L}}(B)
\,-
\\
&
-\,
A\overline{B}
-
\overline{B}\overline{B}\overline{Q}
-
B\overline{B}\overline{B}
-
\overline{B}R\overline{R}
-
BBR
\endaligned
\right),
\endaligned
\]

\[
\!\!\!\!\!\!\!\!\!\!
\!\!\!\!\!\!\!\!\!\!\!\!\!\!\!\!\!\!\!\!
\!\!\!\!\!\!\!\!\!\!\!\!\!\!\!\!\!\!\!\!
\aligned
\Big[
\mathcal{L},\,\,
\Big[
\mathcal{L},\,
\big[\overline{\mathcal{L}},\overline{\mathcal{S}}\big]
\Big]
\Big]
&
=
\mathcal{T}
\cdot
\left(
\aligned
\mathcal{L}\big(\mathcal{L}\big(\overline{P}\big)\big)
&
+
\overline{R}\,\mathcal{L}(P)
+
P\,\mathcal{L}\big(\overline{R}\big)
+
\overline{Q}\,\mathcal{L}(A)
+
A\,\mathcal{L}\big(\overline{Q}\big)
+
\\
&
+
P\,\mathcal{L}\big(\overline{R}\big)
+
A\,\mathcal{L}\big(\overline{Q}\big)
+
\\
&
+
P\overline{P}
+
PQ\overline{R}
+
BP\overline{Q}
+
AR\overline{R}
+
A\overline{B}\overline{Q}
\endaligned
\right)
+
\\
&
+
\mathcal{S}
\cdot
\left(
\aligned
\mathcal{L}\big(\mathcal{L}\big(\overline{R}\big)\big)
&
+
\mathcal{L}\big(\overline{P}\big)
+
\overline{R}\,\mathcal{L}(Q)
+
Q\,\mathcal{L}\big(\overline{R}\big)
+
\overline{Q}\,\mathcal{L}(B)
+
B\,\mathcal{L}\big(\overline{Q}\big)
+
\\
&
+
\mathcal{L}\big(\overline{P}\big)
+
Q\,\mathcal{L}\big(\overline{R}\big)
+
B\,\mathcal{L}\big(\overline{Q}\big)
+
\\
&
+
P\overline{R}
+
A\overline{Q}
+
\overline{P}Q
+
QQ\overline{R}
+
BQ\overline{Q}
+
BR\overline{R}
+
B\overline{B}\overline{Q}
\endaligned
\right)
+
\\
&
+
\overline{\mathcal{S}}
\cdot
\left(
\aligned
\mathcal{L}\big(\mathcal{L}\big(\overline{Q}\big)\big)
&
+
\overline{R}\,\mathcal{L}(R)
+
R\,\mathcal{L}\big(\overline{R}\big)
+
\overline{Q}\,\mathcal{L}\big(\overline{B}\big)
+
\overline{B}\,\mathcal{L}\big(\overline{Q}\big)
+
\\
&
+
R\,\mathcal{L}\big(\overline{R}\big)
+
\overline{B}\,\mathcal{L}\big(\overline{Q}\big)
+
\\
&
+
\overline{P}R
+
QR\overline{R}
+
B\overline{Q}R
+
\overline{B}R\overline{R}
+
\overline{B}\overline{B}\overline{Q}
\endaligned
\right).
\endaligned
\]
Adding these four expressions, one obtains
${\sf J}_{\sf rpl}$ and $K_{\sf rpl}$.
\endproof

\medskip\noindent{\bf Summary.} One has the 10 Lie bracket relations:
\[
\boxed{\,\,
\aligned
\big[\overline{\mathcal{S}},\mathcal{S}\big]
&
=
\overline{K}_{\sf rpl}\cdot\overline{\mathcal{S}}
-
K_{\sf rpl}\cdot\mathcal{S}
-
\isqrt\,J_{\sf rpl}\cdot\mathcal{T},
\\
\big[\overline{\mathcal{S}},\mathcal{T}\big]
&
=
-\,
\overline{F}_{\sf rpl}\cdot\overline{\mathcal{S}}
-
\overline{G}_{\sf rpl}\cdot\mathcal{S}
-
\overline{E}_{\sf rpl}\cdot\mathcal{T},
\\
\big[\overline{\mathcal{S}},\overline{\mathcal{L}}\big]
&
=
-\,
\overline{Q}\cdot\overline{\mathcal{S}}
-
\overline{R}\cdot\mathcal{S}
-
\overline{P}\cdot\mathcal{T},
\\
\big[\overline{\mathcal{S}},\mathcal{L}\big]
&
=
-\,
\overline{B}\cdot\overline{\mathcal{S}}
-
B\cdot\mathcal{S}
-
A\cdot\mathcal{T},
\\
\big[\mathcal{S},\,\mathcal{T}\big]
&
=
-\,
G_{\sf rpl}\cdot\overline{\mathcal{S}}
-
F_{\sf rpl}\cdot\mathcal{S}
-
E_{\sf rpl}\cdot\mathcal{T},
\\
\big[\mathcal{S},\,\overline{\mathcal{L}}\big]
&
=
-\,
\overline{B}\cdot\overline{\mathcal{S}}
-
B\cdot\mathcal{S}
-
A\cdot\mathcal{T},
\\
\big[\mathcal{S},\,\mathcal{L}\big]
&
=
-\,
R\cdot\overline{\mathcal{S}}
-
Q\cdot\mathcal{S}
-
P\cdot\mathcal{T},
\\
\big[\mathcal{T},\,\overline{\mathcal{L}}\big]
&
=
-\,
\overline{\mathcal{S}},
\\
\big[\mathcal{T},\,\mathcal{L}\big]
&
=
\,-\,\mathcal{S},
\\
\big[\overline{\mathcal{L}},\,\mathcal{L}\big]
&
=
\isqrt\,\mathcal{T}.\,\,
\endaligned}
\]

\medskip\noindent{\bf Initial Darboux structure of the 
dual coframe.} The coframe:
\[
\big\{du_3,du_2,du_1,\,dz,\,d\overline{z}\big\}
\]
is clearly dual to the frame: 
\[
\Big\{ 
{\textstyle{\frac{\partial}{\partial u_3}}},\, 
{\textstyle{\frac{\partial}{\partial u_2}}},\,
{\textstyle{\frac{\partial}{\partial u_1}}},\,
{\textstyle{\frac{\partial}{\partial z}}},\,
{\textstyle{\frac{\partial}{\partial\overline{z}}}} 
\Big\}.
\]
Introduce then the coframe:
\[
\big\{
\overline{\sigma_0},\,
\sigma_0,\,
\rho_0,\,
\overline{\zeta_0},\,
\zeta_0\big\}
\]
which is dual to the frame:
\[
\big\{
\overline{\mathcal{S}},\,
\mathcal{S},\,
\mathcal{T},\,
\overline{\mathcal{L}},\,
\mathcal{L}\big\},
\]
namely:
\[
\begin{array}{ccccc}
\overline{\sigma_0}(\overline{\mathcal{S}})=1 \ \ \ & \ \ \
\overline{\sigma_0}(\mathcal{S})=0 \ \ \ & \ \ \
\overline{\sigma_0}(\mathcal{T})=0 \ \ \ & \ \ \
\overline{\sigma_0}(\overline{\mathcal{L}})=0 \ \ \ & \ \ \
\overline{\sigma_0}\big(\mathcal{L}\big)=0,
\\
\sigma_0(\overline{\mathcal{S}})=0 \ \ \ & \ \ \
\sigma_0(\mathcal{S})=1 \ \ \ & \ \ \ \sigma_0(\mathcal{T})=0 \ \ \ &
\ \ \ \sigma_0(\overline{\mathcal{L}})=0 \ \ \ & \ \ \
\sigma_0\big(\mathcal{L}\big)=0,
\\
\rho_0(\overline{\mathcal{S}})=0 \ \ \ & \ \ \ \rho_0(\mathcal{S})=0
\ \ \ & \ \ \ \rho_0(\mathcal{T})=1 \ \ \ & \ \ \
\rho_0(\overline{\mathcal{L}})=0 \ \ \ & \ \ \ \rho_0(\mathcal{L})=0,
\\
\overline{\zeta_0}(\overline{\mathcal{S}})=0 \ \ \ & \ \ \
\overline{\zeta_0}(\mathcal{S})=0 \ \ \ & \ \ \
\overline{\zeta_0}(\mathcal{T})=0 \ \ \ & \ \ \
\overline{\zeta_0}(\overline{\mathcal{L}})=1 \ \ \ & \ \ \
\overline{\zeta_0}\big(\mathcal{L}\big)=0,
\\
\zeta_0(\overline{\mathcal{S}})=0 \ \ \ & \ \ \
\zeta_0(\mathcal{S})=0 \ \ \ & \ \ \ \zeta_0(\mathcal{T})=0 \ \ \ & \
\ \ \zeta_0(\overline{\mathcal{L}})=0 \ \ \ & \ \ \
\zeta_0\big(\mathcal{L}\big)=1.
\end{array}
\]
One has:
\[
\zeta_0 
= 
dz 
\ \ \ \ \ \ \ \ \ \ \ \ \ 
\text{\rm and} 
\ \ \ \ \ \ \ \ \ \ \ \ \ 
\overline{\zeta_0} 
= 
d\overline{z}.
\]

Organize the ten Lie brackets as a 
convenient auxiliary array:
\[
\footnotesize
\begin{array}{cccccccccccc}
& & \overline{\mathcal{S}} & & \mathcal{S} & & \mathcal{T} & &
\overline{\mathcal{L}} & & \mathcal{L}
\\
& & \boxed{d\overline{\sigma_0}} & & \boxed{d\sigma_0} & &
\boxed{d\rho_0} & & \boxed{d\overline{\zeta_0}} & & \boxed{d\zeta_0}
\\
\big[\overline{\mathcal{S}},\,\mathcal{S}\big] & = &
\overline{K}_{\sf rpl}\cdot\overline{\mathcal{S}} & + &
-\,K_{\sf rpl}\cdot\mathcal{S} & + & -\isqrt\,J_{\sf rpl}\cdot\mathcal{T} & + & 0 &
+ & 0 & \boxed{\overline{\sigma_0}\wedge\sigma_0}
\\
\big[\overline{\mathcal{S}},\,\mathcal{T}\big] & = &
-\,\overline{F}_{\sf rpl}\cdot\overline{\mathcal{S}} & + &
-\,\overline{G}_{\sf rpl}\cdot\mathcal{S} & + &
-\,\overline{E}_{\sf rpl}\cdot\mathcal{T} & + & 0 & + & 0 &
\boxed{\overline{\sigma_0}\wedge\rho_0}
\\
\big[\overline{\mathcal{S}},\,\overline{\mathcal{L}}\big] & = &
-\,\overline{Q}\cdot\overline{\mathcal{S}} & + &
-\,\overline{R}\cdot\mathcal{S} & + & -\,\overline{P}\cdot\mathcal{T}
& + & 0 & + & 0 & \boxed{\overline{\sigma_0}\wedge\overline{\zeta_0}}
\\
\big[\overline{\mathcal{S}},\,\mathcal{L}\big] & = &
-\,\overline{B}\cdot\overline{\mathcal{S}} & + &
-\,{B}\cdot\mathcal{S} & + & -A\cdot\mathcal{T} & + & 0 & + & 0 &
\boxed{\overline{\sigma_0}\wedge\zeta_0}
\\
\big[\mathcal{S},\,\mathcal{T}\big] & = &
-\,G_{\sf rpl}\cdot\overline{\mathcal{S}} & + & -\,F_{\sf rpl}\cdot 
\mathcal{S} & +
& -\,E_{\sf rpl}\cdot\mathcal{T} & + & 0 & + & 0 &
\boxed{\sigma_0\wedge\rho_0}
\\
\big[\mathcal{S},\,\overline{\mathcal{L}}\big] & = &
-\,\overline{B}\cdot\overline{\mathcal{S}} & + & -\,B\cdot
\mathcal{S} & + & -\,A\cdot\mathcal{T} & + & 0 & + & 0 &
\boxed{\sigma_0\wedge\overline{\zeta_0}}
\\
\big[\mathcal{S},\,\mathcal{L}\big] & = &
-\,R\cdot\overline{\mathcal{S}} & + & -\,Q\cdot\mathcal{S} & + &
-\,P\cdot\mathcal{T} & + & 0 & + & 0 & \boxed{\sigma_0\wedge\zeta_0}
\\
\big[\mathcal{T},\,\overline{\mathcal{L}}\big] & = &
-\,\overline{\mathcal{S}} & + & 0 & + & 0 & + & 0 & + & 0 &
\boxed{\rho_0\wedge\overline{\zeta_0}}
\\
\big[\mathcal{T},\,\mathcal{L}\big] & = & 0 & + & -\,\mathcal{S} & +
& 0 & + & 0 & + & 0 & \boxed{\rho_0\wedge\zeta_0}
\\
\big[\overline{\mathcal{L}},\,\mathcal{L}\big] & = & 0 & + & 0 & + &
\isqrt\,\mathcal{T} & + & 0 & + & 0 &
\boxed{\overline{\zeta_0}\wedge\zeta_0}
\end{array}
\]

Read {\em vertically} and put an overall minus sign:
\[
\boxed{ 
\aligned d\overline{\sigma}_0 & =
-\,\overline{K}_{\sf rpl} \cdot \overline{\sigma}_0\wedge\sigma_0 +
\overline{F}_{\sf rpl}\cdot \overline{\sigma}_0\wedge\rho_0 +
\overline{Q}\cdot \overline{\sigma}_0\wedge\overline{\zeta}_0 +
\overline{B}\cdot \overline{\sigma}_0\wedge\zeta_0 +
\\
& \ \ \ \ \ + G_{\sf rpl}\cdot \sigma_0\wedge\rho_0 + \overline{B}\cdot
\sigma_0\wedge\overline{\zeta}_0 + R\cdot \sigma_0\wedge\zeta_0 +
\rho_0\wedge\overline{\zeta}_0,
\\
d\sigma_0 & = K_{\sf rpl}\cdot \overline{\sigma}_0\wedge\sigma_0 +
\overline{G}_{\sf rpl}\cdot \overline{\sigma}_0\wedge\rho_0 +
\overline{R}\cdot \overline{\sigma}_0\wedge\overline{\zeta}_0 +
{B}\cdot \overline{\sigma}_0\wedge\zeta_0 +
\\
& \ \ \ \ \ + F_{\sf rpl}\cdot \sigma_0\wedge\rho_0 + B\cdot
\sigma_0\wedge\overline{\zeta}_0 + Q\cdot \sigma_0\wedge\zeta_0 +
\rho_0\wedge\zeta_0,
\\
d\rho_0 & = \isqrt\,J_{\sf rpl}\cdot \overline{\sigma}_0\wedge\sigma_0 +
\overline{E}_{\sf rpl}\cdot \overline{\sigma}_0\wedge\rho_0 +
\overline{P}\cdot \overline{\sigma}_0\wedge\overline{\zeta}_0 +
A\cdot \overline{\sigma}_0\wedge\zeta_0 +
\\
& \ \ \ \ \ + E_{\sf rpl}\cdot \sigma_0\wedge\rho_0 + A\cdot
\sigma_0\wedge\overline{\zeta}_0 + P\cdot \sigma_0\wedge\zeta_0 -
\isqrt\,\overline{\zeta}_0\wedge\zeta_0,
\\
d\overline{\zeta}_0 & = 0,
\\
d\zeta_0 & = 0.
\endaligned}
\]


\bigskip

\section{\sf $M^5 \subset \C^4$ of general class 
$\text{\sf III}_{\text{\sf 2}}$: 
\\
initial frame and coframe in local coordinates}
\label{initial-III-2}
\HEAD{\ref{initial-III-2}.~$M^5 \subset \C^4$ of general class 
$\text{\sf III}_{\text{\sf 2}}$
initial frame and coframe in local coordinates}{
Jo\"el {\sc Merker}, D\'epartement de Math\'ematiques d'Orsay}

\medskip

Next, consider:
\[
\Big(
M^5
\,\subset\,
\C^4
\Big)
\,\,\in\,\,
\text{\sf General Class $\text{\sf III}_{\text{\sf 2}}$}.
\]

As for the class $\text{\sf III}_{\text{\sf 1}}$, represent $M$ in
coordinates:
\[
(z,w_1,w_2,w_3) 
= 
\big(
x+\isqrt\,y,\, 
u_1+\isqrt\,v_1,\,
u_2+\isqrt\,v_2,\,
u_3+\isqrt\,v_3\big),
\]
as a graph:
\[
\aligned
v_1
&
=
\varphi_1(x,y,u_1,u_2,u_3),
\\
v_2
&
=
\varphi_2(x,y,u_1,u_2,u_3),
\\
v_3
&
=
\varphi_3(x,y,u_1,u_2,u_3),
\endaligned
\]

One assumes at every point that the following
biholomorphically invariant geometric condition holds:
\[
\boxed{\,\,
\aligned
{\bf 3}
&
=
\rank_\C\Big(T^{1,0}M,\,T^{0,1}M,\,
\big[T^{1,0}M,\,T^{0,1}M\big]
\Big),
\\
{\bf 4}
&
=
\rank_\C\Big(T^{1,0}M,\,T^{0,1}M,\,
\big[T^{1,0}M,\,T^{0,1}M\big],
\\
&
\ \ \ \ \ \ \ \ \ \ \ \ \ \ \
\big[T^{1,0}M,\,\big[T^{1,0}M,\,T^{0,1}M\big]\big]
\Big),
\\
{\bf 4}
&
=
\rank_\C\Big(T^{1,0}M,\,T^{0,1}M,\,
\big[T^{1,0}M,\,T^{0,1}M\big],
\\
&
\ \ \ \ \ \ \ \ \ \ \ \ \ \ \
\big[T^{1,0}M,\,\big[T^{1,0}M,\,T^{0,1}M\big]\big],
\big[T^{0,1}M,\,\big[T^{1,0}M,\,T^{0,1}M\big]\big]
\Big),
\\
{\bf 5}
&
=
\rank_\C\Big(T^{1,0}M,\,T^{0,1}M,\,
\big[T^{1,0}M,\,T^{0,1}M\big],
\big[T^{1,0}M,\,\big[T^{1,0}M,\,T^{0,1}M\big]\big],\,\,
\\
&
\ \ \ \ \ \ \ \ \ \ \ \ \ \ \
\big[T^{1,0}M,\,\big[T^{1,0}M,\,\big[T^{1,0}M,\,T^{0,1}M\big]\big]\big]
\Big),
\endaligned}
\]
the third rank condition being an exceptional degeneracy
feature, because here, $5$ fields have only rank $4$.

\medskip

There is a unique (local) generator of $T^{1,0}M$ of the form:
\[
\mathcal{L}
=
\frac{\partial}{\partial z}
+
\frac{\Lambda_1}{\Delta}\,\frac{\partial}{\partial u_1}
+
\frac{\Lambda_2}{\Delta}\,\frac{\partial}{\partial u_2}
+
\frac{\Lambda_3}{\Delta}\,\frac{\partial}{\partial u_3},
\]
having conjugate:
\[
\overline{\mathcal{L}}
=
\frac{\partial}{\partial\overline{z}}
+
\frac{\overline{\Lambda}_1}{\overline{\Delta}}\,\frac{\partial}{\partial u_1}
+
\frac{\overline{\Lambda_2}}{\overline{\Delta}}\,\frac{\partial}{\partial u_2}
+
\frac{\overline{\Lambda_3}}{\overline{\Delta}}\,\frac{\partial}{\partial u_3},
\]
with coefficient-functions that have exactly the
same expressions as in the class
$\text{\sf III}_{\text{\sf 1}}$ treated
in the preceding section.

\medskip

Furthermore, the three fields\,\,---\,\,the first
$\overline{\mathcal{ T}} = \mathcal{ T}$ being real\,\,---:
\[
\aligned
\mathcal{T}
&
:=
\isqrt\,\big[\mathcal{L},\overline{\mathcal{L}}\big],
\\
\mathcal{S}
&
:=
\big[\mathcal{L},\mathcal{T}\big]
=
\big[\mathcal{L},\,\isqrt\,\big[\mathcal{L},\overline{\mathcal{L}}\big]\big],
\\
\mathcal{R}
&
:=
\big[\mathcal{L},\mathcal{S}\big]
=
\Big[\mathcal{L},\,
\big[\mathcal{L},\,\isqrt\,\big[\mathcal{L},\overline{\mathcal{L}}\big]\big]
\Big]
\endaligned
\]
are of the form:
\[
\aligned
\mathcal{T}
&
=
\frac{\Upsilon_1}{\Delta^2\,\overline{\Delta}^2}\,
\frac{\partial}{\partial u_1}
+
\frac{\Upsilon_2}{\Delta^2\,\overline{\Delta}^2}\,
\frac{\partial}{\partial u_2}
+
\frac{\Upsilon_3}{\Delta^2\,\overline{\Delta}^2}\,
\frac{\partial}{\partial u_3},
\\
\mathcal{S}
&
=
\frac{\Pi_1}{\Delta^4\,\overline{\Delta}^3}\,
\frac{\partial}{\partial u_1}
+
\frac{\Pi_2}{\Delta^4\,\overline{\Delta}^3}\,
\frac{\partial}{\partial u_2}
+
\frac{\Pi_3}{\Delta^4\,\overline{\Delta}^3}\,
\frac{\partial}{\partial u_3},
\\
\mathcal{R}
&
=
\frac{\Sigma_1}{\Delta^6\,\overline{\Delta}^4}\,
\frac{\partial}{\partial u_1}
+
\frac{\Sigma_2}{\Delta^6\,\overline{\Delta}^4}\,
\frac{\partial}{\partial u_2}
+
\frac{\Sigma_3}{\Delta^6\,\overline{\Delta}^4}\,
\frac{\partial}{\partial u_3},
\endaligned
\]
with coefficient-functions that depend explicitly upon 
$\varphi_1, \varphi_2, \varphi_3$, although
computers seem too weak to do that.

One therefore {\em decides} not to 
keep an explicit track of these dependencies.

One also has the conjugate field:
\[
\aligned
\overline{\mathcal{S}}
&
=
\overline{\big[\mathcal{L},\mathcal{T}\big]}
\\
&
=
\big[\overline{\mathcal{L}},\,\mathcal{T}\big]
\\
&
=
\frac{\overline{\Pi}_1}{\Delta^3\,\overline{\Delta}^4}\,
\frac{\partial}{\partial u_1}
+
\frac{\overline{\Pi}_2}{\Delta^3\,\overline{\Delta}^4}\,
\frac{\partial}{\partial u_2}
+
\frac{\overline{\Pi}_3}{\Delta^3\,\overline{\Delta}^4}\,
\frac{\partial}{\partial u_3}.
\endaligned
\]

\medskip\noindent{\bf Reformulation of the geometric assumptions:}
\[
\aligned
{\bf 3}
&
=
\rank_\C
\big(\mathcal{L},\,\overline{\mathcal{L}},\,
\mathcal{T}\big),
\\
{\bf 4}
&
=
\rank_\C
\big(\mathcal{L},\,\overline{\mathcal{L}},\,
\mathcal{T},\,\mathcal{S}\big),
\\
{\bf 4}
&
=
\rank_\C
\big(\mathcal{L},\,\overline{\mathcal{L}},\,
\mathcal{T},\,\mathcal{S},\,\overline{\mathcal{S}}\big),
\\
{\bf 5}
&
=
\rank_\C
\big(\mathcal{L},\,\overline{\mathcal{L}},\,
\mathcal{T},\,\mathcal{S},\,\mathcal{R}\big).
\endaligned
\]

Consequently, there are $2$ functions $A$, $B$ so that:
\[
\boxed{\,
\aligned
\overline{\mathcal{S}}
&
=
A\cdot\mathcal{T}
+
B\cdot\mathcal{S}\,
\\
&
=
\big[
\overline{\mathcal{L}},\mathcal{T}
\big].
\endaligned}
\]

Hence:
\[
\aligned
\mathcal{S}
=
\big[\mathcal{L},\mathcal{T}\big]
=
\overline{\big[\overline{\mathcal{L}},\mathcal{T}\big]}
&
=
\overline{A}\cdot\mathcal{T}
+
\overline{B}\cdot\underbrace{\overline{\mathcal{S}}}_{\sf replace}
\\
&
=
\overline{A}\cdot\mathcal{T}
+
\overline{B}A\cdot\mathcal{T}
+
\overline{B}B\cdot\mathcal{S}
\\
&
=
\big(\overline{A}+\overline{B}A\big)\cdot\mathcal{T}
+
B\overline{B}\cdot\mathcal{S},
\endaligned
\]
whence by identification:
\[
\boxed{\,\,
\aligned
0
&
=
\overline{A}
+
\overline{B}A,\,\,
\\
1
&
=
B\overline{B}.
\endaligned}
\]

Next, remind that the Jacobi identity gives:
\[
\Big[\mathcal{L},\big[\overline{\mathcal{L}},\,
\isqrt\,\big[\mathcal{L},\overline{\mathcal{L}}\big]\big]\Big]
=
\Big[\overline{\mathcal{L}},\big[\mathcal{L},\,
\isqrt\,\big[\mathcal{L},\overline{\mathcal{L}}\big]\big]\Big],
\]
that is to say:
\[
\big[\overline{\mathcal{L}},\mathcal{S}\big]
=
\big[\mathcal{L},\overline{\mathcal{S}}\big].
\]

Consequently:
\[
\aligned
\big[\overline{\mathcal{L}},\mathcal{S}\big]
=
\big[\mathcal{L},\overline{\mathcal{S}}\big]
&
=
\big[\mathcal{L},A\mathcal{T}+B\mathcal{S}\big]
\\
&
=
\mathcal{L}(A)\cdot\mathcal{T}
+
\mathcal{L}(B)\cdot\mathcal{S}
+
A\underbrace{\big[\mathcal{L},\mathcal{T}\big]}_{\mathcal{S}}
+
B\underbrace{\big[\mathcal{L},\mathcal{S}\big]}_{\mathcal{R}}
\\
&
=
\mathcal{L}(A)\cdot\mathcal{T}
+
\big(\mathcal{L}(B)+A\big)\cdot\mathcal{S}
+
B\cdot\mathcal{R},
\endaligned
\]
{\em i.e.}:
\[
\boxed{\,\,
\big[\mathcal{S},\overline{\mathcal{L}}\big]
=
-\,\mathcal{L}(A)\cdot\mathcal{T}
-
\big(\mathcal{L}(B)+A\big)\cdot\mathcal{S}
-
B\cdot\mathcal{R}.\,\,
}
\]

The $2$ functions $A$, $B$ will be fundamental, plus $3$ other
functions $E$, $F$, $G$ only in:
\[
\boxed{\,\,
\aligned
\big[\mathcal{L},\mathcal{R}\big]
&
=
E\cdot\mathcal{T}
+
F\cdot\mathcal{S}
+
G\cdot\mathcal{R},
\endaligned}
\]
because the other Lie brackets will be expressed in terms of:
\[
A,\ \ \ \ \ 
B,\ \ \ \ \
E,\ \ \ \ \
F,\ \ \ \ \
G.
\]

\medskip

Indeed, between the {\bf 5} elements:
\[
\big\{
\mathcal{R},\mathcal{S},\mathcal{T},\overline{\mathcal{L}},\mathcal{L}
\big\},
\]
we yet have to determine {\bf 4} amongst 
all the {\bf 10} possible Lie brackets
\[
\aligned
\big[\mathcal{R},\mathcal{S}\big]
&
=
\something,
\\
\big[\mathcal{R},\mathcal{T}\big]
&
=
\something,
\\
\big[\mathcal{R},\overline{\mathcal{L}}\big]
&
=
\something,
\\
\big[\mathcal{R},\mathcal{L}\big]
&
=
-\,E\cdot\mathcal{T}
-
F\cdot\mathcal{S}
-
G\cdot\mathcal{R},
\\
\big[\mathcal{S},\mathcal{T}\big]
&
=
\something,
\\
\big[\mathcal{S},\overline{\mathcal{L}}\big]
&
=
-\,\mathcal{L}(A)\cdot\mathcal{T}
-
\big(\mathcal{L}(B)+A\big)\cdot\mathcal{S}
-
B\cdot\mathcal{R},
\\
\big[\mathcal{S},\mathcal{L}\big]
&
=
-\,\mathcal{R},
\\
\big[\mathcal{T},\overline{\mathcal{L}}\big]
&
=
-\,A\cdot\mathcal{T}
-
B\cdot\mathcal{S},
\\
\big[\mathcal{T},\mathcal{L}\big]
&
=
-\,\mathcal{S},
\\
\big[\overline{\mathcal{L}},\mathcal{L}\big]
&
=
i\,\mathcal{T},
\endaligned
\]
and these {\bf 4} remaining will express in terms of:
\[
A,\ \ \ \ \
B,\ \ \ \ \
E,\ \ \ \ \
F,\ \ \ \ \
G.
\]

\medskip\noindent{\bf Preparatory lemma.}
{\em The conjugate $\overline{ R}$ expresses as:}
\[
\aligned
\overline{R}
&
=
\big(
\overline{\mathcal{L}}(A)
+
B\,\mathcal{L}(A)
+
AA
\big)
\cdot
\mathcal{T}
+
\\
&
+
\big(
\overline{\mathcal{L}}(B)
+
B\,\mathcal{L}(B)
+
2\,AB
\big)
\cdot
\mathcal{S}
+
\\
&
+
\big(BB\big)
\cdot
\mathcal{R}.
\endaligned
\]

\proof
Indeed, from:
\[
\mathcal{R}
=
\big[\mathcal{L},\mathcal{S}\big],
\]
by conjugating:
\[
\aligned
\overline{\mathcal{R}}
&
=
\big[\overline{\mathcal{L}},\,\overline{\mathcal{S}}\big]
\\
&
=
\big[\overline{\mathcal{L}},\,\,
A\,\mathcal{T}
+
B\,\mathcal{S}
\big]
\\
&
=
\overline{\mathcal{L}}(A)
\cdot
\mathcal{T}
+
\overline{\mathcal{L}}(B)
\cdot
\mathcal{S}
+
\\
&
\ \ \ \ \ 
+
A\,
\underbrace{
\big[\overline{\mathcal{L}},\,\mathcal{T}\big]}_{
A\mathcal{T}+B\mathcal{S}}
+
B\,
\underbrace{
\big[\overline{\mathcal{L}},\,\mathcal{S}\big]}_{
\substack{
\mathcal{L}(A)\mathcal{T}+
\\
(\mathcal{L}(B)+A)\mathcal{S}
\\
+
BR}},
\endaligned
\]
which, collecting, is so.
\endproof

\noindent{\bf Lemma.}
{\em The {\bf 3} functions in the Lie bracket:}
\[
\big[\overline{\mathcal{L}},\,\mathcal{R}\big]
\,=:\,
H_{\sf rpl}\cdot\mathcal{T}
+
J_{\sf rpl}\cdot\mathcal{S}
+
K_{\sf rpl}\cdot\mathcal{R}
\]
express in terms of $A$, $B$, $E$, $F$, $G$ as:
\[
\aligned
K_{\sf rpl}
&
=
2\,B\,\overline{\mathcal{L}}\big(\overline{B}\big)
+
BB\,\mathcal{L}\big(\overline{B}\big)
+
B\,\overline{\mathcal{L}}\big(\overline{B}\big)
+
\\
&
\ \ \ \ \
+
2\,AB
+
\overline{G},
\endaligned
\]
\[
\aligned
J_{\sf rpl}
&
=
B\,\overline{\mathcal{L}}\big(\mathcal{L}\big(\overline{B}\big)\big)
+
\overline{\mathcal{L}}\big(\overline{\mathcal{L}}
\big(\overline{B}\big)\big)
\,-
\\
&
\ \ \ \ \
-\,2\,BB\,\mathcal{L}\big(\overline{B}\big)\,
\overline{\mathcal{L}}\big(\overline{B}\big)
-
BBB\,\mathcal{L}\big(\overline{B}\big)\,
\mathcal{L}\big(\overline{B}\big)
-
2\,BB\,\mathcal{L}\big(\overline{B}\big)\,
\overline{\mathcal{L}}\big(\overline{B}\big)
\,-
\\
&
\ \ \ \ \
-\,
2\,ABB\,\mathcal{L}\big(\overline{B}\big)
-
B\,\mathcal{L}\big(\overline{B}\big)\,\overline{G}
-
2\,B\,\overline{\mathcal{L}}\big(\overline{B}\big)\,
\overline{\mathcal{L}}\big(\overline{B}\big)
-
\overline{G}\,\overline{\mathcal{L}}\big(\overline{B}\big)
\,-
\\
&
\ \ \ \ \ 
-\,
4\,AB\,\overline{\mathcal{L}}\big(\overline{B}\big)
-
2\,ABB\,\mathcal{L}\big(\overline{B}\big)
-
2\,AB\,\overline{\mathcal{L}}\big(\overline{B}\big)
+
3\,\overline{\mathcal{L}}(A)
+
B\,\mathcal{L}(A)
+
\\
&
\ \ \ \ \
+
\overline{B}\overline{F}
-
3\,AAB
-
2\,A\overline{G},
\endaligned
\]
and, without full replacements:
\[
\aligned
H_{\sf rpl}
&
=
-\,A\,J_{\sf rpl}
-
K_{\sf rpl}\,
\mathcal{L}(A)
-
K_{\sf rpl}\,\overline{B}\,\overline{\mathcal{L}}(A)
+
K_{\sf rpl}\,AA
+
\\
&
\ \ \ \ \
+
\mathcal{L}\big(\overline{\mathcal{L}}(A)\big)
+
B\,\mathcal{L}\big(\mathcal{L}(A)\big)
+
\\
&
\ \ \ \ \
+
\mathcal{L}(B)\,\mathcal{L}(A)
+
2\,A\,\mathcal{L}(A)
+
BBE.
\endaligned
\]

\proof
Start with computing:
\[
\!\!\!\!\!\!\!\!\!\!\!\!\!\!\!\!\!\!\!\!
\!\!\!\!\!\!\!\!\!\!\!\!\!\!\!\!\!\!\!\!
\aligned
\big[\mathcal{L},\overline{\mathcal{R}}\big]
=
\Big[
\mathcal{L},\,\,
&
\big(
\overline{\mathcal{L}}(A)
+
B\,\mathcal{L}(A)
+
AA
\big)
\cdot
\mathcal{T}
+
\\
&
+
\big(
\overline{\mathcal{L}}(B)
+
B\,\mathcal{L}(B)
+
2\,AB
\big)
\cdot
\mathcal{S}
+
\\
&
+
\big(BB\big)
\cdot
\mathcal{R}
\ \ \ \ \ \ \ \ \ \ \ \ \ \ \ \ \ \ \ \ \ \ \ \ \ \ \ \ \ 
\ \ \ \ \ \ \ \ \ \ \
\Big]
\\
&
=
\Big(
\mathcal{L}\big(\overline{\mathcal{L}}(A)\big)
+
B\,\mathcal{L}\big(\mathcal{L}(A)\big)
+
\mathcal{L}(B)\,\mathcal{L}(A)
+
2\,A\,\mathcal{L}(A)
\Big)
\cdot\mathcal{T}
+
\\
&
\ \ \ \ \
+
\Big(
\mathcal{L}\big(\overline{\mathcal{L}}(B)\big)
+
B\,\mathcal{L}\big(\mathcal{L}(B)\big)
+
\mathcal{L}(B)\,\mathcal{L}(B)
+
2\,B\,\mathcal{L}(A)
+
2\,A\,\mathcal{L}(B)
\Big)
\cdot\mathcal{S}
+
\\
&
\ \ \ \ \
+
\big(2\,B\,\mathcal{L}(B)\big)
\cdot
\mathcal{R}
+
\\
&
\ \ \ \ \
+
\big(
\overline{\mathcal{L}}(A)
+
B\,\mathcal{L}(A)
+
AA
\big)
\cdot
\underbrace{
\big[\mathcal{L},\mathcal{T}\big]}_{\mathcal{S}}
+
\\
&
\ \ \ \ \
+
\big(\overline{\mathcal{L}}(B)
+
B\,\mathcal{L}(B)
+
2\,AB
\big)
\cdot
\underbrace{
\big[\mathcal{L},\mathcal{S}\big]}_{
\mathcal{R}}
+
\\
&
\ \ \ \ \
+
BB
\cdot\!\!\!\!\!
\underbrace{
\big[\mathcal{L},\mathcal{R}\big]}_{
E\mathcal{T}+F\mathcal{S}+G\mathcal{R}},
\endaligned
\]
which, collecting, becomes:

\[
\aligned
\big[\mathcal{L},\overline{\mathcal{R}}\big]
&
=
\mathcal{T}
\cdot
\left(
\aligned
&
\mathcal{L}\big(\overline{\mathcal{L}}(A)\big)
+
B\,\mathcal{L}\big(\mathcal{L}(A)\big)
+
\\
&
+
\mathcal{L}(B)\,\mathcal{L}(A)
+
2\,A\,\mathcal{L}(A)
+
BBE
\endaligned
\right)
+
\\
&
+
\mathcal{S}
\cdot
\left(
\aligned
&
\mathcal{L}\big(\overline{\mathcal{L}}(B)\big)
+
B\,\mathcal{L}\big(\mathcal{L}(B)\big)
+
\\
&
+
\mathcal{L}(B)\,\mathcal{L}(B)
+
2\,B\,\mathcal{L}(A)
+
2\,A\,\mathcal{L}(B)
+
\overline{\mathcal{L}}(A)
+
B\,\mathcal{L}(A)
+
\\
&
+
AA
+
BBF
\endaligned
\right)
+
\\
&
+
\mathcal{R}
\cdot
\left(
\aligned
&
2\,B\,\mathcal{L}(B)
+
\overline{\mathcal{L}}(B)
+
B\,\mathcal{L}(B)
+
\\
&
+
2\,AB
+
BBG
\endaligned
\right).
\endaligned
\]

On the other hand, conjugating:
\[
\big[
\overline{\mathcal{L}},\,\mathcal{R}
\big]
\,=:\,
H_{\sf rpl}\cdot\mathcal{T}
+
J_{\sf rpl}\cdot\mathcal{S}
+
K_{\sf rpl}\cdot\mathcal{R},
\]
one receives:
\[
\aligned
\big[\mathcal{L},\overline{R}\big]
&
=
\overline{
\big[\overline{\mathcal{L}},\mathcal{R}\big]}
\\
&
=
\overline{H}_{\sf rpl}
\cdot
\mathcal{T}
+
\overline{J}_{\sf rpl}
\cdot
\overline{\mathcal{S}}
+
\overline{K}_{\sf rpl}
\cdot
\overline{\mathcal{R}}
\\
&
=
\overline{H}_{\sf rpl}
\cdot
\mathcal{T}
+
\overline{J}_{\sf rpl}\,
\big(
A\cdot\mathcal{T}
+
B\cdot\mathcal{S}
\big)
+
\\
&
\ \ \ \ \
+
\overline{K}_{\sf rpl}
\Big\{
\big(
\overline{\mathcal{L}}(A)
+
B\,\mathcal{L}(A)
+
AA
\big)
\cdot
\mathcal{T}
+
\\
&
\ \ \ \ \ \ \ \ \ \ \ \ \ \ \ \ \ \ \
+
\big(
\overline{\mathcal{L}}(B)
+
B\,\mathcal{L}(B)
+
2\,AB
\big)
\cdot
\mathcal{S}
+
\\
&
\ \ \ \ \ \ \ \ \ \ \ \ \ \ \ \ \ \ \
+
\big(BB\big)\cdot\mathcal{R}
\ \ \ \ \ \ \ \ \ \ \ \ \ \ \ \ \ \ \ \ \ \ \ \ \ \ \ \ \ \ \ \ \ \ \
\ \ \ \ \ \
\Big\},
\endaligned
\]
that is to say:
\[
\aligned
\big[
\overline{\mathcal{L}},\,\mathcal{R}
\big]
&
=
\mathcal{T}
\cdot
\Big(
\overline{H}_{\sf rpl}
+
A\,\overline{J}_{\sf rpl}
+
\overline{K}_{\sf rpl}\,
\overline{\mathcal{L}}(A)
+
\overline{K}_{\sf rpl}\,
B\,\mathcal{L}(A)
+
\overline{K}_{\sf rpl}\,
AA
\Big)
+
\\
&
\ \ \ \ \
+
\mathcal{S}
\cdot
\Big(
B\,\overline{J}_{\sf rpl}
+
\overline{K}_{\sf rpl}\,
\overline{\mathcal{L}}(B)
+
\overline{K}_{\sf rpl}\,B\,\mathcal{L}(B)
+
2\,\overline{K}_{\sf rpl}\,AB
\Big)
+
\\
&
\ \ \ \ \
+
\mathcal{R}
\cdot
\Big(
BB\,\overline{K}_{\sf rpl}
\Big).
\endaligned
\]

Equate the coefficients of $\mathcal{ R}$
in these two expressions of $\big[ \mathcal{L}, 
\overline{ \mathcal{ R}} \big]$:
\[
\aligned
BB\,\overline{K}_{\sf rpl}
&
=
2\,B\,\mathcal{L}(B)
+
\overline{\mathcal{L}}(B)
+
B\,\mathcal{L}(B)
+
\\
&
\ \ \ \ \
+
2\,AB
+
BBG.
\endaligned
\] 
Recall:
\[
B\overline{B}
\equiv
1,
\]
hence multiply both sides by $\overline{ B} \overline{ B}$ to get:
\[
\aligned
\overline{K}_{\sf rpl}
&
=
2\,\overline{B}\,\mathcal{L}(B)
+
\overline{B}\overline{B}\,
\overline{\mathcal{L}}(B)
+
\overline{B}\,\mathcal{L}(B)
+
\\
&
\ \ \ \ \
+
2\,A\overline{B}
+
G.
\endaligned
\]
Conjugate this, and get $K_{\sf rpl}$ as stated.

Next, identify the coefficients of $\mathcal{ S}$
in the two expressions of $\big[ \mathcal{L},
\overline{ \mathcal{ R}} \big]$, and get:
\[
\aligned
B\,\overline{J}_{\sf rpl}
&
=
-\,
\overline{K}_{\sf rpl}\,\overline{\mathcal{L}}(B)
-
\overline{K}_{\sf rpl}\,B\,\mathcal{L}(B)
-
2\,\overline{K}_{\sf rpl}\,AB
+
\\
&
\ \ \ \ \
+
\mathcal{L}\big(\overline{\mathcal{L}}(B)\big)
+
B\,\mathcal{L}\big(\mathcal{L}(B)\big)
+
\\
&
\ \ \ \ \
+
\mathcal{L}(B)\,\mathcal{L}(B)
+
2\,B\,\mathcal{L}(A)
+
2\,A\,\mathcal{L}(B)
+
\overline{\mathcal{L}}(A)
+
B\,\mathcal{L}(A)
+
\\
&
\ \ \ \ \
+
AA
+
BBF.
\endaligned
\]
Multiply both sides by $\overline{ B}$, conjugate,
and get $J_{\sf rpl}$ are stated.

Lastly,
to get $H_{\sf rpl}$, identify the
conjugated coefficients of $\mathcal{ R}$.
\endproof

Thanks to all this:
\[
\!\!\!\!\!\!\!\!\!\!\!\!\!\!\!\!\!\!\!\!
\aligned
\big[\mathcal{S},\mathcal{T}\big]
&
=
\big[\mathcal{S},\,
\isqrt\,\big[\mathcal{L},\overline{\mathcal{L}}\big]
\\
&
=
\isqrt\,\big[\overline{\mathcal{L}},\big[\mathcal{L},\mathcal{S}\big]\big]
-
\isqrt\,\big[\mathcal{L},\big[\overline{\mathcal{L}},\mathcal{S}\big]\big]
\\
&
=
\isqrt\,\big[\overline{\mathcal{L}},\mathcal{R}\big]
-
\isqrt\,\big[\mathcal{L},\,
\mathcal{L}(A)\mathcal{T}
+
\big(\mathcal{L}(B)+A\big)\mathcal{S}
+
B\mathcal{R}
\big]
\\
&
=
\isqrt\,H_{\sf rpl}\cdot\mathcal{T}
+
\isqrt\,J_{\sf rpl}\cdot\mathcal{S}
+
\isqrt\,K_{\sf rpl}\cdot\mathcal{R}
-
\\
&
\ \ \ \ \
-\,
\isqrt\,\mathcal{L}\big(\mathcal{L}(A)\big)\cdot\mathcal{T}
-
\isqrt\,\mathcal{L}\big(\mathcal{L}(B)\big)\cdot\mathcal{S}
-
\isqrt\,\mathcal{L}(A)\cdot\mathcal{S}
-
\isqrt\,\mathcal{L}(B)\cdot\mathcal{R}
-
\\
&
\ \ \ \ \
-\,\isqrt\,\mathcal{L}(A)
\underbrace{\big[\mathcal{L},\mathcal{T}\big]}_{
\mathcal{S}}
-
\isqrt\,\big(\mathcal{L}(B)+A\big)
\underbrace{\big[\mathcal{L},\mathcal{R}\big]}_{
\mathcal{R}}
-
\isqrt\,B
\underbrace{\big[\mathcal{L},\mathcal{R}\big]}_{
E\mathcal{T}+F\mathcal{S}+G\mathcal{R}},
\endaligned
\]
which gives:
\[
\boxed{\,\,
\aligned
\big[\mathcal{S},\mathcal{T}\big]
&
=
\mathcal{T}\cdot\big(
-\,
\isqrt\,\mathcal{L}(\mathcal{L}(A))
+
\isqrt\,H_{\sf rpl}
-
\isqrt\,BE
\big)
+
\\
&
+
\mathcal{S}\cdot\big(
-
\isqrt\,\mathcal{L}(\mathcal{L}(B)\big)
-
2\isqrt\,\mathcal{L}(A)
+
\isqrt\,J_{\sf rpl}
-
\isqrt\,BF
\big)
+\,\,
\\
&
+
\mathcal{R}\cdot\big(
-
2\isqrt\,\mathcal{L}(B)
+
\isqrt\,K_{\sf rpl}
-
\,\isqrt\,A
-
\isqrt\,BG\big).
\endaligned}
\]

By similar reasonings based on the Jacobi identity
(exercise):
\[
\!\!\!\!\!\!\!\!\!\!\!\!\!\!\!\!\!\!\!\!
\boxed{\,\,
\aligned
\big[\mathcal{R},\mathcal{T}\big]
&
=
\mathcal{T}\cdot\big(
\isqrt\,\overline{\mathcal{L}}(E)
-
\isqrt\,\mathcal{L}\big(K_{\sf rpl}\big)
+
\isqrt\,F\mathcal{L}(A)
+
\isqrt\,AE
+
\\
&
\ \ \ \ \
+
\isqrt\,GH_{\sf rpl}
-
\isqrt\,EK_{\sf rpl}
\big)
+
\\
&
+
\mathcal{S}\cdot\big(
\isqrt\,\overline{\mathcal{L}}(F)
-
\isqrt\,L\big(J_{\sf rpl}\big)
+
\isqrt\,F\mathcal{L}(B)
+
\isqrt\,BE
+
\isqrt\,AF
+
\\
&
\ \ \ \ \ 
+
\isqrt\,GJ_{\sf rpl}
-
\isqrt\,H_{\sf rpl}
-
\isqrt\,FK_{\sf rpl}
\big)
+
\\
&
+
\mathcal{R}\cdot\big(
\isqrt\,\overline{\mathcal{L}}(G)
-
\isqrt\,\mathcal{L}\big(K_{\sf rpl}\big)
+
\isqrt\,BF
+
\isqrt\,GH_{\sf rpl}
-
\isqrt\,J_{\sf rpl}
-
\isqrt\,GH_{\sf rpl}
\big).
\endaligned}
\]

It is advisable to abbreviate the
just computed complicated Lie brackets
by naming functions:
\[
\aligned
&
L_{\sf rpl},\ \ \ \ \ \ \ \ \ \ \
M_{\sf rpl},\ \ \ \ \ \ \ \ \ \ \
N_{\sf rpl},
\\
&
O_{\sf rpl},\ \ \ \ \ \ \ \ \ \ \
P_{\sf rpl},\ \ \ \ \ \ \ \ \ \ \
Q_{\sf rpl},
\\
&
R_{\sf rpl},\ \ \ \ \ \ \ \ \ \ \
S_{\sf rpl},\ \ \ \ \ \ \ \ \ \ \
T_{\sf rpl},
\endaligned
\]
with:
\[
\aligned
\big[\mathcal{S},\mathcal{T}\big]
&
=
-\,
L\cdot\mathcal{T}
-
M\cdot\mathcal{S}
-
N\cdot\mathcal{R},
\\
\big[\mathcal{R},\mathcal{T}\big]
&
=
-\,O_{\sf rpl}\cdot\mathcal{T}
-
P_{\sf rpl}\cdot\mathcal{S}
-
Q_{\sf rpl}\cdot\mathcal{R},
\\
\big[\mathcal{R},\mathcal{S}\big]
&
=
-\,R_{\sf rpl}\cdot\mathcal{T}
-
S_{\sf rpl}\cdot\mathcal{S}
-
T_{\sf rpl}\cdot\mathcal{R};
\endaligned
\]
in fact lastly:
\[
\aligned
\big[\mathcal{R},\mathcal{S}\big]
&
=
\Big[
\big[\mathcal{L},\,
\big[\mathcal{L},\,\isqrt\,
\big[\mathcal{L},\overline{\mathcal{L}}\big]
\big]\big],
\,\,\,\,
\big[
\mathcal{L},\,\isqrt\,
\big[\mathcal{L},\overline{\mathcal{L}}\big]\big]
\Big]
\\
&
=
\big[
\mathcal{R},\,\big[\mathcal{L},\mathcal{T}\big]
\big]
\\
&
=
-\,
\big[\mathcal{T},\,
\big[\mathcal{R},\mathcal{L}\big]\big]
-
\big[\mathcal{L},\,
\big[\mathcal{T},\mathcal{R}\big]\big]
\\
&
=
\big[
\mathcal{L},\,
E\,\mathcal{T}+F\,\mathcal{S}+G\,\mathcal{R}
\big]
\,-
\\
&
\ \ \ \ \ 
-\,
\big[
\mathcal{L},\,
O_{\sf rpl}\,\mathcal{T}
+
P_{\sf rpl}\,\mathcal{S}
+
Q_{\sf rpl}\,\mathcal{R}
\big]
\\
&
=
\mathcal{T}(E)\cdot\mathcal{T}
+
\mathcal{T}(F)\cdot\mathcal{S}
+
\mathcal{T}(G)\cdot\mathcal{R}
+
\\
&
\ \ \ \ \
+
E\,
\zero{\big[\mathcal{T},\mathcal{T}\big]}
+
F\,
\underbrace{
\big[\mathcal{T},\mathcal{S}\big]}_{\sf replace}
+
G\,
\underbrace{
\big[\mathcal{T},\mathcal{R}\big]}_{\sf replace}
\,-
\\
&
-\,
\mathcal{L}\big(O_{\sf rpl}\big)
\cdot
\mathcal{T}
-
\mathcal{L}\big(P_{\sf rpl}\big)
\cdot
\mathcal{S}
-
\mathcal{L}\big(Q_{\sf rpl}\big)
\cdot
\mathcal{R}
\,-
\\
&
-\,
O_{\sf rpl}\,
\underbrace{
\big[\mathcal{L},\mathcal{S}\big]}_{\mathcal S}
-
P_{\sf rpl}\,
\underbrace{
\big[\mathcal{L},\mathcal{S}\big]}_{\mathcal{R}}
-\,
Q_{\sf rpl}
\!\!\!\!
\underbrace{
\big[\mathcal{L},\mathcal{R}\big]}_{
E\mathcal{T}+F\mathcal{S}+G\mathcal{R}},
\endaligned
\]
so that in:
\[
\big[\mathcal{R},\mathcal{S}\big]
=
-\,
R_{\sf rpl}\cdot\mathcal{T}
-
S_{\sf rpl}\cdot\mathcal{S}
-
T_{\sf rpl}\cdot\mathcal{R},
\]
the three appearing coefficient-functions indeed
express\,\,---\,\,lengthily if expanded which is not done here\,\,---\,\,in
terms of the five
fundamental functions $A$, $B$, $E$, $F$, $G$
and their $\big\{ \mathcal{L},
\overline{ \mathcal{L}} \big\}$-derivatives.

\medskip\noindent{\bf Summary.}
In the 10 Lie brackets, "${}_{\sf rpl}$" functions are thought of as
{\em replaced} by their values in terms of the five fundamental
functions $A$, $B$, $E$, $F$, $G$:
\[
\boxed{\,\,
\aligned
\big[\mathcal{R},\mathcal{S}\big]
&
=
-\,R_{\sf rpl}\cdot\mathcal{T}
-
S_{\sf rpl}\cdot\mathcal{S}
-
T_{\sf rpl}\cdot\mathcal{R},
\\
\big[\mathcal{R},\mathcal{T}\big]
&
=
-\,O_{\sf rpl}\cdot\mathcal{T}
-
P_{\sf rpl}\cdot\mathcal{S}
-
Q_{\sf rpl}\cdot\mathcal{R},
\\
\big[\mathcal{R},\overline{\mathcal{L}}\big]
&
=
-\,H_{\sf rpl}\cdot\mathcal{T}
-
J_{\sf rpl}\cdot\mathcal{S}
-
K_{\sf rpl}\cdot\mathcal{R},
\\
\big[\mathcal{R},\mathcal{L}\big]
&
=
-\,E\cdot\mathcal{T}
-
F\cdot\mathcal{S}
-
G\cdot\mathcal{R},
\\
\big[\mathcal{S},\mathcal{T}\big]
&
=
-\,L_{\sf rpl}\cdot\mathcal{T}
-
M_{\sf rpl}\cdot\mathcal{S}
-
N_{\sf rpl}\cdot\mathcal{R},
\\
\big[\mathcal{S},\overline{\mathcal{L}}\big]
&
=
-\,\mathcal{L}(A)\cdot\mathcal{T}
-
\big(\mathcal{L}(B)+A\big)\cdot\mathcal{S}
-
B\cdot\mathcal{R},\,\,
\\
\big[\mathcal{S},\mathcal{L}\big]
&
=
-\,\mathcal{R},
\\
\big[\mathcal{T},\overline{\mathcal{L}}\big]
&
=
-\,A\cdot\mathcal{T}
-
B\cdot\mathcal{S},
\\
\big[\mathcal{T},\mathcal{L}\big]
&
=
-\,\mathcal{S},
\\
\big[\overline{\mathcal{L}},\mathcal{L}\big]
&
=
\isqrt\,\mathcal{T}.
\endaligned}
\]

\medskip\noindent{\bf Initial Darboux structure of the
dual coframe.}
Dualize and introduce a dual coframe:
\[
\big\{
\mathcal{R},\,
\mathcal{S},\,\mathcal{T},\,\overline{\mathcal{L}},\,\mathcal{L}\big\}
\,\overset{\text{\sf dual}}{\longleftrightarrow}\,
\big\{
\tau_0,\,\sigma_0,\,\rho_0,\,\overline{\zeta}_0,\,\zeta_0
\big\}.
\]

\medskip

The convenient auxiliary array that one should read vertically then is:
\[
\footnotesize
\begin{array}{cccccccccccc}
& & \mathcal{R} & & \mathcal{S} & & \mathcal{T} & &
\overline{\mathcal{L}} & & \mathcal{L}
\\
& & \boxed{d\tau_0} & & \boxed{d\sigma_0} & & \boxed{d\rho_0} & &
\boxed{d\overline{\zeta_0}} & & \boxed{d\zeta_0}
\\
\big[\mathcal{R},\,\mathcal{S}\big]
& = & -\,T_{\sf rpl} & + & -\,S_{\sf rpl} & + & -\,R_{\sf rpl} & + & 
0 & + & 0 &
\boxed{\tau_0\wedge\sigma_0}
\\
\big[\mathcal{R},\,\mathcal{T}\big]
& = & -\,Q_{\sf rpl} & + & -\,P_{\sf rpl} & + & -\,O_{\sf rpl} & + & 
0 & + & 0 &
\boxed{\tau_0\wedge\rho_0}
\\
\big[\mathcal{R},\,\overline{\mathcal{L}}\big]
& = & -\,K_{\sf rpl} & + & -\,J_{\sf rpl} & + & -\,H_{\sf rpl}
& + & 0 & + & 0 &
\boxed{\tau_0\wedge\overline{\zeta_0}}
\\
\big[\mathcal{R},\,\mathcal{L}\big]
& = & -\,G & + & -\,F & + & -\,E & + & 0 & + & 0 &
\boxed{\tau_0\wedge\zeta_0}
\\
\big[\mathcal{S},\,\mathcal{T}\big]
& = & -\,N_{\sf rpl} & + & -\,M_{\sf rpl} & + & -\,L_{\sf rpl} & + & 
0 & + & 0 &
\boxed{\sigma_0\wedge\rho_0}
\\
\big[\mathcal{S},\,\overline{\mathcal{L}}\big]
& = & -\,B & + & -\,\mathcal{L}(B)-A & + & -\,\mathcal{L}(A) & + & 0 & + & 0 &
\boxed{\sigma_0\wedge\overline{\zeta}_0}
\\
\big[\mathcal{S},\,\mathcal{L}\big]
& = & -\,1 & + & 0 & + & 0 & + & 0 & + & 0 &
\boxed{\sigma_0\wedge\zeta_0}
\\
\big[\mathcal{T},\,\overline{\mathcal{L}}\big]
& = & 0 & + & -\,B & + & -\,A & + & 0 & + & 0 &
\boxed{\rho_0\wedge\overline{\zeta_0}}
\\
\big[\mathcal{T},\,\mathcal{L}\big]
& = & 0 & + & -\,1 & + & 0 & + & 0 & + & 0 &
\boxed{\rho_0\wedge\zeta_0}
\\
\big[\overline{\mathcal{L}},\,\mathcal{L}\big]
& = & 0 & + & 0 & + & i & + & 0 & + & 0 &
\boxed{\overline{\zeta}_0\wedge\zeta_0}
\end{array}
\]
and one obtains the initial Darboux structure for 
class $\text{\sf III}_{\text{\sf 2}}$ CR manifolds:
\[
\aligned
d\tau_0
&
=
T_{\sf rpl}\,\tau_0\wedge\sigma_0
+
Q_{\sf rpl}\,\tau_0\wedge\rho_0
+
K_{\sf rpl}\,\tau_0\wedge\overline{\zeta}_0
+
G\,\tau_0\wedge\zeta_0
+
\\
&
\ \ \ \ \
+
N_{\sf rpl}\,\sigma_0\wedge\rho_0
+
B\,\sigma_0\wedge\overline{\zeta}_0
+
\sigma_0\wedge\zeta_0,
\\
d\sigma_0
&
=
S_{\sf rpl}\,\tau_0\wedge\sigma_0
+
P_{\sf rpl}\,\tau_0\wedge\rho_0
+
J_{\sf rpl}\,\tau_0\wedge\overline{\zeta}_0
+
F\,\tau_0\wedge\zeta_0
+
M_{\sf rpl}\,\sigma_0\wedge\rho_0
+
\\
&
\ \ \ \ \
+
\big(\mathcal{L}(B)+A\big)\sigma_0\wedge\overline{\zeta}_0
+
B\,\rho_0\wedge\overline{\zeta}_0
+
\rho_0\wedge\zeta_0,
\\
d\rho_0
&
=
R_{\sf rpl}\,\tau_0\wedge\sigma_0
+
O_{\sf rpl}\,\tau_0\wedge\rho_0
+
H_{\sf rpl}\,\tau_0\wedge\overline{\zeta}_0
+
E\,\tau_0\wedge\zeta_0
+
L_{\sf rpl}\,\sigma_0\wedge\rho_0
+
\\
&
\ \ \ \ \ 
+
\mathcal{L}(A)\,\sigma_0\wedge\overline{\zeta}_0
+
A\,\rho_0\wedge\overline{\zeta}_0
+
\isqrt\,\zeta_0\wedge\overline{\zeta}_0,
\\
d\overline{\zeta}_0
&
=
0,
\\
d\zeta_0
&
=
0.
\endaligned
\]


\bigskip

\section{\sf $M^5 \subset \C^3$ of general class 
$\text{\sf IV}_{\text{\sf 1}}$: 
\\
initial frame and coframe in local coordinates}
\label{initial-IV-1}
\HEAD{\ref{initial-IV-1}.~$M^5 \subset \C^3$ of general class 
$\text{\sf IV}_{\text{\sf 1}}$
initial frame and coframe in local coordinates}{
Jo\"el {\sc Merker}, D\'epartement de Math\'ematiques d'Orsay}

\medskip

Consider:
\[
\Big(
M^5
\,\subset\,
\C^3
\Big)
\,\,\in\,\,
\text{\sf General Class $\text{\sf IV}_{\text{\sf 2}}$}.
\]
Represent $M$ as a (local) graph:
\[
v
=
\varphi\big(x_1,x_2,y_1,y_2,u\big),
\]
with:
\[
\varphi(0)
=
0,
\]
without (necessarily) requiring that $T_0 M = \{ v = 0\}$.

Provided only that:
\[
0
\neq
\isqrt+\varphi_u(0),
\]
two local generators for $T^{1, 0} M$ are:
\[
\aligned
\mathcal{L}_1
&
=
\frac{\partial}{\partial z_1}
-
\frac{\varphi_{z_1}}{\isqrt+\varphi_u}\,
\frac{\partial}{\partial u},
\\
\mathcal{L}_2
&
=
\frac{\partial}{\partial z_2}
-
\frac{\varphi_{z_2}}{\isqrt+\varphi_u}\,
\frac{\partial}{\partial u},
\endaligned
\]
having conjugates:
\[
\aligned
\overline{\mathcal{L}}_1
&
=
\frac{\partial}{\partial\overline{z}_1}
-
\frac{\varphi_{\overline{z}_1}}{-\isqrt+\varphi_u}\,
\frac{\partial}{\partial u},
\\
\overline{\mathcal{L}}_2
&
=
\frac{\partial}{\partial\overline{z}_2}
-
\frac{\varphi_{\overline{z}_2}}{-\isqrt+\varphi_u}\,
\frac{\partial}{\partial u}.
\endaligned
\]

Set as before in~\cite{ Merker-Pocchiola-Sabzevari-5-CR-II}:
\[
\aligned
A_1
&
:=
\frac{-\,\varphi_{z_1}}{\isqrt+\varphi_u},
\\
A_2
&
:=
\frac{-\,\varphi_{z_2}}{\isqrt+\varphi_u},
\endaligned
\]
two functions which are hence smooth near $0$.

Choose:
\[
\rho_0
:=
du
-
A_1\,dz_1
-
A_2\,dz_2
-
\overline{A}_1\,d\overline{z}_1
-
\overline{A}_2\,d\overline{z}_2,
\]
so that, in $\C \otimes_\R TM$:
\[
T^{1,0}M
\oplus
T^{0,1}M
\,=\,
\big\{
\rho_0=0
\big\}.
\]

Compute:
\[
\aligned
\mathcal{T}
:=
&\,
\isqrt\,\big[\mathcal{L}_1,\overline{\mathcal{L}}_1\big]
\\
=
&\,
\isqrt\,
\bigg[
\frac{\partial}{\partial z}_1
+
A_1\,\frac{\partial}{\partial u},
\,\,
\frac{\partial}{\partial\overline{z}}_1
+
\overline{A}_1\,\frac{\partial}{\partial u},
\bigg]
\\
&
=
\isqrt\,
\Big(
\mathcal{L}_1\big(\overline{A}_1\big)
-
\overline{\mathcal{L}}_1\big(A_1\big)
\Big)\,
\frac{\partial}{\partial u}.
\endaligned
\]

After a possible ${\sf GL}_2 ( \C)$-change of coordinates:
\[
(z_1,z_2)
\,\longmapsto\,
\big(\alpha\,z_1+\beta\,z_2,\,\,
\gamma\,z_1+\delta\,z_2
\big),
\]
one may assume that at the origin:
\[
0
\,\neq\,
\isqrt\,\Big(
\mathcal{L}_1\big(\overline{A}_1\big)
-
\overline{\mathcal{L}}_1\big(A_1\big)
\Big)
(0),
\]
hence in a neighborhood too.

Set:
\[
\ell_1
:=
\isqrt\,\Big(
\mathcal{L}_1\big(\overline{A}_1\big)
-
\overline{\mathcal{L}}_1\big(A_1\big)
\Big),
\]
so that:
\[
\overline{\ell}_{11}
=
\ell_{11}.
\]

Dually to the frame:
\[
\big\{
\mathcal{T},\,
\overline{\mathcal{L}}_1,\,
\overline{\mathcal{L}}_2,\,
\mathcal{L}_1,\,
\mathcal{L}_2
\big\},
\]
one has the coframe:
\[
\big\{
\rho_0,\,
\overline{\zeta}_{01},\,
\overline{\zeta}_{02},\,
\zeta_{01},\,
\zeta_{02}
\big\},
\]
where:
\[
\aligned
\rho_0
&
:=
\frac{du-\overline{A}_1\,d\overline{z}_1-\overline{A}_2\,d\overline{z}_2
-A_1\,dz_1-A_2\,dz_2}{\ell_{11}},
\\
d\overline{\zeta}_{01}
&
:=
d\overline{z}_1,
\\
d\overline{\zeta}_{02}
&
:=
d\overline{z}_2,
\\
d\zeta_{01}
&
:=
dz_1,
\\
d\zeta_{02}
&
:=
dz_2.
\endaligned
\]

Next, determine the Lie structure, namely:
\[
\aligned
&
\big[\mathcal{T},\overline{\mathcal{L}}_1\big],\ \ \ \ \
\big[\mathcal{T},\overline{\mathcal{L}}_2\big],\ \ \ \ \
\big[\mathcal{T},\mathcal{L}_1\big],\ \ \ \ \
\big[\mathcal{T},\mathcal{L}_2\big],
\\
&
\big[\overline{\mathcal{L}}_1,\overline{\mathcal{L}}_2\big],\ \ \ \ \
\big[\overline{\mathcal{L}}_1,\mathcal{L}_1\big],\ \ \ \ \
\big[\overline{\mathcal{L}}_1,\mathcal{L}_2\big],
\\
&
\big[\overline{\mathcal{L}}_2,\mathcal{L}_1\big],\ \ \ \ \
\big[\overline{\mathcal{L}}_2,\mathcal{L}_2\big],
\\
&
\big[\mathcal{L}_1,\mathcal{L}_2\big].
\endaligned
\]

Compute:
\[
\aligned
\big[\mathcal{T},\overline{\mathcal{L}}_1\big]
&
=
\bigg[
\ell_{11}\,\frac{\partial}{\partial u},\,\,
\frac{\partial}{\partial\overline{z}_1}
+
\overline{A}_1\,\frac{\partial}{\partial u}
\bigg]
\\
&
=
-\,
\Big(
\ell_{11,\overline{z}_1}
+
\overline{A}_1\,\ell_{11,u}
-
\ell_{11}\,\overline{A}_{1,u}
\Big)\,
\frac{\partial}{\partial u}
\\
&
=
-\,
\bigg(
\frac{
\ell_{11,\overline{z}_1}
+
\overline{A}_1\,\ell_{11,u}
-
\ell_{11}\,\overline{A}_{1,u}
}{
\ell_{11}}
\bigg)\,
\mathcal{T}
\\
&
=:
-\,
\overline{P}_1
\cdot
\mathcal{T},
\endaligned
\]
with (remind reality of $\ell_{ 11}$):
\[
P_1
:=
\frac{\ell_{11,z_1}+A_1\,\ell_{11,u}-\ell_{11}\,A_{1,u}}{\ell_{11}}.
\]

Similarly:
\[
\aligned
\big[\mathcal{T},\overline{\mathcal{L}}_2\big]
&
=
\bigg[
\ell_{11}\,\frac{\partial}{\partial u},\,\,
\frac{\partial}{\partial\overline{z}_2}
+
\overline{A}_2\,\frac{\partial}{\partial u}
\bigg]
\\
&
=
-\,
\Big(
\ell_{11,\overline{z}_2}
+
\overline{A}_2\,\ell_{11,u}
-
\ell_{11}\,\overline{A}_{2,u}
\Big)\,
\frac{\partial}{\partial u}
\\
&
=
-\,
\bigg(
\frac{
\ell_{11,\overline{z}_2}
+
\overline{A}_2\,\ell_{11,u}
-
\ell_{11}\,\overline{A}_{2,u}
}{
\ell_{11}}
\bigg)\,
\mathcal{T}
\\
&
=:
-\,
\overline{P}_2
\cdot
\mathcal{T},
\endaligned
\]
with:
\[
P_2
:=
\frac{\ell_{11,z_2}+A_2\,\ell_{11,u}-\ell_{11}\,A_{2,u}}{\ell_{11}}.
\]

Conjugating:
\[
\aligned
\big[\mathcal{T},\mathcal{L}_1\big]
&
=
-\,P_1\cdot\mathcal{T},
\\
\big[\mathcal{T},\mathcal{L}_2\big]
&
=
-\,P_2\cdot\mathcal{T}.
\endaligned
\]

Next, since:
\[
\aligned
\big[T^{1,0}M,\,T^{1,0}M\big]
&
\,\subset\,
T^{1,0}M,
\\
\big[T^{0,1}M,\,T^{0,1}M\big]
&
\,\subset\,
T^{0,1}M,
\endaligned
\]
one has here ({\em cf.} the
Scholium on pp.~26--28 in~\cite{Merker-Pocchiola-Sabzevari-5-CR-II}):
\[
\aligned
\big[\overline{\mathcal{L}}_1,\overline{\mathcal{L}}_2\big]
&
=
0,
\\
\big[\mathcal{L}_1,\mathcal{L}_2\big]
&
=
0.
\endaligned
\]

Also:
\[
\aligned
\big[\overline{\mathcal{L}}_1,\mathcal{L}_2\big]
&
=
\Big(
\overline{\mathcal{L}}_1\big(A_2\big)
-
\mathcal{L}_2\big(\overline{A}_1\big)
\Big)\,
\frac{\partial}{\partial u},
\\
\big[\overline{\mathcal{L}}_2,\mathcal{L}_1\big]
&
=
\Big(
\overline{\mathcal{L}}_2\big(A_1\big)
-
\mathcal{L}_1\big(\overline{A}_2\big)
\Big)\,
\frac{\partial}{\partial u},
\\
\big[\overline{\mathcal{L}}_2,\mathcal{L}_2\big]
&
=
\Big(
\overline{\mathcal{L}}_2\big(A_2\big)
-
\mathcal{L}_2\big(\overline{A}_2\big)
\Big)\,
\frac{\partial}{\partial u}.
\endaligned
\]

Setting:
\[
\aligned
\ell_{11}
&
:=
\isqrt\,
\Big(
\mathcal{L}_1\big(\overline{A}_1\big)
-
\overline{\mathcal{L}}_1\big(A_1\big)
\Big),
\ \ \ \ \ \ \ \ \ \ \ \ \ \ \ \ \ \
\ell_{12}
:=
\isqrt\,
\Big(
\mathcal{L}_2\big(\overline{A}_1\big)
-
\overline{\mathcal{L}}_1\big(A_2\big)
\Big),
\\ 
\ell_{21}
&
:=
\isqrt\,
\Big(
\mathcal{L}_1\big(\overline{A}_2\big)
-
\overline{\mathcal{L}}_2\big(A_1\big)
\Big),
\ \ \ \ \ \ \ \ \ \ \ \ \ \ \ \ \ \
\ell_{22}
:=
\isqrt\,
\Big(
\mathcal{L}_2\big(\overline{A}_2\big)
-
\overline{\mathcal{L}}_2\big(A_2\big)
\Big),
\endaligned
\]
observing:
\[
\aligned
\overline{\ell}_{11}
&
=
\ell_{11},
\\
\overline{\ell}_{12}
&
=
\ell_{21},
\\
\overline{\ell}_{22}
&
=
\ell_{22},
\endaligned
\]
one therefore has:
\[
\aligned
\big[\overline{\mathcal{L}}_1,\mathcal{L}_1\big]
&
=
\isqrt\,
\ell_{11}\,
\frac{\partial}{\partial u}
\\
&
=
\isqrt\,\mathcal{T},
\endaligned
\]
\[
\aligned
\big[\overline{\mathcal{L}}_1,\mathcal{L}_2\big]
&
=
\isqrt\,
\ell_{12}\,
\frac{\partial}{\partial u}
\\
&
=
\isqrt\,
\frac{\ell_{12}}{\ell_{11}}\,
\mathcal{T},
\endaligned
\]
\[
\aligned
\big[\overline{\mathcal{L}}_2,\mathcal{L}_1\big]
&
=
\isqrt\,
\ell_{21}\,
\frac{\partial}{\partial u}
\\
&
=
\isqrt\,
\frac{\ell_{21}}{\ell_{11}}\,
\mathcal{T},
\endaligned
\]
\[
\aligned
\big[\overline{\mathcal{L}}_2,\mathcal{L}_2\big]
&
=
\isqrt\,
\ell_{22}\,
\frac{\partial}{\partial u}
\\
&
=
\isqrt\,
\frac{\ell_{22}}{\ell_{11}}\,
\mathcal{T}.
\endaligned
\]

Abbreviate:
\[
\aligned
B
&
:=
\frac{\ell_{21}}{\ell_{11}},
\\
A
&
:=
\frac{\ell_{22}}{\ell_{11}},
\endaligned
\]
so that:
\[
\aligned
\big[\overline{\mathcal{L}}_1,\mathcal{L}_1\big]
&
=
\isqrt\,
\cdot
\mathcal{T},
\ \ \ \  
\endaligned
\]
\[
\aligned
\big[\overline{\mathcal{L}}_1,\mathcal{L}_2\big]
&
=
\isqrt\,B\,
\cdot
\mathcal{T}.
\endaligned
\]
\[
\aligned
\big[\overline{\mathcal{L}}_2,\mathcal{L}_1\big]
&
=
\isqrt\,\overline{B}\,
\cdot
\mathcal{T}.
\endaligned
\]
\[
\aligned
\big[\overline{\mathcal{L}}_2,\mathcal{L}_2\big]
&
=
\isqrt\,A\,
\cdot
\mathcal{T}.
\endaligned
\]

\medskip\noindent{\bf Summary.}
The $10$ Lie brackets gather as:
\[
\aligned
\big[\mathcal{T},\overline{\mathcal{L}}_1\big]
&
=
-\,\overline{P}_1\cdot\mathcal{T},
\\
\big[\mathcal{T},\overline{\mathcal{L}}_2\big]
&
=
-\,\overline{P}_2\cdot\mathcal{T},
\\
\big[\mathcal{T},\mathcal{L}_1\big]
&
=
-\,P_1\cdot\mathcal{T},
\\
\big[\mathcal{T},\mathcal{L}_2\big]
&
=
-\,P_2\cdot\mathcal{T},
\\
\big[\overline{\mathcal{L}}_1,\overline{\mathcal{L}}_2\big]
&
=
0,
\\
\big[\overline{\mathcal{L}}_1,\mathcal{L}_1\big]
&
=
\isqrt
\cdot
\mathcal{T},
\\
\big[\overline{\mathcal{L}}_1,\mathcal{L}_2\big]
&
=
\isqrt\,B
\cdot
\mathcal{T},
\\
\big[\overline{\mathcal{L}}_2,\mathcal{L}_1\big]
&
=
\isqrt\,\overline{B}
\cdot
\mathcal{T},
\\
\big[\overline{\mathcal{L}}_2,\mathcal{L}_2\big]
&
=
\isqrt\,A
\cdot
\mathcal{T},
\\
\big[\mathcal{L}_1,\mathcal{L}_2\big]
&
=
0,
\endaligned
\]
and incorporate $4$ fundamental functions:
\[
A,\ \ \ \ \ \ \ \ \
B,\ \ \ \ \ \ \ \ \
P_1,\ \ \ \ \ \ \ \ \
P_2.
\]

Organize then the ten Lie brackets as a 
convenient auxiliary array:
\[
\footnotesize
\begin{array}{cccccccccccc}
& & \mathcal{T} & & \overline{\mathcal{L}}_1 & &
\overline{\mathcal{L}_2} & &
\mathcal{L}_1 & & \mathcal{L}_2
\\
& & \boxed{d\rho_0} & & \boxed{d\overline{\zeta}_{01}} & &
\boxed{d\overline{\zeta}_{02}} & & 
\boxed{d\zeta_{01}} & & \boxed{d\zeta_{02}}\medskip
\\
\big[\mathcal{T},\,\overline{\mathcal{L}}_1\big] & = &
-\,\overline{P}_1\cdot\mathcal{T} & + &
0 & + & 0 & + & 0 &
+ & 0 & \boxed{\rho_0\wedge\overline{\zeta}_{01}}
\\
\big[\mathcal{T},\,\overline{\mathcal{L}_2}\big] & = &
-\,\overline{P}_1\cdot\mathcal{T} & + &
0 & + &
0 & + & 0 & + & 0 &
\boxed{\rho_0\wedge\overline{\zeta}_{02}}
\\
\big[\mathcal{T},\,\mathcal{L}_1\big] & = &
-\,P_1\cdot\mathcal{T} & + &
0 & + & 0
& + & 0 & + & 0 & \boxed{\rho_0\wedge\zeta_{01}}
\\
\big[\mathcal{T},\,\mathcal{L}_2\big] & = &
-\,P_2\cdot\mathcal{T} & + &
0 & + & 0 & + & 0 & + & 0 &
\boxed{\rho_0\wedge\zeta_{02}}
\\
\big[\overline{\mathcal{L}}_1,\,\overline{\mathcal{L}_2}\big] & = &
0 & + & 0 & +
& 0 & + & 0 & + & 0 &
\boxed{\overline{\zeta}_{01}\wedge\overline{\zeta}_{02}}
\\
\big[\overline{\mathcal{L}}_1,\,\mathcal{L}_1\big] & = &
\isqrt\cdot\mathcal{T} & + & 0 & + & 0 & + & 0 & + & 0 &
\boxed{\overline{\zeta}_{01}\wedge\zeta_{01}}
\\
\big[\overline{\mathcal{L}}_1,\,\mathcal{L}_2\big] & = &
\isqrt\,B\cdot\mathcal{T} & + & 0 & + &
0 & + & 0 & + & 0 & 
\boxed{\overline{\zeta}_{01}\wedge\zeta_{02}}
\\
\big[\overline{\mathcal{L}_2},\,\mathcal{L}_1\big] & = &
\isqrt\,\overline{B}\cdot\mathcal{T} & + & 0 & + & 0 & + & 0 & + & 0 &
\boxed{\overline{\zeta}_{02}\wedge\zeta_{01}}
\\
\big[\overline{\mathcal{L}_2},\,\mathcal{L}_2\big] & = & 
\isqrt\,A\cdot\mathcal{T}  & + & 
& +
& 0 & + & 0 & + & 0 & 
\boxed{\overline{\zeta}_{02}\wedge\zeta_{02}}
\\
\big[\mathcal{L}_1,\,\mathcal{L}_2\big] & = & 0 & + & 0 & + &
0 & + & 0 & + & 0 &
\boxed{\zeta_{01}\wedge\zeta_{02}}
\end{array}
\]

Hence the Darboux structure is:
\[
\aligned
d\rho_0
&
=
\overline{P}_1
\cdot
\rho_0\wedge\overline{\zeta}_{01}
+
\overline{P}_2
\cdot
\rho_0\wedge\overline{\zeta}_{02}
+
P_1
\cdot
\rho_0\wedge\zeta_{01}
+
P_2
\cdot
\rho_0\wedge\zeta_{02}
+
\\
&
\ \ \ \ \
+
\isqrt
\cdot
\zeta_{01}\wedge\overline{\zeta}_{01}
+
\isqrt\,B
\cdot
\zeta_{02}\wedge\overline{\zeta}_{01}
+
\isqrt\,\overline{B}
\cdot
\zeta_{01}\wedge\overline{\zeta}_{02}
+
\isqrt\,A
\cdot
\zeta_{02}\wedge\overline{\zeta}_{02},
\\
d\overline{\zeta}_{01}
&
=
0,
\\
d\overline{\zeta}_{02}
&
=
0,
\\
d\zeta_{01}
&
=
0,
\\
d\zeta_{02}
&
=
0.
\endaligned
\]


\bigskip

\section{\sf $M^5 \subset \C^3$ of general class 
$\text{\sf IV}_{\text{\sf 2}}$: 
\\
initial frame and coframe in local coordinates}
\label{initial-IV-2}
\HEAD{\ref{initial-IV-2}.~$M^5 \subset \C^3$ of general class 
$\text{\sf IV}_{\text{\sf 2}}$
initial frame and coframe in local coordinates}{
Jo\"el {\sc Merker}, D\'epartement de Math\'ematiques d'Orsay}

\medskip

Consider:
\[
\Big(
M^5
\,\subset\,
\C^3
\Big)
\,\,\in\,\,
\text{\sf General Class $\text{\sf IV}_{\text{\sf 2}}$}.
\]
Represent $M$ as a (local) graph:
\[
v
=
\varphi\big(x_1,x_2,y_1,y_2,u\big),
\]
with:
\[
\varphi(0)
=
0,
\]
without (necessarily) requiring that $T_0 M = \{ v = 0\}$.

Provided only that:
\[
0
\neq
\isqrt+\varphi_u(0),
\]
two local generators for $T^{1, 0} M$ are:
\[
\aligned
\mathcal{L}_1
&
=
\frac{\partial}{\partial z_1}
-
\frac{\varphi_{z_1}}{\isqrt+\varphi_u}\,
\frac{\partial}{\partial u},
\\
\mathcal{L}_2
&
=
\frac{\partial}{\partial z_2}
-
\frac{\varphi_{z_2}}{\isqrt+\varphi_u}\,
\frac{\partial}{\partial u},
\endaligned
\]
having conjugates:
\[
\aligned
\overline{\mathcal{L}}_1
&
=
\frac{\partial}{\partial\overline{z}_1}
-
\frac{\varphi_{\overline{z}_1}}{-\isqrt+\varphi_u}\,
\frac{\partial}{\partial u},
\\
\overline{\mathcal{L}}_2
&
=
\frac{\partial}{\partial\overline{z}_2}
-
\frac{\varphi_{\overline{z}_2}}{-\isqrt+\varphi_u}\,
\frac{\partial}{\partial u}.
\endaligned
\]

Set as before in~\cite{ Merker-Pocchiola-Sabzevari-5-CR-II}:
\[
\aligned
A_1
&
:=
\frac{-\,\varphi_{z_1}}{\isqrt+\varphi_u},
\\
A_2
&
:=
\frac{-\,\varphi_{z_2}}{\isqrt+\varphi_u},
\endaligned
\]
two functions which are hence smooth near $0$.

Choose:
\[
\rho_0
:=
du
-
A_1\,dz_1
-
A_2\,dz_2
-
\overline{A}_1\,d\overline{z}_1
-
\overline{A}_2\,d\overline{z}_2,
\]
so that, in $\C \otimes_\R TM$:
\[
T^{1,0}M
\oplus
T^{0,1}M
\,=\,
\big\{
\rho_0=0
\big\}.
\]

Hence (\cite{ Merker-Pocchiola-Sabzevari-5-CR-II}):
\[
\text{\sf Levi-Matrix}_{\mathcal{L},\overline{\mathcal{L}}}^M
\,=\,
\left(\!
\begin{array}{cc}
\isqrt\big(\mathcal{L}_1\big(\overline{A}_1\big)
-
\overline{\mathcal{L}}_1\big(A_1\big)\big)
&
\isqrt\big(\mathcal{L}_2\big(\overline{A}_1\big)
-
\overline{\mathcal{L}}_1\big(A_1\big)\big)
\\
\isqrt\big(\mathcal{L}_1\big(\overline{A}_2\big)
-
\overline{\mathcal{L}}_2\big(A_1\big)\big)
&
\isqrt\big(\mathcal{L}_2\big(\overline{A}_2\big)
-
\overline{\mathcal{L}}_2\big(A_1\big)\big)
\end{array}
\!\right).
\]

By hypothesis, this matrix is everywhere of rank $1$.

After a possible ${\sf GL}_2 ( \C)$-change of coordinates:
\[
(z_1,z_2)
\,\longmapsto\,
\big(\alpha\,z_1+\beta\,z_2,\,\,
\gamma\,z_1+\delta\,z_2
\big),
\]
one may assume that at the origin:
\[
0
\,\neq\,
\isqrt\big(
\mathcal{L}_1\big(\overline{A}_1\big)
-
\overline{\mathcal{L}}_1\big(A_1\big)
\big)
(0),
\]
hence in a neighborhood too.

So:
\[
{\bf 1}
\,=\,
\rank_\C
\Big(
\text{\sf Levi-Matrix}_{\mathcal{L},\overline{\mathcal{L}}}^M
\Big),
\]
which now means:
\[
0
\,\equiv\,
\left\vert\!
\begin{array}{cc}
\mathcal{L}_1\big(\overline{A}_1\big)
-
\overline{\mathcal{L}}_1\big(A_1\big)
&
\mathcal{L}_2\big(\overline{A}_1\big)
-
\overline{\mathcal{L}}_1\big(A_2\big)
\\
\mathcal{L}_1\big(\overline{A}_2\big)
-
\overline{\mathcal{L}}_2\big(A_1\big)
&
\mathcal{L}_2\big(\overline{A}_2\big)
-
\overline{\mathcal{L}}_2\big(A_2\big)
\end{array}
\!\right\vert.
\]

Introduce the very fundamental function which is
the negative of the quotient of the two entries of the first line:
\[
\boxed{\,
k
\,:=\,
-\,
\frac{
\mathcal{L}_2\big(\overline{A}_1\big)
-
\overline{\mathcal{L}}_1\big(A_2\big)}{
\mathcal{L}_1\big(\overline{A}_1\big)
-
\overline{\mathcal{L}}_1\big(A_1\big)}.\,}
\]

A (mental) exercise convinces that the kernel
of the Levi matrix is generated by the vector-valued
function:
\[
\left(\!\!
\begin{array}{c}
k
\\
1
\end{array}
\!\!\right)
\]
namely:
\[
\left(\!\!
\begin{array}{c}
0
\\
0
\end{array}
\!\!\right)
\,=\,
\left(\!\!
\begin{array}{cc}
\mathcal{L}_1\big(\overline{A}_1\big)
-
\overline{\mathcal{L}}_1\big(A_1\big)
&
\mathcal{L}_2\big(\overline{A}_1\big)
-
\overline{\mathcal{L}}_1\big(A_2\big)
\\
\mathcal{L}_1\big(\overline{A}_2\big)
-
\overline{\mathcal{L}}_2\big(A_1\big)
&
\mathcal{L}_2\big(\overline{A}_2\big)
-
\overline{\mathcal{L}}_2\big(A_2\big)
\end{array}
\!\!\right)
\left(\!\!
\begin{array}{c}
k
\\
1
\end{array}
\!\!\right).
\]

Introduce the $(1, 0)$ vector field:
\[
\boxed{\,
\mathcal{K}
\,:=\,
k\,\mathcal{L}_1+\mathcal{L}_2.\,}
\]

By~\cite{ Merker-Pocchiola-Sabzevari-5-CR-II}, 
$\mathcal{K}$ is a generator for the Levi-kernel
subbundle:
\[
K^{1,0}M
\,\subset\,
T^{1,0}M,
\]
and it is invariant through biholomorphisms, with conjugate:
\[
K^{0,1}M
\,\subset\,
T^{0,1}M.
\]

Furthermore (\cite{ Merker-Pocchiola-Sabzevari-5-CR-II}, pp.~72--73):
\[
\aligned
\big[K^{1,0}M,\,K^{1,0}M\big]
&
\,\subset\,
K^{1,0}M,
\\
\big[K^{0,1}M,\,K^{0,1}M\big]
&
\,\subset\,
K^{0,1}M,
\\
\big[K^{1,0}M,\,K^{0,1}M\big]
&
\,\subset\,
K^{1,0}M
\oplus
K^{0,1}M.
\endaligned
\]

Now, since:
\[
\aligned
\mathcal{K}
&
=
k\,\frac{\partial}{\partial z_1}
+
\frac{\partial}{\partial z_2}
+
\big(k\,A_1+A_2\big)\,
\frac{\partial}{\partial u},
\\
\overline{\mathcal{K}}
&
=
\overline{k}\,\frac{\partial}{\partial\overline{z}_1}
+
\frac{\partial}{\partial\overline{z}_2}
+
\big(\overline{k}\,\overline{A}_1+\overline{A}_2\big)\,
\frac{\partial}{\partial u},
\endaligned
\]
it is visible that the Lie bracket:
\[
\big[
\mathcal{K},\,\overline{\mathcal{K}}
\big]
\,=\,
\function\cdot\mathcal{K}
+
\function\cdot\overline{\mathcal{K}}
\]
contains {\em no} $\frac{ \partial}{ \partial z_2}$ and
{\em no} $\frac{ \partial}{ \partial \overline{z }_2}$, 
whence necessarily:
\[
\big[
\mathcal{K},\,\overline{\mathcal{K}}
\big]
=
0,
\]
which gives:
\[
\aligned
0
&
=
\bigg[
k\,\frac{\partial}{\partial z_1}
+
\frac{\partial}{\partial z_2}
+
\big(k\,A_1+A_2\big)\,
\frac{\partial}{\partial u},
\,\,\,
\overline{k}\,\frac{\partial}{\partial\overline{z}_1}
+
\frac{\partial}{\partial\overline{z}_2}
+
\big(\overline{k}\,\overline{A}_1
+
\overline{A}_2\big)\,
\frac{\partial}{\partial u}
\bigg]
\\
&
=
\mathcal{K}\big(\overline{k}\big)\,
\frac{\partial}{\partial\overline{z}_1}
-
\overline{\mathcal{K}}(k)\,
\frac{\partial}{\partial z_1}
+
\something\cdot
\frac{\partial}{\partial u},
\endaligned
\]
and yields:
\[
\aligned
0
&
\equiv
\mathcal{K}\big(\overline{k}\big),
\\
0
&
\equiv
\overline{\mathcal{K}}(k).
\endaligned
\]

Now, the natural choice of a frame for $\C \otimes_\R TM$ 
includes $\mathcal{ K}$:
\[
\big\{
\mathcal{T},\,
\overline{\mathcal{L}}_1,\,
\overline{\mathcal{K}},\,
\mathcal{L}_1,\,
\mathcal{K}
\big\},
\]
where:
\[
\mathcal{T}
:=
\isqrt\,
\big[\mathcal{L}_1,\overline{\mathcal{L}}_1\big]
\]
is, as always, real.

To determine the Lie structure of this frame,
$10$ bracket are extant:
\[
\aligned
\big[\mathcal{K},\mathcal{L}_1\big],\ \ \ \ \
\big[\mathcal{K},\overline{\mathcal{K}}\big],\ \ \ \ \
\big[\mathcal{K},\overline{\mathcal{L}}_1\big],\ \ \ \ \
\big[\mathcal{K},\mathcal{T}\big],
\\
\big[\mathcal{L}_1,\overline{\mathcal{K}}\big],\ \ \ \ \
\big[\mathcal{L}_1,\overline{\mathcal{L}}_1\big],\ \ \ \ \
\big[\mathcal{L}_1,\mathcal{T}\big],
\\
\big[\overline{\mathcal{K}},\overline{\mathcal{L}}_1\big],\ \ \ \ \
\big[\overline{\mathcal{K}},\mathcal{T}\big],
\\
\big[\overline{\mathcal{L}}_1,\mathcal{T}\big].
\endaligned
\]

Firstly:
\[
\aligned
\big[\mathcal{K},\mathcal{L}_1\big]
&
=
\big[k\,\mathcal{L}_1+\mathcal{L}_2,\,\,
\mathcal{L}_1\big]
\\
&
-\,
\mathcal{L}_1(k)
\cdot
\mathcal{L}_1.
\endaligned
\]

Secondly, as seen:
\[
\big[
\mathcal{K},\overline{\mathcal{K}}
\big]
=
0.
\]

Thirdly:
\[
\aligned
\big[\mathcal{K},\overline{\mathcal{L}}_1\big]
&
=
\big[k\,\mathcal{L}_1+\mathcal{L}_2,\,\,
\overline{\mathcal{L}}_1\big]
\\
&
\ \ \ \ \
-\,
\overline{\mathcal{L}}_1(k)\cdot\mathcal{L}_1
+
\\
&
\ \ \ \ \
+
\zero{k\,\big[\mathcal{L}_1,\overline{\mathcal{L}}_1\big]
+
\big[\mathcal{L}_2,\overline{\mathcal{L}}_1\big]},
\endaligned
\]
the last two terms disappearing by the definition of $k$, for:
\[
\aligned
\big[\mathcal{L}_1,\overline{\mathcal{L}}_1\big]
&
=
\Big(
\mathcal{L}_1\big(\overline{A}_1\big)
-
\overline{\mathcal{L}}_1\big(A_1\big)
\Big)\,
\frac{\partial}{\partial u},
\\
\big[\mathcal{L}_2,\overline{\mathcal{L}}_1\big]
&
=
\Big(
\mathcal{L}_2\big(\overline{A}_1\big)
-
\overline{\mathcal{L}}_1\big(A_2\big)
\Big)\,
\frac{\partial}{\partial u}.
\endaligned
\] 

Fourthly, using Jacobi:
\[
\aligned
\big[\mathcal{K},\mathcal{T}\big]
&
=
\isqrt\,
\Big[
\mathcal{K},\,
\big[\mathcal{L}_1,\overline{\mathcal{L}}_1\big]
\Big]
\\
&
=
-\,\isqrt\,
\Big[
\overline{\mathcal{L}}_1,\,
\underbrace{
\big[\mathcal{K},\mathcal{L}_1\big]}_{
-\,\mathcal{L}_1(k)\cdot\mathcal{L}_1}
\Big]
-
\isqrt\,
\Big[
\mathcal{L}_1,\,
\underbrace{
\big[\overline{\mathcal{L}}_1,\mathcal{K}\big]}_{
\overline{\mathcal{L}}_1(k)\cdot\mathcal{L}_1}
\Big]
\\
&
=
\isqrt\,
\big[
\overline{\mathcal{L}}_1,\,
\mathcal{L}_1(k)\,\mathcal{L}_1\big]
-
\isqrt\,
\big[\mathcal{L}_1,\,
\overline{\mathcal{L}}_1(k)\,\mathcal{L}_1\big]
\\
&
=
\isqrt\,
\Big(
\overline{\mathcal{L}}_1\big(\mathcal{L}_1(k)\big)
-
\mathcal{L}_1\big(\overline{\mathcal{L}}_1(k)\big)
\Big)
\cdot
\mathcal{L}_1
+
\\
&
\ \ \ \ \ 
+
\mathcal{L}_1(k)\,\isqrt\,
\big[\overline{\mathcal{L}}_1,\mathcal{L}_1\big]
-
\isqrt\,\overline{\mathcal{L}}_1(k)\,
\zero{\big[\mathcal{L}_1,\mathcal{L}_1\big]}
\\
&
=
-\,
\mathcal{T}(k)\cdot\mathcal{L}_1
-
\mathcal{L}_1(k)
\cdot
\mathcal{T},
\endaligned
\]
that is to say:
\[
\big[\mathcal{K},\mathcal{T}\big]
=
-\,
\mathcal{T}(k)\cdot\mathcal{L}_1
-
\mathcal{L}_1(k)
\cdot
\mathcal{T}.
\]

Fifthly, conjugating an equation that precedes:
\[
\big[\mathcal{L}_1,\overline{\mathcal{K}}\big]
=
\mathcal{L}_1\big(\overline{k}\big)
\cdot
\overline{\mathcal{L}}_1.
\]

Sixthly:
\[
\big[
\mathcal{L}_1,\overline{\mathcal{L}}_1
\big]
=
-\,\isqrt\,\mathcal{T}.
\]

Seventhly, abbreviating:
\[
\mathcal{T}
=
\ell\,\frac{\partial}{\partial u},
\]
with:
\[
\ell
:=
\isqrt\,\big(
\mathcal{L}_1\big(\overline{A}_1\big)
-
\overline{\mathcal{L}}_1\big(A_1\big)
\big),
\]
compute:
\[
\aligned
\big[\mathcal{L}_1,\mathcal{T}\big]
&
=
\bigg[
\frac{\partial}{\partial z_1}
+
A_1\,\frac{\partial}{\partial u},\,\,
\ell\,\frac{\partial}{\partial u}
\bigg]
\\
&
=
\Big(
\ell_{z_1}
+
A_1\,\ell_u
-
\ell\,A_{1,u}
\Big)\,
\frac{\partial}{\partial u}
\\
&
=
\frac{\ell_{z_1}+A_1\,\ell_u-\ell\,A_{1,u}}{\ell}\,
\mathcal{T}.
\endaligned
\]
One therefore comes to the second and last fundamental
function:
\[
\boxed{\,
P
:=
\frac{\ell_{z_1}+A_1\,\ell_u-\ell\,A_{1,u}}{\ell},\,}
\]
and gets:
\[
\big[\mathcal{L}_1,\mathcal{T}\big]
=
P\cdot\mathcal{T}.
\]

Eighthly:
\[
\big[
\overline{\mathcal{K}},\overline{\mathcal{L}}_1
\big]
=
-\,\overline{\mathcal{L}}_1\big(\overline{k}\big)
\cdot
\overline{\mathcal{L}}_1.
\]

Ninthly:
\[
\big[
\overline{\mathcal{K}},\mathcal{T}\big]
=
-\,\overline{\mathcal{T}}\big(\overline{k}\big)
\cdot
\overline{\mathcal{L}}_1
-
\overline{\mathcal{L}}_1\big(\overline{k}\big)
\cdot
\overline{\mathcal{T}}.
\]

Tenthly:
\[
\big[\overline{\mathcal{L}}_1,\mathcal{T}\big]
=
\overline{P}\cdot\mathcal{T}.
\]

\medskip\noindent{\bf Summary.} One has the 10 Lie bracket relations:
\[
\boxed{\,\,
\aligned
\big[\mathcal{T},\overline{\mathcal{L}}_1\big]
&
=
-\,\overline{P}\cdot\mathcal{T},
\\
\big[\mathcal{T},\overline{\mathcal{K}}\big]
&
=
\overline{\mathcal{L}}_1\big(\overline{k}\big)
\cdot
\mathcal{T}
+
\mathcal{T}
\big(\overline{k}\big)
\cdot
\overline{\mathcal{L}}_1,
\\
\big[\mathcal{T},\mathcal{L}_1\big]
&
=
-\,P\cdot\mathcal{T},
\\
\big[\mathcal{T},\overline{\mathcal{K}}\big]
&
=
\mathcal{L}_1(k)\cdot\mathcal{T}
+
\mathcal{T}(k)\cdot\mathcal{L}_1,
\\
\big[\overline{\mathcal{L}}_1,\overline{\mathcal{K}}\big]
&
=
\overline{\mathcal{L}}_1\big(\overline{k}\big)\cdot
\overline{\mathcal{L}}_1,
\\
\big[\overline{\mathcal{L}}_1,\mathcal{L}_1\big]
&
=
\isqrt\,\mathcal{T},
\\
\big[\overline{\mathcal{L}}_1,\mathcal{K}\big]
&
=
\overline{\mathcal{L}}_1(k)
\cdot\mathcal{L}_1,
\\
\big[\overline{\mathcal{K}},\mathcal{L}_1\big]
&
=
-\,\mathcal{L}_1\big(\overline{k}\big)
\cdot\overline{\mathcal{L}}_1,
\\
\big[\overline{\mathcal{K}},\mathcal{K}\big]
&
=
0,
\\
\big[\mathcal{L}_1,\mathcal{K}\big]
&
=
\mathcal{L}_1(k)\cdot\mathcal{L}_1.
\endaligned}
\]

\medskip\noindent{\bf Initial Darboux structure of the 
dual coframe.} 
Introduce then the coframe:
\[
\big\{
\rho_0,\,\overline{\kappa}_0,\,\overline{\zeta}_0,\,
\kappa_0,\,\zeta_0
\big\}
\]
which is dual to the frame:
\[
\big\{
\mathcal{T},\,
\overline{\mathcal{L}}_1,\,
\overline{\mathcal{K}},\,
\mathcal{L}_1,\,
\mathcal{K}
\big\},
\]
namely:
\[
\begin{array}{ccccc}
\rho_0(\mathcal{T})=1 \ \ \ & \ \ \
\rho_0(\overline{\mathcal{L}}_1)=0 \ \ \ & \ \ \
\rho_0(\overline{\mathcal{K}})=0 \ \ \ & \ \ \
\rho_0(\mathcal{L}_1)=0 \ \ \ & \ \ \
\rho_0(\mathcal{K})=0,
\\
\overline{\kappa}_0(\mathcal{T})=1 \ \ \ & \ \ \
\overline{\kappa}_0(\overline{\mathcal{L}}_1)=0 \ \ \ & \ \ \
\overline{\kappa}_0(\overline{\mathcal{K}})=0 \ \ \ & \ \ \
\overline{\kappa}_0(\mathcal{L}_1)=0 \ \ \ & \ \ \
\overline{\kappa}_0(\mathcal{K})=0,
\\
\overline{\zeta}_0(\mathcal{T})=1 \ \ \ & \ \ \
\overline{\zeta}_0(\overline{\mathcal{L}}_1)=0 \ \ \ & \ \ \
\overline{\zeta}_0(\overline{\mathcal{K}})=0 \ \ \ & \ \ \
\overline{\zeta}_0(\mathcal{L}_1)=0 \ \ \ & \ \ \
\overline{\zeta}_0(\mathcal{K})=0,
\\
\kappa_0(\mathcal{T})=1 \ \ \ & \ \ \
\kappa_0(\overline{\mathcal{L}}_1)=0 \ \ \ & \ \ \
\kappa_0(\overline{\mathcal{K}})=0 \ \ \ & \ \ \
\kappa_0(\mathcal{L}_1)=0 \ \ \ & \ \ \
\kappa_0(\mathcal{K})=0,
\\
\zeta_0(\mathcal{T})=1 \ \ \ & \ \ \
\zeta_0(\overline{\mathcal{L}}_1)=0 \ \ \ & \ \ \
\zeta_0(\overline{\mathcal{K}})=0 \ \ \ & \ \ \
\zeta_0(\mathcal{L}_1)=0 \ \ \ & \ \ \
\zeta_0(\mathcal{K})=0.
\end{array}
\]
One has:
\[
\aligned
\rho_0
&
=
\frac{du-A_1dz_1-A_2dz_2
-\overline{A}_1d\overline{z}_1
-\overline{A}_2d\overline{z}_2}{\ell},
\\
\kappa_0
&
=
dz_1-k\,dz_2,
\\
\zeta_0
&
=
dz_2.
\endaligned
\]

Organize the ten Lie brackets as a 
convenient auxiliary array:
\[
\footnotesize
\begin{array}{cccccccccccc}
& & \mathcal{T} & & \overline{\mathcal{L}}_1 & &
\overline{\mathcal{K}} & &
\mathcal{L}_1 & & \mathcal{K}
\\
& & \boxed{d\rho_0} & & \boxed{d\overline{\kappa_0}} & &
\boxed{d\overline{\zeta}_0} & & 
\boxed{d\kappa_0} & & \boxed{d\zeta_0}\medskip
\\
\big[\mathcal{T},\,\overline{\mathcal{L}}_1\big] & = &
-\,\overline{P}\cdot\mathcal{T} & + &
0 & + & 0 & + & 0 &
+ & 0 & \boxed{\rho_0\wedge\overline{\kappa}_0}
\\
\big[\mathcal{T},\,\overline{\mathcal{K}}\big] & = &
\overline{\mathcal{L}}_1\big(\overline{k}\big)\cdot\mathcal{T} & + &
\mathcal{T}\big(\overline{k}\big)\cdot\overline{\mathcal{L}}_1 & + &
0 & + & 0 & + & 0 &
\boxed{\rho_0\wedge\overline{\zeta}_0}
\\
\big[\mathcal{T},\,\mathcal{L}_1\big] & = &
-\,P\cdot\mathcal{T} & + &
0 & + & 0
& + & 0 & + & 0 & \boxed{\rho_0\wedge\kappa_0}
\\
\big[\mathcal{T},\,\mathcal{K}\big] & = &
\mathcal{L}_1(k)\cdot\mathcal{T} & + &
0 & + & 0 & + & \mathcal{T}(k)\cdot\mathcal{L}_1 & + & 0 &
\boxed{\rho_0\wedge\zeta_0}
\\
\big[\overline{\mathcal{L}}_1,\,\overline{\mathcal{K}}\big] & = &
0 & + & 
\overline{\mathcal{L}}_1\big(\overline{k}\big)
\cdot\overline{\mathcal{L}}_1 & +
& 0 & + & 0 & + & 0 &
\boxed{\overline{\kappa}_0\wedge\overline{\zeta}_0}
\\
\big[\overline{\mathcal{L}}_1,\,\mathcal{L}_1\big] & = &
\isqrt\,\cdot\mathcal{T} & + & 0 & + & 0 & + & 0 & + & 0 &
\boxed{\overline{\kappa}_0\wedge\kappa_0}
\\
\big[\overline{\mathcal{L}}_1,\,\mathcal{K}\big] & = &
0 & + & 0 & + &
0 & + & \overline{\mathcal{L}}_1(k)\cdot\mathcal{L}_1 & + & 0 & 
\boxed{\overline{\kappa}_0\wedge\zeta_0}
\\
\big[\overline{\mathcal{K}},\,\mathcal{L}_1\big] & = &
0 & + & -\,\mathcal{L}_1\big(\overline{k}\big)\cdot
\overline{\mathcal{L}}_1 & + & 0 & + & 0 & + & 0 &
\boxed{\overline{\zeta_0}\wedge\kappa_0}
\\
\big[\overline{\mathcal{K}},\,\mathcal{K}\big] & = & 0 & + & 
0 & +
& 0 & + & 0 & + & 0 & 
\boxed{\overline{\zeta}_0\wedge\zeta_0}
\\
\big[\mathcal{L}_1,\,\mathcal{K}\big] & = & 0 & + & 0 & + &
0 & + & \mathcal{L}_1(k)\cdot\mathcal{L}_1 & + & 0 &
\boxed{\kappa_0\wedge\zeta_0}
\end{array}
\]

Read {\em vertically} and put an overall minus sign:
\[
\boxed{ 
\aligned 
d\rho_0
&
=
\overline{P}
\cdot
\rho_0\wedge\overline{\kappa}_0
-
\overline{\mathcal{L}}_1\big(\overline{k}\big)
\cdot
\rho_0\wedge\overline{\zeta}_0
+
P
\cdot
\rho_0\wedge\kappa_0
-
\mathcal{L}_1(k)
\cdot
\rho_0\wedge\zeta_0
+
\isqrt\,
\kappa_0\wedge\overline{\kappa}_0,
\\
d\overline{\kappa}_0
&
=
-\,\mathcal{T}\big(\overline{k}\big)
\cdot
\rho_0\wedge\overline{\zeta}_0
-
\overline{\mathcal{L}}_1\big(\overline{k}\big)
\cdot
\overline{\kappa}_0
\wedge
\overline{\zeta}_0
+
\mathcal{L}_1\big(\overline{k}\big)
\cdot
\overline{\zeta}_0
\wedge
\kappa_0,
\\
d\overline{\zeta}_0
&
=
0,
\\
d\kappa_0
&
=
-\,\mathcal{T}(k)
\cdot
\rho_0\wedge\zeta_0
-
\overline{\mathcal{L}}_1(k)
\cdot
\overline{\kappa}_0\wedge\zeta_0
-
\mathcal{L}_1(k)
\cdot
\kappa_0\wedge\zeta_0,
\\
d\zeta_0
&
=
0.
\endaligned}
\]

\medskip

This together with the related matrix ambiguity group
(\cite{ Merker-5-CR-IV}):
\[
\aligned
{\sf G}_{\text{\sf IV}_{\text{\sf 2}}}^{\sf initial}
\,:=\,
\left\{
\left(\!\!
\begin{array}{ccccc}
{\sf c} & 0 & 0 & 0 & 0
\\
{\sf b} & {\sf a} & 0 & 0 & 0
\\
0 & 0 & \overline{\sf c} & 0 & 0
\\
0 & 0 & \overline{\sf b} & \overline{\sf a} & 0
\\
{\sf e} & {\sf d} & \overline{\sf e} & \overline{\sf d} & 
{\sf a}\overline{\sf a}
\end{array}
\!\!\right)
\,\in\,
\mathcal{M}_{5\times 5}(\C)\,\colon\,\,
{\sf a},\,{\sf c}\,\in\,\C\backslash\{0\},\,\,
{\sf b},\,{\sf d},\,{\sf e}
\,\in\,\C
\right\}
\endaligned
\]
initializes (launches) an explicit application of Cartan's equivalence 
method achieved by Pocchiola in~\cite{ Pocchiola-2013}.


\vfill\end{document}